\font\tenmib=cmmib10
\font\sevenmib=cmmib10 scaled 800
\font\titolo=cmbx12
\font\titolone=cmbx10 scaled\magstep 2

\font\cs=cmcsc10
\font\sc=cmcsc10
\font\css=cmcsc8

\font\ninerm=cmr9
\font\ottorm=cmr8
\textfont5=\tenmib\scriptfont5=\sevenmib\scriptscriptfont5=\fivei

\font\euftw=eufm9 scaled\magstep1
\font\euftww=eufm7 scaled\magstep1
\font\msytw=msbm9 scaled\magstep1

\font\msytwww=msbm7 scaled\magstep1
\font\msytwwww=msbm5 scaled\magstep1
\font\msytwwwww=msbm4 scaled\magstep1
\font\indbf=cmbx10 scaled\magstep2

\font\ottorm=cmr8\font\ottoi=cmmi8\font\ottosy=cmsy8%
\font\ottobf=cmbx8\font\ottott=cmtt8%
\font\ottocss=cmcsc8%
\font\ottosl=cmsl8\font\ottoit=cmti8%
\font\sixrm=cmr6\font\sixbf=cmbx6\font\sixi=cmmi6\font\sixsy=cmsy6%
\font\fiverm=cmr5\font\fivesy=cmsy5\font\fivei=cmmi5\font\fivebf=cmbx5%

\def\ottopunti{\def\rm{\fam0\ottorm}%
\textfont0=\ottorm\scriptfont0=\sixrm\scriptscriptfont0=\fiverm%
\textfont1=\ottoi\scriptfont1=\sixi\scriptscriptfont1=\fivei%
\textfont2=\ottosy\scriptfont2=\sixsy\scriptscriptfont2=\fivesy%
\textfont3=\tenex\scriptfont3=\tenex\scriptscriptfont3=\tenex%
\textfont4=\ottocss\scriptfont4=\sc\scriptscriptfont4=\sc%
\textfont\itfam=\ottoit\def\it{\fam\itfam\ottoit}%
\textfont\slfam=\ottosl\def\sl{\fam\slfam\ottosl}%
\textfont\ttfam=\ottott\def\tt{\fam\ttfam\ottott}%
\textfont\bffam=\ottobf\scriptfont\bffam=\sixbf%
\scriptscriptfont\bffam=\fivebf\def\bf{\fam\bffam\ottobf}%
\setbox\strutbox=\hbox{\vrule height7pt depth2pt width0pt}%
\normalbaselineskip=9pt\let\sc=\sixrm\normalbaselines\rm}
\let\nota=\ottopunti%

%
%
%
%
%
%
%

\global\newcount\numsec\global\newcount\numapp
\global\newcount\numfor\global\newcount\numfig
\global\newcount\numsub
\numsec=0\numapp=0\numfig=0
\def\veroparagrafo{\number\numsec}\def\veraformula{\number\numfor}
\def\veraappendice{\number\numapp}\def\verasub{\number\numsub}
\def\verafigura{\number\numfig}

\def\section(#1,#2){\advance\numsec by 1\numfor=1\numsub=1\numfig=1%
\SIA p,#1,{\veroparagrafo} %
\write15{\string\Fp (#1){\secc(#1)}}%
\write16{ sec. #1 ==> \secc(#1)  }%
\hbox to \hsize{\titolo\hfill \number\numsec. #2\hfill%
\expandafter{\alato(sec. #1)}}\*}

\def\appendix(#1,#2){\advance\numapp by 1\numfor=1\numsub=1\numfig=1%
\SIA p,#1,{A\veraappendice} %
\write15{\string\Fp (#1){\secc(#1)}}%
\write16{ app. #1 ==> \secc(#1)  }%
\hbox to \hsize{\titolo\hfill Appendix A\number\numapp. #2\hfill%
\expandafter{\alato(app. #1)}}\*}

\def\senondefinito#1{\expandafter\ifx\csname#1\endcsname\relax}

\def\SIA #1,#2,#3 {\senondefinito{#1#2}%
\expandafter\xdef\csname #1#2\endcsname{#3}\else
\write16{???? ma #1#2 e' gia' stato definito !!!!} \fi}

\def \Fe(#1)#2{\SIA fe,#1,#2 }
\def \Fp(#1)#2{\SIA fp,#1,#2 }
\def \Fg(#1)#2{\SIA fg,#1,#2 }

\def\etichetta(#1){(\veroparagrafo.\veraformula)%
\SIA e,#1,(\veroparagrafo.\veraformula) %
\global\advance\numfor by 1%
\write15{\string\Fe (#1){\equ(#1)}}%
\write16{ EQ #1 ==> \equ(#1)  }}

\def\etichettaa(#1){(A\veraappendice.\veraformula)%
\SIA e,#1,(A\veraappendice.\veraformula) %
\global\advance\numfor by 1%
\write15{\string\Fe (#1){\equ(#1)}}%
\write16{ EQ #1 ==> \equ(#1) }}

\def\getichetta(#1){\veroparagrafo.\verafigura%
\SIA g,#1,{\veroparagrafo.\verafigura} %
\global\advance\numfig by 1%
\write15{\string\Fg (#1){\graf(#1)}}%
\write16{ Fig. #1 ==> \graf(#1) }}

\def\getichettaa(#1){A\veraappendice.\verafigura%
\SIA g,#1,{A\veraappendice.\verafigura} %
\global\advance\numfig by 1%
\write15{\string\Fg (#1){\graf(#1)}}%
\write16{ Fig. #1 ==> \graf(#1) }}

\def\etichettap(#1){\veroparagrafo.\verasub%
\SIA p,#1,{\veroparagrafo.\verasub} %
\global\advance\numsub by 1%
\write15{\string\Fp (#1){\secc(#1)}}%
\write16{ par #1 ==> \secc(#1)  }}

\def\etichettapa(#1){A\veraappendice.\verasub%
\SIA p,#1,{A\veraappendice.\verasub} %
\global\advance\numsub by 1%
\write15{\string\Fp (#1){\secc(#1)}}%
\write16{ par #1 ==> \secc(#1)  }}

\def\Eq(#1){\eqno{\etichetta(#1)\alato(#1)}}
\def\eq(#1){\etichetta(#1)\alato(#1)}
\def\Eqa(#1){\eqno{\etichettaa(#1)\alato(#1)}}
\def\eqa(#1){\etichettaa(#1)\alato(#1)}
\def\eqg(#1){\getichetta(#1)\alato(fig. #1)}
\def\sub(#1){\0\palato(p. #1){\bf \etichettap(#1).}}
\def\asub(#1){\0\palato(p. #1){\bf \etichettapa(#1).}}

\def\equv(#1){\senondefinito{fe#1}$\clubsuit$#1%
\write16{eq. #1 non e' (ancora) definita}%
\else\csname fe#1\endcsname\fi}
\def\grafv(#1){\senondefinito{fg#1}$\clubsuit$#1%
\write16{fig. #1 non e' (ancora) definito}%
\else\csname fg#1\endcsname\fi}
\def\secv(#1){\senondefinito{fp#1}$\clubsuit$#1%
\write16{par. #1 non e' (ancora) definito}%
\else\csname fp#1\endcsname\fi}

\def\equ(#1){\senondefinito{e#1}\equv(#1)\else\csname e#1\endcsname\fi}
\def\graf(#1){\senondefinito{g#1}\grafv(#1)\else\csname g#1\endcsname\fi}
\def\figura(#1){{\css Figure} \getichetta(#1)}
\def\figuraa(#1){{\css Figure} \getichettaa(#1)}
\def\secc(#1){\senondefinito{p#1}\secv(#1)\else\csname p#1\endcsname\fi}
\def\sec(#1){{\S\secc(#1)}}
\def\refe(#1){{[\secc(#1)]}}

\def\BOZZA{\bz=1
\def\alato(##1){\rlap{\kern-\hsize\kern-1.2truecm{$\scriptstyle##1$}}}
\def\palato(##1){\rlap{\kern-1.2truecm{$\scriptstyle##1$}}}
}

\def\alato(#1){}
\def\galato(#1){}
\def\palato(#1){}


{\count255=\time\divide\count255 by 60 \xdef\hourmin{\number\count255}
        \multiply\count255 by-60\advance\count255 by\time
   \xdef\hourmin{\hourmin:\ifnum\count255<10 0\fi\the\count255}}

\def\oramin{\hourmin }

\def\data{\number\day/\ifcase\month\or gennaio \or febbraio \or marzo \or
aprile \or maggio \or giugno \or luglio \or agosto \or settembre
\or ottobre \or novembre \or dicembre \fi/\number\year;\ \oramin}
\footline={\rlap{\hbox{\copy200}}\tenrm\hss \number\pageno\hss}


\newcount\driver 
\newdimen\xshift \newdimen\xwidth
\def\ins#1#2#3{\vbox to0pt{\kern-#2 \hbox{\kern#1
#3}\vss}\nointerlineskip}

\def\insertplot#1#2#3#4#5{\par%
\xwidth=#1 \xshift=\hsize \advance\xshift by-\xwidth \divide\xshift by 2%
\yshift=#2 \divide\yshift by 2%
\line{\hskip\xshift \vbox to #2{\vfil%
\ifnum\driver=0 #3
\special{ps: plotfile #4.eps} 
\fi
\ifnum\driver=1 #3 \includegraphics{#4.eps}\fi
\ifnum\driver=2 #3
\ifnum\mgnf=0\special{#4.eps 1. 1. scale} \fi
\ifnum\mgnf=1\special{#4.eps 1.2 1.2 scale}\fi
\fi }\hfill \raise\yshift\hbox{#5}}}

\def\insertplotttt#1#2#3#4{\par%
\xwidth=#1 \xshift=\hsize \advance\xshift by-\xwidth \divide\xshift by 2%
\yshift=#2 \divide\yshift by 2%
\line{\hskip\xshift \vbox to #2{\vfil%
\ifnum\driver=0 #3
\special{ps: plotfile #4.eps} 
\fi
\ifnum\driver=1 #3 \includegraphics{#4.eps}\fi
\ifnum\driver=2 #3
\ifnum\mgnf=0\special{#4.eps 1. 1. scale} \fi
\ifnum\mgnf=1\special{#4.eps 1.2 1.2 scale}\fi
\fi }\hfill}}

\def\insertplott#1#2#3{\par%
\xwidth=#1 \xshift=\hsize \advance\xshift by-\xwidth \divide\xshift by 2%
\yshift=#2 \divide\yshift by 2%
\line{\hskip\xshift \vbox to #2{\vfil%
\ifnum\driver=0
\special{ps: plotfile #3.eps} 
 \fi
\ifnum\driver=1 \includegraphics{#3.eps}\fi
\ifnum\driver=2
\ifnum\mgnf=0\special{#3.eps 1. 1. scale} \fi 
\ifnum\mgnf=1\special{#3.eps 1.2 1.2 scale}\fi
\fi }\hfill}}

\newdimen\xshift \newdimen\xwidth \newdimen\yshift
\def\eqfig#1#2#3#4#5{
\par\xwidth=#1 \xshift=\hsize \advance\xshift
by-\xwidth \divide\xshift by 2
\yshift=#2 \divide\yshift by 2
\line{\hglue\xshift \vbox to #2{\vfil
\ifnum\driver=0 #3
\special{ps: plotfile #4.ps} 
\fi
\ifnum\driver=1 #3 \includegraphics{#4.ps}\fi
\ifnum\driver=2 #3 \special{
\ifnum\mgnf=0 #4.ps 1. 1. scale \fi
\ifnum\mgnf=1 #4.ps 1.2 1.2 scale\fi}
\fi}\hfill\raise\yshift\hbox{#5}}}

\let\a=\alpha \let\b=\beta  \let\g=\gamma  \let\d=\delta \let\e=\varepsilon
  \let\h=\eta   \let\th=\theta \let\k=\kappa \let\l=\lambda
\let\m=\mu    \let\n=\nu    \let\x=\xi     \let\p=\pi    
\let\s=\sigma \let\t=\tau   \let\f=\varphi 
   \let\o=\omega
\let\G=\Gamma \let\D=\Delta  \let\Th=\Theta 
     \let\F=\Phi    
\let\O=\Omega 

\def\\{\hfill\break} \let\==\equiv

\let\io=\infty 

\let\0=\noindent

\def\ie{\hbox{\it i.e.\ }}
\let\dpr=\partial

\def\mm{\hbox{\euftw m}}\def\mmmm{\hbox{\euftww m}}

\def\sn{{\rm sn}}\def\cn{{\rm cn}}\def\dn{{\rm dn}}

\def\tende#1{\,\vtop{\ialign{##\crcr\rightarrowfill\crcr
 \noalign{\kern-1pt\nointerlineskip}
 \hskip3.pt${\scriptstyle #1}$\hskip3.pt\crcr}}\,}
\def\otto{\,{\kern-1.truept\leftarrow\kern-5.truept\to\kern-1.truept}\,}

\def\EE{{\cal E}}\def\MM{{\cal M}} \def\VV{{\cal V}}

\def\TT{{\cal T}}\def\III{{\cal I}}
\def\RR{{\cal R}}\def\LL{{\cal L}}\def\JJ{{\cal J}} 
\def\SS{{\cal S}}

\def\T#1{{#1_{\kern-3pt\lower7pt\hbox{$\widetilde{}$}}\kern3pt}}
\def\VVV#1{{\underline #1}_{\kern-3pt
\lower7pt\hbox{$\widetilde{}$}}\kern3pt\,}
\def\W#1{#1_{\kern-3pt\lower7.5pt\hbox{$\widetilde{}$}}\kern2pt\,}

\def\lis{\overline}

\def\indica{\leaders \hbox to 0.5cm{\hss.\hss}\hfill}
\def\guida{\leaders\hbox to 1em{\hss.\hss}\hfill}

   \def\V0{{\bf 0}}
   \def\PP{{\bf P}} 
    
\def\II{{\bf I}} \def\LLL{{\bf L}}
\def\ul{\underline}
\def\wt{\widetilde}
\def\oo{{\omega}}

\mathchardef\aa   = "050B
\mathchardef\bb   = "050C
\mathchardef\xxx  = "0518
\mathchardef\hhh  = "0511
\mathchardef\zzzzz= "0510
\mathchardef\oo   = "0521
\mathchardef\tt   = "051C
\mathchardef\lll  = "0515
\mathchardef\mmm  = "0516
\mathchardef\Dp   = "0540
\mathchardef\H    = "0548
\mathchardef\FFF  = "0546
\mathchardef\ppp  = "0570
\mathchardef\nnn  = "0517
\mathchardef\pps  = "0520
\mathchardef\XXX  = "0504
\mathchardef\FFF  = "0508
\mathchardef\nnnnn= "056E
\mathchardef\minore  = "053C
\mathchardef\maggiore  = "053E

\def\to{\rightarrow}
\def\la{\left\langle}
\def\ra{\right\rangle}

\let\ciao=\bye
\def\qed{\hfill\raise1pt\hbox{\vrule height5pt width5pt depth0pt}}

 \def\Val{{\rm Val}}
\def\indic{\hbox{\raise-2pt \hbox{\indbf 1}}}

\def\RRR{\hbox{\msytw R}}

\def\NNN{\hbox{\msytw N}} 
\def\nnn{\hbox{\msytwww N}}
\def\ZZZ{\hbox{\msytw Z}} 
\def\zzz{\hbox{\msytwww Z}}

\def\vvv{\hbox{\msytwwww V}} \def\vvvv{\hbox{\msytwwwww V}}
\def\www{\hbox{\msytwwww W}} \def\wwww{\hbox{\msytwwwww W}}
\def\GGG{\hbox{\msytw G}}


\newcount\mgnf  
\mgnf=0 

\ifnum\mgnf=0
\def\openone{\leavevmode\hbox{\ninerm 1\kern-3.3pt\tenrm1}}%
\def\*{\vglue0.3truecm}\fi
\ifnum\mgnf=1
\def\openone{\leavevmode\hbox{\ninerm 1\kern-3.63pt\tenrm1}}%
\def\*{\vglue0.5truecm}\fi


\newcount\tipobib\newcount\bz\bz=0\newcount\aux\aux=1
\newdimen\bibskip\newdimen\maxit\maxit=0pt


\tipobib=0
\def\9#1{\ifnum\aux=1#1\else\relax\fi}

\newwrite\bib
\immediate\openout\bib=\jobname.bib
\global\newcount\bibex
\bibex=0
\def\verabib{\number\bibex}

\ifnum\tipobib=0
\def\cita#1{\expandafter\ifx\csname c#1\endcsname\relax
\hbox{$\clubsuit$}#1\write16{Manca #1 !}%
\else\csname c#1\endcsname\fi}
\def\rife#1#2#3{\immediate\write\bib{\string\raf{#2}{#3}{#1}}
\immediate\write15{\string\C(#1){[#2]}}
\setbox199=\hbox{#2}\ifnum\maxit < \wd199 \maxit=\wd199\fi}
\else
\def\cita#1{%
\expandafter\ifx\csname d#1\endcsname\relax%
\expandafter\ifx\csname c#1\endcsname\relax%
\hbox{$\clubsuit$}#1\write16{Manca #1 !}%
\else\probib(ref. numero )(#1)%
\csname c#1\endcsname%
\fi\else\csname d#1\endcsname\fi}%
\def\rife#1#2#3{\immediate\write15{\string\Cp(#1){%
\string\immediate\string\write\string\bib{\string\string\string\raf%
{\string\verabib}{#3}{#1}}%
\string\Cn(#1){[\string\verabib]}%
\string\CCc(#1)%
}}}%
\fi

\def\Cn(#1)#2{\expandafter\xdef\csname d#1\endcsname{#2}}
\def\CCc(#1){\csname d#1\endcsname}
\def\probib(#1)(#2){\global\advance\bibex+1%
\9{\immediate\write16{#1\verabib => #2}}%
}

\def\C(#1)#2{\SIA c,#1,{#2}}
\def\Cp(#1)#2{\SIAnx c,#1,{#2}}

\def\SIAnx #1,#2,#3 {\senondefinito{#1#2}%
\expandafter\def\csname#1#2\endcsname{#3}\else%
\write16{???? ma #1,#2 e' gia' stato definito !!!!}\fi}

\bibskip=10truept
\def\hboxto{\hbox to}

\catcode`\{=12\catcode`\}=12
\catcode`\<=1\catcode`\>=2
\immediate\write\bib<
        \string\halign{\string\hboxto \string\maxit%
        {\string #\string\hfill}&%
        \string\vtop{\string\parindent=0pt\string\advance\string\hsize%
        by -1.55truecm%
        \string#\string\vskip \bibskip
        }\string\cr%
>
\catcode`\{=1\catcode`\}=2
\catcode`\<=12\catcode`\>=12

\def\raf#1#2#3{\ifnum \bz=0 [#1]&#2 \cr\else
\llap{${}_{\rm #3}$}[#1]&#2\cr\fi}

\newread\bibin

\catcode`\{=12\catcode`\}=12
\catcode`\<=1\catcode`\>=2
\def\chiudibib<
\catcode`\{=12\catcode`\}=12
\catcode`\<=1\catcode`\>=2
\immediate\write\bib<}>
\catcode`\{=1\catcode`\}=2
\catcode`\<=12\catcode`\>=12
>
\catcode`\{=1\catcode`\}=2
\catcode`\<=12\catcode`\>=12

\def\makebiblio{
\ifnum\tipobib=0
\advance \maxit by 10pt
\else
\maxit=1.truecm
\fi
\chiudibib
\immediate \closeout\bib
\openin\bibin=\jobname.bib
\ifeof\bibin\relax\else
\raggedbottom
\input \jobname.bib
\fi
}

\openin13=#1.aux \ifeof13 \relax \else
\input #1.aux \closein13\fi
\openin14=\jobname.aux \ifeof14 \relax \else
\input \jobname.aux \closein14 \fi
\immediate\openout15=\jobname.aux

\def\biblio{\*\*\centerline{\titolo References}\*\nobreak\makebiblio}


\ifnum\mgnf=0
   \magnification=\magstep0
   \hsize=15truecm\vsize=21.0truecm\voffset2.truecm\hoffset.5truecm
   \parindent=0.3cm\baselineskip=0.45cm\fi
\ifnum\mgnf=1
   \magnification=\magstep1\hoffset=0.truecm
   \hsize=15truecm\vsize=21.0truecm
   \baselineskip=18truept plus0.1pt minus0.1pt \parindent=0.9truecm
   \lineskip=0.5truecm\lineskiplimit=0.1pt      \parskip=0.1pt plus1pt\fi


\mgnf=0   
\driver=1 
\openin14=\jobname.aux \ifeof14 \relax \else
\input \jobname.aux \closein14 \fi
\openout15=\jobname.aux




\centerline{\titolone Periodic solutions for completely resonant}
\vskip.2truecm
\centerline{\titolone nonlinear wave equations}

\vskip1.truecm \centerline{{\titolo
Guido Gentile$^{\dagger}$,
Vieri Mastropietro$^{*}$,
Michela Procesi$^{\star}$}}
\vskip.2truecm
\centerline{{}$^{\dagger}$ Dipartimento di Matematica,
Universit\`a di Roma Tre, Roma, I-00146}
\vskip.1truecm
\centerline{{}$^{*}$ Dipartimento di Matematica,
Universit\`a di Roma ``Tor Vergata'', Roma, I-00133}
\vskip.1truecm
\centerline{{}$^{\star}$ SISSA, Trieste, I-34014}
\vskip1truecm


\line{\vtop{
\line{\hskip1.1truecm\vbox{\advance \hsize by -2.3 truecm
\0{\cs Abstract.}
{\it 
We consider the nonlinear string equation with Dirichlet boundary
conditions $u_{xx}-u_{tt}=\f(u)$, with $\f(u)=\F u^{3} + O(u^{5})$
odd and analytic, $\F\neq0$, and we construct
small amplitude periodic solutions with frequency $\o$
for a large Lebesgue measure set of $\o$ close to $1$.
This extends previous results where only a zero-measure set of frequencies
could be treated (the ones for which no small divisors appear).
The proof is based on combining the Lyapunov-Schmidt decomposition,
which leads to two separate sets of equations dealing
with the resonant and nonresonant Fourier components,
respectively the Q and the P equations,
with resummation techniques of divergent powers series,
allowing us to control the small divisors problem.
The main difficulty with respect the nonlinear wave equations
$u_{xx}-u_{tt}+ M u = \f(u)$, $M\neq0$, is that not only
the P equation but also the Q equation is infinite-dimensional.
}} \hfill} }}
\vskip1.truecm

\*\*
\section(1,Introduction)

\0We consider the {\it nonlinear wave equation} in $d=1$ given by
$$ \cases{
u_{tt} - u_{xx} = \f(u), & \cr
u(0,t) = u(\p,t) = 0 , & \cr}
\Eq(1.1)$$
where Dirichlet boundary conditions allow us to use as a basis in
$L^{2}([0,\p])$ the set of functions $\{ \sin m x ,
m\in\NNN\}$, and $\f(u)$ is any odd analytic function $\f(u) = \F u^{3}
+ O(u^{5})$ with $\F\neq0$.
We shall consider the problem of existence of periodic solutions for
\equ(1.1), which  represents a completely resonant
case for the nonlinear wave equation as in the absence of
nonlinearities all the frequencies are resonant. 

In the finite dimensional case the problem has its analogous
in the study of periodic orbits close to elliptic equilibrium points:
results of existence have been obtained in such a case by Lyapunov \cita{L}
in the nonresonant case, by Birkhoff and Lewis \cita{BL}
in case of resonances of order greater than four,
and by Weinstein \cita{We} in case of any kind of resonances.
Systems with infinitely many degrees of freedom (as
the nonlinear wave equation, the nonlinear Schr\"odinger equation and
other PDE systems) have been studied much more recently;
the problem is much more difficult because of the presence of a
small divisors problem, which is absent in the finite dimensional case. 
For the nonlinear wave equations $u_{tt}-u_{xx}+M u=\f(u)$, 
with mass $M$ {\it strictly} positive, existence of periodic solutions
has been proved by Craig and Wayne \cita{CW1}, by P\"oschel \cita{P}
(by adapting the analogous result found by
Kuksin and P\"oschel \cita{KP} for the nonlinear Schr\"odinger equation)
and by Bourgain \cita{Bo2} (see also the review \cita{Cr}). 
In order to solve the small divisors problem one has to
require that the amplitude and frequency of the solution must belong to a
Cantor set, and the main difficulty is to prove 
that such a set can be chosen with {\it non-zero} Lebesgue measure.
We recall that for such systems also quasi-periodic solutions have been
proved to exist in \cita{KP}, \cita{P}, \cita{Bo3}
(in many other papers the case in which the coefficient $M$
of the linear term is replaced by a function depending
on parameters is considered; see for instance
\cita{Wa}, \cita{Bo1} and the reviews \cita{K1}, \cita{K2}). 

In all the quoted papers 
only non-resonant cases are considered.
Some cases with some low-order resonances between the
frequencies have been studied by Craig and Wayne \cita{CW2}.
The completely resonant case \equ(1.1) has been
studied with variational methods starting from
Rabinowitz
\cita{R1}, \cita{R2}, \cita{BN},
\cita{BCN}, \cita{FR}, where periodic solutions with
period which is a rational multiple of $\p$ have been obtained;
such solutions correspond to a zero-measure set of
values of the amplitudes. The case of irrational periods, 
which in principle could provide a large measure of values,
has been studied so far only under strong Diophantine conditions
(as the ones introduced in \cita{Ba}) which essentially remove
the small divisors problem leaving in fact again a zero-measure
set of values \cita{LS}, \cita{BP}, \cita{BB1}. 
It is however conjectured that also for $M=0$
periodic solutions should exist for a large measure set
of values of the amplitudes, see for instance \cita{K2},
and indeed we prove in this paper that this is actually the case:
the unperturbed periodic solutions with periods $\o_{j}=2\p/j$ can be
continued into periodic solutions with period $\o_{\e,j}$ close to $\o_{j}$.
We rely on the Renormalization Group approach proposed in \cita{GM},
which consists in a Lyapunov-Schmidt decomposition followed
by a tree expansion of the solution which allows us
to control the small divisors problem.

Note that a naive implementation of the Craig and Wayne approach to the
resonant case encounters some obvious difficulties.
Also in such an approach the first step consists in
a Lyapunov-Schmidt decomposition, which leads to two separate
sets of equations dealing with the resonant and nonresonant
Fourier components, respectively the Q and the P equations.
In the nonresonant case, by calling $v$ the solution of the Q equation
and $w$ the solution of the P equation, if one suppose to fix $v$
in the P equation then one gets a solution $w=w(v)$,
which inserted into the Q equation gives a closed finite-dimensional
equation for $v$. In the resonant case the Q equation is
infinite-dimensional. We can again find the solution $w(v)$
of the P equation by assuming that the Fourier coefficients of
$v$ are exponentially decreasing. However one ``consumes'' part of the
exponential decay to bound the small divisors, hence
the function $w(v)$ has Fourier coefficients decaying
with a smaller rate, and when inserted into
the Q equation the $v$ solution has not the properties
assumed at the beginning. Such a problem can be
avoided if there are no small divisors, as in such a case even a power
law decay is sufficient to find $w(v)$ and there
is no loss of smoothness; this is what is done in the
quoted papers \cita{LS}, \cita{BP} and \cita{BB1}:
the drawback is that only a zero-measure set is found.
The advantage of our approach is that we treat the P and Q
equation simultaneously, solving them together.
As in \cita{BP} and \cita{BB2} we also consider the problem of finding
how many solutions can be obtained with
given period, and we study their minimal period. 

\*

If $\f=0$ every real solution of \equ(1.1) can be written as
$$ u(x,t) = \sum_{n=1}^{\io} U_{n} \sin nx \cos (\o_n t+\th_n) ,
\Eq(1.2) $$
where $\o_{n}=n$ and $U_{n}\in\RRR$ for all $n\in\NNN$.

For $\e>0$ we set $\F=\s F$, with $\s={\rm sgn}\F$ and $F>0$,
and rescale $u \to \sqrt{\e/F} u$ in \equ(1.1), so obtaining
$$ \cases{
u_{tt} - u_{xx} = \s \e u^3 + O(\e^{2}) , & \cr
u(0,t) = u(\p,t) = 0 , & \cr}
\Eq(1.3)$$
where $O(\e^{2})$ denotes an analytic function of $u$ and $\e$ of
order at least $2$ in $\e$, and we define $\o_{\e}=\sqrt{1-\l\e}$,
with $\l\in\RRR$, so that $\o_{\e}=1$ for $\e=0$.

As the nonlinearity $\f$ is odd the solution of \equ(1.3) can be
extended in the $x$ variable to an odd $2\pi$ periodic function
(even in the variable $t$). 
We shall consider $\e$ small and we shall show
that there exists a solution of \equ(1.3), which is
$2\pi/\o_{\e}$-periodic in $t$ and $\e$-close to the function
$$ u_{0}(x,\o_{\e}t) = a_0(\o_{\e} t + x) - a_0(\o_{\e} t - x) ,
\Eq(1.4) $$
provided that $\e$ is in an appropriate Cantor set and $a_0(\x)$
is the odd $2\pi$-periodic solution of the integro-differential equation
$$ \l \ddot a_{0} = - 3 \la a_{0}^{2} \ra a_{0} - a_{0}^{3} ,
\Eq(1.5) $$
where the dot denotes derivative with respect to $\x$, and,
given any periodic function $F(\x)$ with period $T$, we denote by
$$ \la F \ra = {1\over T} \int_{0}^{T} {\rm d}\x \, F(\x)
\Eq(1.6) $$
its average. Then a $2\p/\o_{\e}$-periodic solution of \equ(1.1)
is simply obtained by scaling back the solution of \equ(1.3).

The equation \equ(1.5) has odd $2\pi$-periodic solutions,
provided that one sets $\l>0$; we shall choose $\l=1$ in the following.
An explicit computation gives \cita{BP}
$$ a_{0}(\x) = V_{\mmmm} \, \sn (\O_{\mmmm}\x,\mm)
\Eq(1.7) $$
for $\mm$ a suitable negative constant ($\mm \approx -0.2554$),
with $\O_{\mmmm}=2K(\mm)/\p$ and $V_{\mmmm}=\sqrt{-2\mm}\O_{\mmmm}$, where
$\sn(\O_{\mmmm}\x,\mm)$ is the sine-amplitude function and $K(\mm)$ is the
elliptic integral of the first kind, with modulus $\sqrt{\mm}$ \cita{GrR};
see Appendix A1 for further details.
Call $2\k$ the width of the analyticity strip of
the function $a_{0}(\x)$ and $\a$ the maximum value
it can assume in such a strip; then one has
$$ \left| a_{0,n} \right| \le \a e^{-2k|n|} .
\Eq(1.15) $$
Our result (including also the cases
of frequencies which are multiples of $\o_{\e}$)
can be more precisely stated as follows.

\*

\0{\bf Theorem.} {\it Consider the equation \equ(1.1),
where $\f(u)=\F u^{3} + O(u^{5})$ is an odd analytic function,
with $F=|\F| \neq 0$. Define $u_{0}(x,t)=a_{0}(t+x)-a_{0}(t-x)$,
with $a_{0}(\x)$ the odd  $2\pi$-periodic solution of \equ(1.5).
There are a positive constant $\e_{0}$ and for all $j\in\NNN$
a set $\EE \in[0,\e_{0}/j^{2}]$ satisfying
$$ \lim_{\e\to0} { {\rm meas}(\EE \cap [0,\e] ) \over \e} = 1 ,
\Eq(1.8) $$
such that for all $\e\in\EE$, by setting $\o_{\e}=\sqrt{1-\e}$ and
$$ \|f(x,t)\|_r= \sum_{(n,m)\in \zzz^2}f_{n,m}e^{r(|n|+|m|)},
\Eq(1.9) $$
for analytic $2\pi$-periodic functions,
there exist $2\p/j\o_{\e}$-periodic solutions
$u_{\e,j}(x,t)$ of \equ(1.1), analytic in $(t,x)$, with
$$ \left\| u_{\e,j}(x,t) - j \sqrt{\e/F} u_{0}(j x,j\o_{\e}t)
\right\|_{\k'} \le C \, j \, \e\sqrt{\e} ,
\Eq(1.10) $$
for some constants $C>0$ and $0<\k'<\k$.}

\*

Note that such a result provides a solution of the
open problem 7.4 in \cita{K2},
as far as periodic solutions are concerned.

As we shall see for $\f(u)=F u^{3}$ for
all $j\in\NNN$ one can take the set $\EE=[0,\e_{0}]$,
independently of $j$, so that for fixed $\e\in\EE$
no restriction on $j$ have to be imposed.

\*

We look for a solution of \equ(1.3) of the form
$$ \eqalign{
u(x,t) & = \sum_{(n,m) \in \zzz^{2}}
e^{in j \o t + i j m x} u_{n,m}
= v(x,t) + w(x,t) , \cr
v(x,t) & = a(\x) - a(\x') , \qquad
\x= \o t + x, \quad \x'= \o t - x , \cr a(\x) & = \sum_{n\in  \zzz} 
e^{i n \x} a_{n} , \cr
w(x,t) & = \sum_{\scriptstyle (n,m) \in  \zzz^{2} \atop
\scriptstyle |n| \neq |m|}
e^{in j \o t +i j m x} w_{n,m} , \cr}
\Eq(1.11) $$
with $\o=\o_{\e}$, such that one has $w(x,t)=0$
and $a(\x)=a_{0}(\x)$ for $\e=0$.
Of course by the symmetry of \equ(1.1), hence of \equ(1.4),
we can look for solutions (if any) which verify
$$ u_{n,m} = - u_{n,-m} = u_{-n,m} 
\Eq(1.12) $$
for all $n,m\in\ZZZ$.

Inserting \equ(1.11) into \equ(1.3) gives two sets of equations, called
the Q and P equations \cita{CW1}, which are given, respectively, by
$$ \eqalign{
\hbox{Q} \qquad & \cases{
n^{2} a_{n} = \left[ \f(v+w) \right]_{n,n} , & \cr
- n^{2} a_{n} = \left[ \f(v+w) \right]_{n,-n} , & \cr} \cr
\hbox{P} \qquad & \left( - \o^2 n^2 + m^2 \right) w_{n,m} = \e \left[
\f(v+w) \right]_{n,m} , \qquad |m| \neq |n| , \cr}
\Eq(1.13) $$
where we denote by $[F]_{n,m}$ the Fourier component of the function
$F(x,t)$ with labels $(n,m)$, so that
$$ F(x,t) = \sum_{(n,m)\in\zzz^{2}} e^{i n\o t + m x} [F]_{n,m} .
\Eq(1.14) $$
In the same way we shall call $[F]_{n}$ the Fourier component of
the function $F(\x)$ with label $n$;
in particular one has $[F]_{0}=\la F\ra$. Note also that the two
equations Q are in fact the same, by the symmetry property 
$\left[ \f(v+w) \right]_{n,m}=-\left[ \f(v+w) \right]_{n,-m}$, 
which follows from \equ(1.12).

We start by considering the case $\f(u)=u^{3}$ and $j=1$,
for simplicity. We shall discuss at the end how the other cases
can be dealt with; see Section \secc(8).

\*\*
\section(2,Lindstedt series expansion)

\0One could try to write a power series
expansion in $\e$ for $u(x,t)$, using \equ(1.13)
to get recursive equations for the coefficients.
However by proceeding in this way one finds that
the coefficient of order $k$ is given by a sum of terms
some of which of order $O(k!^\a)$, for some constant $\a$.
This is the same phenomenon occurring
in the Lindstedt series for invariant KAM tori
in the case of quasi-integrable Hamiltonian systems;
in such a case however one can show that
there are {\it cancellations} between the terms
contributing to the coefficient of order $k$,
which at the end admits a bound $C^k$, for a suitable constant $C$.
On the contrary such cancellations are absent in the present case
and we have to proceed in a different way, equivalent to a resummation
(see \cita{GM} where such procedure was applied to the same
nonlinear wave equation with a mass term, $u_{xx}-u_{tt}+M u=\f(u)$).

\*

\0{\bf Definition 1. }{\it Given a sequence
$\{\n_{m}(\e)\}_{|m|\ge 1}$, such that $\n_{m}=\n_{-m} $,
we define the {\rm renormalized frequencies} as
$$ \tilde\o_m^2 \= \o_m^2 + \n_m , \qquad \o_{m}=|m| ,
\Eq(2.1) $$
and the quantities $\n_{m}$ will be called the {\rm counterterms}.}

\*

By the above definition and the parity properties \equ(1.12)
the P equation in \equ(1.13) can be rewritten as
$$ \eqalign{
w_{n,m} \left( - \o^2 n^{2} + \tilde\o_m^2 \right)
& = \n_{m} w_{n,m} + \e [\f(v+w)]_{n,m} \cr
& = \n_{m}^{(a)} w_{n,m} + \n_{m}^{(b)} w_{n,-m} + \e [\f(v+w)]_{n,m}, \cr}
\Eq(2.2) $$
where
$$\n_{m}^{(a)}-\n_{m}^{(b)}=\n_{m} .
\Eq(2.3) $$

With the notations of \equ(1.14), and recalling that we are
considering $\f(u)=u^{3}$, we can write
$$ \eqalign{
\left[ (v+w)^3 \right]_{n,n} & = [v^{3}]_{n,n}+[w^{3}]_{n,n}+
3 [v^{2} w]_{n,n} + 3 [w^{2} v]_{n,n} \cr
& \equiv [v^{3}]_{n,n} + [g(v,w)]_{n,n} , \cr}
\Eq(2.5) $$
where, again by using the parity properties \equ(1.12),
$$ [v^{3}]_{n,n} = [a^{3}]_{n} + 3 \la a^{2} \ra a_{n} .
\Eq(2.6) $$
Then the first Q equation in \equ(1.12) can be rewritten as
$$ n^{2} a_{n} = [a^{3}]_{n} + 3 \la a^{2} \ra a_{n} +
\left[ g(v,w) \right]_{n,n} ,
\Eq(2.7) $$
so that $a_{n}$ is the Fourier coefficient of the
$2\p$-periodic solution of the equation
$$  \ddot a = - \left( a^{3} + 3 \la a^{2} \ra a + G(v,w) \right) ,
\Eq(2.8) $$ 
where we have introduced the function
$$ G(v,w) = \sum_{n\in\zzz}
e^{in\x} \left[ g(v,w) \right]_{n,n} .
\Eq(2.9) $$
To study the equations \equ(2.2) and \equ(2.7)
we introduce an auxiliary parameter $\m$, which
at the end will be set equal to 1, by writing \equ(2.2) as
$$ w_{n,m} \left( - \o^2 n^{2} + \tilde\o_m^2 \right)
= \m \n_{m}^{(a)} w_{n,m} + \m \n_{m}^{(b)} w_{n,-m} +
\m \e [\f(v+w)]_{n,m} ,
\Eq(2.2a) $$
and we shall look for $u_{n,m}$ in the form of a 
power series expansion in $\m$,
$$ u_{n,m} = \sum_{k=0}^{\io} \m^{k} u_{n,m}^{(k)} ,
\Eq(2.4) $$
with $u_{n,m}^{(k)}$ depending on $\e$ and on
the parameters $\n_{m}^{(c)}$, with $c=a,b$ and $|m|\ge 1$.

In \equ(2.4) $k=0$ requires $u_{n,\pm n}=\pm a_{0,n}$
and $u_{n,m}=0$ for $|n|\neq|m|$, while for $k\ge 1$,
as we shall see later on, the dependence on the
parameters $\n_{m'}^{(c')}$ will be of the form
$$ \prod_{m'=2}^{\io} \prod_{c'=a,b}
\left( \n_{m'}^{(c')} \right)^{k_{m'}^{(c')}} ,
\Eq(2.5a) $$
with $|\ul{k}|=k_{1}^{(a)}+k_{1}^{(b)}+k_{2}^{(a)}+k_{2}^{(b)}+\ldots
\le k-1$. Of course we are using the symmetry property to restrict
the dependence only on the positive labels $m'$.

We derive recursive equations for the coefficients $u_{n,m}^{(k)}$
of the expansion. We start from the coefficients with $|n|=|m|$.

By \equ(1.11) and \equ(2.4) we can write
$$  a = a_{0} + \sum_{k=1}^{\io} \m^{k} A^{(k)} ,
\Eq(2.10) $$
and inserting this expression into \equ(2.8)
we obtain for $A^{(k)}$ the equation
$$ \ddot A^{(k)}  = 
- 3\left(a_0^2 A^{(k)} + 
\la a_0^{2} \ra A^{(k)} 
+ 2 \la a_0 A^{(k)} \ra a_0\right) + f^{(k)} ,
\Eq(2.11) $$ 
with 
$$ f^{(k)} = - \sum_{\scriptstyle k_{1}+k_{2}+k_{3}=k \atop
\scriptstyle k_{i}=k \to |n_{i}| \neq |m_{i}|}
\sum_{\scriptstyle n_{1}+n_{2}+n_{3} = n \atop
\scriptstyle m_{1}+m_{2}+m_{3} = m}
u^{(k_{1})} _{n_{1},m_{1}}
u^{(k_{2})} _{n_{2},m_{2}}
u^{(k_{3})} _{n_{3},m_{3}} ,
\Eq(2.12)$$
where we have used the notations
$$ u^{(k)}_{n,m} = \cases{
v^{(k)}_{n,m} , & if $|n|=|m|$ , \cr
w^{(k)}_{n,m} , & if $|n|\neq |m|$ , \cr}
\Eq(2.13) $$
with
$$ v^{(k)}_{n,n}= \cases{
A^{(k)}_{n} , & if $k\not=0$, \cr
a_{0,n} , & if $k= 0$, \cr} \quad
v^{(k)}_{n,-n}= \cases{
- A^{(k)}_{n} , & if $k\not=0$, \cr
- a_{0,n} , & if $k= 0$ . \cr}
\Eq(2.14) $$

Before studying how to find the solution of this equation
we introduce some preliminary definitions.
To shorten notations we write
$$ c(\x)\=\cn(\O_{\mmmm}\x,\mm),\quad 
s(\x)\=\sn(\O_{\mmmm}\x,\mm)\quad d(\x)\=\dn(\O_{\mmmm}\x,\mm) ,
\Eq(2.15) $$
and set $cd(\x)=\cn(\O_{\mmmm}\x,\mm)\,\dn(\O_{\mmmm}\x,\mm)$.
Moreover given an analytic periodic function $F(\x)$ we define
$$\PP[F](\x)  = F(\x)-\la F \ra .
\Eq(2.16) $$
and we introduce a linear operator $\II$ acting
on $2\p$-periodic zero-mean functions
and defined by its action on the basis
$e_{n}(\x)=e^{i n\x}$, $n\in\ZZZ\setminus\{0\}$,
$$\II[e_n](\x)= {e_n(\x)\over i n} .
\Eq(2.17) $$
Note that if $\PP[F]=0$ then $\PP[\II [F]]=0$
($\II[F]$ is simply the zero-mean primitive of $F$);
moreover $\II$ switches parities. 

In order to find an odd solution of \equ(2.11) we replace first 
$\la a_0 A^{(k)}\ra$ with a parameter $C^{(k)}$,
and we study the modified equation
$$ \ddot A^{(k)}  = 
- 3\left(a_0^2 A^{(k)} + 
\la a_0^{2} \ra A^{(k)} 
+ 2 C^{(k)} a_0 \right) + f^{(k)} .
\Eq(2.18) $$ 
Then we have the following result (proved in Appendix \secc(A1)).

\*

\0{\bf Lemma 1.} {\it Given an odd analytic $2\pi$-periodic
function $h(\x)$, the equation
$$ \ddot y = - 3 \left( a_0^2 +\la a_0^{2} \ra \right) y + h ,
\Eq(2.19)$$
admits one and only one odd analytic $2\pi$-periodic
solution $y(\x)$, given by
$$ y = \LLL[h] \=
B_{\mmmm} \Big( \O_{\mmmm}^{-2}D_{\mmmm}^{2} s \la s\, h\ra  +
\O_{\mmmm}^{-1}D_{\mmmm} \left( s \, \II[cd\, h]-
cd \,\II [\PP[ s\, h]] \right) + cd \, \II [\II [cd \, h]] \Big) ,
\Eq(2.20) $$
with $B_{\mmmm}=-\mm/(1-\mm)$ and $D_{\mmmm}=-1/\mm$.}

\*

As $a_{0}$ is analytic and odd, we find immediately,
by induction on $k$ and using lemma 1,
that $f^{(k)}$ is analytic and odd, and that
the solution of the equation \equ(2.18) is odd and given by
$$ \tilde A^{(k)} = 
\LLL[ - 6 C^{(k)} a_0 + f^{(k)} ] .
\Eq(2.21) $$  

The function $\tilde A^{(k)}$ so found depends of
course on the parameter $C^{(k)}$; in order
to obtain $\tilde A^{(k)}= A^{(k)}$,
we have to impose the constraint
$$ C^{(k)} = \la a_0 A^{(k)} \ra ,
\Eq(2.22) $$
and by \equ(2.21) this gives
$$ C^{(k)} = - 6 C^{(k)}  
\la a_0 L[a_0]\ra +  \la a_0 \LLL[f^{(k)} ] \ra ,
\Eq(2.23) $$
which can be rewritten as
$$ \left( 1 + 6 \la a_0 \LLL[a_0]\ra \right) C^{(k)} =
\la a_0 \LLL[f^{(k)} ]\ra .
\Eq(2.24) $$ 
An explicit computation (see Appendix \secc(A3)) gives
$$\la a_0 \LLL[a_0]\ra=  
{1 \over 2} V_{\mmmm}^{2} \O_{\mmmm}^{-2} B_{\mmmm}
 \left( \left( 2D_{\mmmm} - {1\over2} \right)
\langle s^{4} \rangle  + \left( 2D_{\mmmm} \left( D_{\mmmm}-1 \right)
+ {1\over 2} \right) \langle s^{2} \rangle^{2} \right) ,
\Eq(2.25) $$
which yields $r_0 = (1+6 \langle a_0 \LLL[a_0]\rangle ) \neq 0$.
At the end we obtain the recursive definition
$$ \cases{
A^{(k)} = \LLL[ f^{(k)} - 6 C^{(k)} a_0 ] , & \cr
& \cr
C^{(k)} = r_0^{-1} \la a_0 \LLL[ f^{(k)} ] \ra . & \cr}
\Eq(2.26) $$
In Fourier space the first of \equ(2.26) becomes
$$ \eqalign{
A^{(k)}_{n} & = B_{\mmmm} \O_{\mmmm}^{-2} D_{\mmmm}^{2}
s_{n} \mathop{\sum}_{n_1+n_2=0}
s_{n_{1}} \left( f^{(k)}_{n_{2}} -
6 C^{(k)} a_{0,n_{2}} \right) \cr
& \qquad + B_{\mmmm} {\mathop{\sum}_{n_{1}+n_{2}+n_{3}=n}}^{*}
{1 \over i^2 (n_{2}+n_{3})^2}
cd_{n_{1}} cd_{n_{2}} \left( f^{(k)}_{n_{3}} -
6 C^{(k)} a_{0,n_{3}} \right)  \cr
& \qquad + B_{\mmmm}\O_{\mmmm}^{-1}D_{\mmmm}
{\mathop{\sum}_{n_{1}+n_{2}+n_{3}=n}}^{*}
{1 \over i (n_{2}+n_{3})} s_{n_{1}} cd_{n_{2}}
\left( f^{(k)}_{n_{3}} -
6 C^{(k)} a_{0,n_{3}} \right)  \cr
& \qquad - B_{\mmmm} \O_{\mmmm}^{-1}D_{\mmmm} 
{\mathop{\sum}_{n_{1}+n_{2}+n_{3}=n}}^{*} 
{1 \over i (n_{2}+n_{3})} 
cd_{n_{1}} s_{n_{2}} \left( f^{(k)}_{n_{3}} -
6 C^{(k)} a_{0,n_{3}} \right)  \cr
& \= \sum_{n'} \LLL_{n n'} \left( f^{(k)}_{n'} -
6 C^{(k)} a_{0,n'} \right)  , \cr}
\Eq(2.27) $$
where the constants $B_{\mmmm}$ and $D_{\mmmm}$ are defined after \equ(2.20),
and the $*$ in the sums means that one has the constraint
$n_{2}+n_{3}\neq 0$, while the second of \equ(2.26) can be written as
$$ C^{(k)} = r_0^{-1} \sum_{n,n'\in\zzz}
a_{0,-n} \LLL_{n,n'} f^{(k)}_{n'} .
\Eq(2.28) $$

Now we consider the coefficients $u^{(k)}_{n,m}$
with $|n|\neq |m|$. The coefficients $w^{(k)}_{n,m}$
verify the recursive equations
$$ w^{(k)}_{n,m} \left[ - \o^2 n^{2} + \tilde\o_m^2 \right] = 
\n_{m}^{(a)} w_{n,m}^{(k-1)} + \n_{m}^{(b)} w_{n,-m}^{(k-1)} +
[(v+w)^3]^{(k-1)}_{n,m} ,
\Eq(2.29) $$
where
$$ [(v+w)^3]^{(k)}_{n,m}  = \sum_{k_{1}+k_{2}+k_{3}=k}
\sum_{\scriptstyle n_{1}+n_{2}+n_{3} = n \atop
\scriptstyle m_{1}+m_{2}+m_{3} = m} u^{(k_{1})}_{n_{1},m_{1}}
u^{(k_{2})}_{n_{2},m_{2}} u^{(k_{3})}_{n_{3},m_{3}} ,
\Eq(2.30) $$
if we use the same notations \equ(2.13) and \equ(2.14) as in \equ(2.12).

The equations \equ(2.27) and \equ(2.29), together with
\equ(2.30), \equ(2.12), \equ(2.28) and \equ(2.30),
define recursively the coefficients $u^{(k)}_{n,m}$.

To prove the theorem we shall proceed in two steps.
The first step consists in looking for the solution
of the equations \equ(2.27) and \equ(2.29)
by considering $\tilde\o=\{\tilde\o_{m}\}_{|m|\ge 1}$
as a given set of parameters satisfying the Diophantine conditions
(called respectively the first and the second Mel$^{\prime}$nikov conditions)
$$ \eqalignno{
& \left| \o n \pm \tilde \o_{m} \right|
\ge C_{0} |n|^{-\t}
\qquad \forall n\in\ZZZ \setminus 
\{ 0 \} \hbox{ and }
\forall m \in \ZZZ \setminus\{0\} 
\hbox{ such that } |m| \neq |n| , \cr
& \left| \o  n \pm \left( \tilde\o_{m} \pm \tilde\o_{m'}
\right) \right| \ge C_{0} |n|^{-\t}
\qquad
& \eq(2.31) \cr
& \qquad \qquad \qquad \qquad
\forall n\in\ZZZ \setminus \{ 0 \} \hbox{ and }
\forall m,m' \in \ZZZ \setminus\{0\}  
\hbox{ such that } |n|\neq |m \pm m'| , \cr} $$
with positive constants $C_{0},\t$.
We shall prove in Sections \secc(3) to \secc(5) the following result.

\*

\0{\bf Proposition 1.} {\it Let be given a sequence
$\tilde\o=\{\tilde\o_m\}_{|m| \ge 1}$ verifying \equ(2.31),
with $\o=\o_{\e}=\sqrt{1-\e}$
and such that $|\tilde\o_{m}-|m||\le C\e/|m|$ for some constant $C$.
For all $\m_{0}>0$ there exists $\e_{0}>0$ such that
for $|\m|\le \m_{0}$ and $0<\e<\e_{0}$ there is
a sequence $\n(\tilde\o,\e;\m) = \{\n_{m}(\tilde\o,\e;\m)\}_{|m| \ge 1}$,
where each $\n_{m}(\tilde\o,\e;\m)$ is analytic in $\m$,
such that the coefficients $u^{(k)}_{n,m}$
which solve \equ(2.27) and \equ(2.29) define via \equ(2.4)
a function $u(x,t;\tilde\o,\e;\m)$ which is analytic in $\m$,
analytic in $(x,t)$ and $2\pi$-periodic in $t$ and solves
$$ \eqalign{
\hbox{} \qquad &
\cases{
n^{2} a_{n} = \left[ a^{3} \right]_{n,n} +
3 \la a^{2} \ra a_{n} + \left[ g(v,w) \right]_{n,n} , & \cr
- n^{2} a_{n} = \left[ a^{3} \right]_{n,-n} +
3 \la a^{2} \ra a_{-n} + \left[ g(v,w) \right]_{n,-n} , & \cr} \cr
\hbox{} \qquad & \left( - \o^2 n^2 + \tilde\o_{m}^2 \right)
w_{n,m} = \m \n_{m}(\tilde\o,\e;\m) \, w_{n,m} +
\m \e \left[ \f(v+w) \right]_{n,m} , \qquad |m| \neq |n| , \cr}
\Eq(2.32) $$
with the same notations as in \equ(1.13).}

\*

If $\t \le 2$ then one can require only the first Mel'nikov conditions
in \equ(2.31), as we shall show in Section \secc(7).

Then in Proposition 1 one can fix $\m_{0}=1$, so that
one can choose $\m=1$ and set
$u(x,t;\tilde\o,\e)=u(x,t;\tilde\o,\e;1)$ and
$\n_{m}(\tilde\o,\e)=\n_{m}(\tilde\o,\e;1)$.

The second step, to be proved in Section \secc(6),
consists in inverting \equ(2.1), with $\n_{m}=\n_{m}(\tilde\o,\e)$
and $\tilde\o$ verifying \equ(2.31). This requires
some preliminary conditions on $\e$, given by
the Diophantine conditions
$$ \left| \o n \pm m \right|
\ge C_{1} |n|^{-\t_{0}}
\qquad \forall n\in\ZZZ \setminus 
\{ 0 \} \hbox{ and }
\forall m \in \ZZZ \setminus\{0\} 
\hbox{ such that } |m| \neq |n| ,
\Eq(2.33)$$
with positive constants $C_{1}$ and $\t_{0}>1$.

This allows to solve iteratively \equ(2.1),
by imposing further non-resonance conditions besides \equ(2.33),
provided that one takes $C_{1}=2C_{0}$ and $\t_{0}<\t-1$,
which requires $\t>2$. At each iterative step one has to exclude
some further values of $\e$, and at the end the left values 
fill a Cantor set $\EE$ with large relative measure in $[0,\e_{0}]$
and $\tilde\o$ verify \equ(2.33).

If $1<\t \le 2$ the first Mel'nikov conditions, which, as we said
above, become sufficient to prove Proposition 1,
can be obtained by requiring \equ(2.33) with $\t_{0}=\t$;
again this leaves a large measure set of allowed values of $\e$.
This is discussed in Section \secc(7).

The result of this second step can be summarized as follows.

\*

\0{\bf Proposition 2.} {\it
There are $\d>0$ and a set $\EE\subset [0,\e_{0}]$
with complement of relative Lebesgue measure of order $\e_{0}^{\d}$
such that for all $\e\in\EE$ there exists
$\tilde\o=\tilde\o(\e)$ which solves \equ(2.1)
and satisfy the Diophantine conditions \equ(2.31) with
$|\tilde\o_{m}-|m||\le C\e/|m|$ for some constant $C$.}

\*

As we said, our approach is based on constructing the periodic solution
of the string equation by a perturbative expansion
which is the analogue of the {\it Lindstedt series}
for (maximal) KAM invariant tori
in finite-dimensional Hamiltonian systems. 
Such an approach immediately encounters a difficulty;
while the invariant KAM tori are analytic in the perturbative
parameter $\e$, the periodic solutions we are looking for are {\it not}
analytic; hence a power series construction seems at first sight
hopeless. Nevertheless it turns out that the Fourier
coefficients of the periodic solution have the form 
$u_{n,m}(\tilde\o(\o,\e),\e;\m)$; while such functions are
not analytic in $\e$,  it turns out to be analytic
in $\m$, provided that $\tilde\o$ satisfies the
condition \equ(2.31) and $\e$ is small enough;
this is the content of Proposition 1. 
As far as we know, in the case $M\neq0$, besides \cita{GM},
the smoothness in $\e$ at fixed $\tilde\o$ was only discussed
in a paper by Bourgain \cita{Bo2}. This smoothness is what
allow us to write as a series expansion $u_{n,m}(\tilde\o,\e;\m)$;
this strategy was already applied in \cita{GM} in the massive case.

\*\*
\section(3,Tree expansion: the diagrammatic rules)

\0A (connected) graph $\GGG$ is a collection of points (vertices)
and lines connecting all of them. The points of a graph
are most commonly known as graph vertices, but may also be
called nodes or points. Similarly, the lines connecting
the vertices of a graph are most commonly
known as graph edges, but may also be called branches or
simply lines, as we shall do. We denote with
$P(\GGG)$ and $L(\GGG)$ the set of vertices and the set of lines,
respectively. A path between two vertices is a subset of $L(\GGG)$
connecting the two vertices. A graph is planar if it can be drawn
in a plane without graph lines crossing
(i.e. it has graph crossing number 0).

\*

\0{\bf Definition 2.}
{\it A {\rm tree} is a planar graph $\GGG$ containing no closed loops
(cycles); in other words, it is a connected acyclic graph.
One can consider a tree $\GGG$ with a single special vertex $\vvv_{0}$:
this introduces a natural partial ordering on the set
of lines and vertices, and one can imagine that each line
carries an arrow pointing toward the vertex $\vvv_{0}$.
We can add an extra (oriented) line $\ell_{0}$
connecting the special vertex $\vvv_{0}$ to another point which
will be called the {\rm root} of the tree; the added line will
be called the {\rm root line}. In this way we obtain
a {\rm rooted tree} $\th$ defined by $P(\th)=P(\GGG)$
and $L(\th)=L(\GGG)\cup\ell_{0}$. A {\rm labeled tree} is a
rooted tree $\th$ together with a label function defined on
the sets $L(\th)$ and $P(\th)$.}

\*

Note that the definition of rooted tree
given above is sligthly different from the one which is
usually adopted in literature \cita{GoR}, \cita{HP}
according to which a rooted tree is just
a tree with a privileged vertex, without any extra line.
However the modified definition that we gave will
be more convenient for our purposes.
In the following we shall denote with the symbol $\th$ both
rooted trees and labels rooted trees, when no confusion arises

\*

We shall call equivalent two rooted trees which can be transformed
into each other by continuously deforming the lines in the plane
in such a way that the latter do not cross each other
(i.e. without destroying the graph structure).
We can extend the notion of equivalence also to labeled trees,
simply by considering equivalent two labeled trees if they
can be transformed into each other in such a way that also
the labels match.

Given two points $\vvv,\www\in P(\th)$, we say that $\vvv \prec \www$
if $\vvv$ is on the path connecting $\www$ to the root line.
We can identify a line with the points it connects;
given a line $\ell=(\vvv,\www)$ we say that $\ell$
enters $\vvv$ and comes out of $\www$.

In the following we shall deal mostly with labeled trees:
for simplicity, where no confusion can arise, we shall call them just
trees. We consider the following {\it diagramamtic rules}
to construct the trees we have to deal with; this
will implicitly define also the label function.

\*

\0(1) We call {\it nodes} the vertices such that there is at least
one line entering them. We call {\it end-points} the vertices
which have no entering line. We denote with $L(\th)$, $V(\th)$
and $E(\th)$ the set of lines, nodes and end-points, respectively.
Of course $P(\th)=V(\th)\cup E(\th)$.

\*

\0(2) There can be two types of lines, $w$-lines and $v$-lines,
so we can associate to each line $\ell\in L(\th)$ a {\it badge}
label $\g_{\ell}\in\{v,w\}$ and a {\it momentum}
$(n_{\ell},m_{\ell})\in\ZZZ^{2}$,
to be defined in item (8) below. If $\g_{\ell}=v$ one
has $|n_{\ell}|=|m_{\ell}|$, while if $\g_{\ell}=w$
one has $|n_{\ell}|\neq |m_{\ell}|$. One can have
$(n_{\ell},m_{\ell})=(0,0)$ only if $\ell$ is a $v$-line.
To the $v$-lines $\ell$ with $n_{\ell}\neq 0$ we also associate
a label $\d_{\ell}\in\{1,2\}$. All the lines coming out from
the end-points are $v$-lines.

\*

\0(3) To each line $\ell$ coming out from a node
we associate a {\it propagator}
$$ g_{\ell} \= g (\o n_{\ell},m_{\ell}) =
\cases{
\displaystyle
{1\over - \o^{2} n_{\ell}^2 +\tilde\o_{m_{\ell}}^2} , &
if $\g_{\ell}=w$, \cr
& \cr
\displaystyle {1 \over (i n_{\ell})^{\d_{\ell}} } , &
if $\g_{\ell}=v$, $n_{\ell}\neq 0$ , \cr
& \cr
\displaystyle 1 , &
if $\g_{\ell}=v$, $n_{\ell} = 0$ , \cr}
\Eq(3.1) $$
with momentum $(n_{\ell},m_{\ell})$.
We can associate also a propagator to the lines $\ell$ coming out
from end-points, simply by setting $g_{\ell}=1$.

\*

\0(4) Given any node $\vvv\in V(\th)$ denote with $s_{\vvvv}$ the
number of entering lines ({\it branching number}): one can have
only either $s_{\vvvv}=1$ or $s_{\vvvv}=3$.
Also the nodes $\vvv$ can be of $w$-type and $v$-type: we say
that a node is of $v$-type if the line $\ell$
coming out from it has label $\g_{\ell}=v$;
analogously the nodes of $w$-type are defined.
We can write $V(\th)=V_{v}(\th)\cup V_{w}(\th)$, 
with obvious meaning of the symbols; we also call
$V_{w}^{s}(\th)$, $s=1,3$, the set of nodes in $V_{w}(\th)$ with $s$
entering lines, and analogously we define $V_{v}^{s}(\th)$, $s=1,3$.
If $\vvv\in V_{v}^{3}$ and two entering lines come out of end points
then the remaining line entering $\vvv$ has to be a $w$-line.
If $\vvv\in V_{w}^{1}$ then the line entering $\vvv$ has to be
a $w$-line. If $\vvv\in V_{v}^{1}$ then its entering
line comes out of an end-node.

\*

\0(5) To the nodes $\vvv$ of $v$-type we associate a label
$j_{\vvv}\in \{1,2,3,4\}$ and, if $s_{\vvvv}=1$,
an {\it order} label $k_{\vvvv}$,
with $k_{\vvvv}\ge 1$.
Moreover we associate to each node $\vvv$ of $v$-type two {\it mode} labels
$(n_{\vvvv}',m_{\vvvv}')$, with $m_{\vvvv}'=\pm n_{\vvvv}'$,
and $(n_{\vvvv},m_{\vvvv})$, with $m_{\vvvv}=\pm n_{\vvvv}$,
and such that one has
$$ { m_{\vvvv}' \over n_{\vvvv}' } = { m_{\vvvv} \over n_{\vvvv}}
= { \displaystyle \sum_{i=1}^{s_{\vvvv}} m_{\ell_{i}} \over
\displaystyle \sum_{i=1}^{s_{\vvvv}} n_{\ell_{i}} } ,
\Eq(3.1a) $$
where $\ell_{i}$ are the lines entering $\vvv$.
We shall refer to them as the {\it first mode} label
and the {\it second mode} label, respectively.
To a node $\vvv$ of $v$-type we associate also
a {\it node factor} $\h_{\vvvv}$ defined as
$$ \h_{\vvvv} = \cases{
- B_{\mmmm} \O_{\mmmm}^{-2} D_{\mmmm}^{2} s_{n_{\vvvv}'} s_{n_{\vvvv}} , &
if $j_{\vvvv}=1$ and $s_{\vvvv}=3$ , \cr
- 6 B_{\mmmm} \O_{\mmmm}^{-2} D_{\mmmm}^{2} s_{n_{\vvvv}'} s_{n_{\vvvv}}
C^{(k_{\vvvv})} , &
if $j_{\vvvv}=1$ and $s_{\vvvv}=1$ , \cr
- B_{\mmmm} cd_{n_{\vvvv}'} cd_{n_{\vvvv}} , &
if $j_{\vvvv}=2$ and $s_{\vvvv}=3$ , \cr
- 6 B_{\mmmm} cd_{n_{\vvvv}'} cd_{n_{\vvvv}}
C^{(k_{\vvvv})} , &
if $j_{\vvvv}=2$ and $s_{\vvvv}=1$ , \cr
- B_{\mmmm} \O_{\mmmm}^{-1} D_{\mmmm} s_{n_{\vvvv}'} cd_{n_{\vvvv}} , &
if $j_{\vvvv}=3$ and $s_{\vvvv}=3$ , \cr
- 6 B_{\mmmm} \O_{\mmmm}^{-1} D_{\mmmm} s_{n_{\vvvv}'} cd_{n_{\vvvv}}
C^{(k_{\vvvv})} , &
if $j_{\vvvv}=3$ and $s_{\vvvv}=1$ , \cr
B_{\mmmm} \O_{\mmmm}^{-1} D_{\mmmm} cd_{n_{\vvvv}'} s_{n_{\vvvv}} , &
if $j_{\vvvv}=4$ and $s_{\vvvv}=3$ , \cr
6 B_{\mmmm} \O_{\mmmm}^{-1} D_{\mmmm} cd_{n_{\vvvv}'} s_{n_{\vvvv}}
C^{(k_{\vvvv})}, &
if $j_{\vvvv}=4$ and $s_{\vvvv}=1$ . \cr}
\Eq(3.2) $$
Note that the factors $C^{(k_{\vvvv})}=r_{0}^{1}\la a_0
\LLL[f^{(k_{\vvvv})}]\ra$
depend on the coefficients $u^{(k')}_{n,m}$, with $k'<k$,
so that they have to been defined iteratively.
The label $\d_{\ell}$ of the line $\ell$ coming out from
a node $\vvv$ of $v$-type is related to the label $j_{\vvvv}$ of $v$:
if $j_{\vvvv}=1$ then $n_{\ell}=0$, while if $j_{\vvvv}>1$ then
$n_{\ell} \neq 0$ and $\d_{\ell}=1+\d_{j_{\vvvv},2}$,
where $\d_{j_{\vvvv},1}$ denotes the Kronecker delta
(so that $\d_{\ell}=2$ if $j_{\vvvv}+2$ and $\d_{\ell}=1$ otherwise).

\*

\0(6) To the nodes $\vvv$ of $w$-type we simply associate
a {\it node factor} $\h_{\vvvv}$ given by
$$ \h_{\vvvv} = \cases{
\e, & if $s_{\vvvv}=3$, \cr
\n^{(c)}_{m_{\ell}} , & if $s_{\vvvv}=1$. \cr}
\Eq(3.3) $$
In the latter case $(n_{\ell},m_{\ell})$
is the momentum of the line coming out from $\vvv$,
and one has $c=a$ if the momentum of the entering line is
$(n_{\ell},m_{\ell})$ and 
$c=b$ if the momentum of the entering line is
$(n_{\ell},-m_{\ell})$; we shall call $\n$-vertices
the nodes $\vvv$ of $w$-type with $s_{\vvvv}=1$.
In order to unify notations we can associate also
to the nodes $\vvv$ of $w$-type two mode labels, by setting
$(n_{\vvvv}',m_{\vvvv}')=(0,0)$ and $(n_{\vvvv},m_{\vvvv})=(0,0)$.

\*

\0(7) To the end-points $\vvv$ we associate only a first mode
label $(n_{\vvvv}',m_{\vvvv}')$, with $|m_{\vvvv}'|=|n_{\vvvv}'|$,
and an {\it end-point factor}
$$ V_{\vvvv} = (-1)^{1+\d_{n_{\vvvv}',m_{\vvvv}'}} a_{0,n_{\vvvv}'}
= a_{0,m_{\vvvv}} .
\Eq(3.4) $$
The line coming out from an end-point has to be a $v$-line.

\*

\0(8) The momentum $(n_{\ell},m_{\ell})$ of a line $\ell$ is
related to the mode labels of the nodes preceding $\ell$;
if a line $\ell$ comes out from a node $v$ one has
$$ \eqalign{
n_{\ell} & = n_{\vvvv} + \sum_{\scriptstyle \wwww \in V(\th)
\atop \scriptstyle \wwww \prec \vvvv} \left( n_{\wwww}' + n_{\wwww} \right)
+ \sum_{\scriptstyle \wwww \in E(\th) \atop
\scriptstyle \wwww \prec \vvvv} n_{\wwww}' , \cr
m_{\ell} & = m_{\vvvv} + \sum_{\scriptstyle \wwww \in V(\th)
\atop \scriptstyle \wwww \prec \vvvv}
\left( m_{\wwww}' + m_{\wwww} \right) +
\sum_{\scriptstyle \wwww \in E(\th)
\atop \scriptstyle \wwww \prec \vvvv} m_{\wwww}' , \cr}
\Eq(3.5) $$
where some of the mode labels can be vanishing according
to the notations introduced above. If $\ell$ comes out
from an end-point we set $(n_{\ell},m_{\ell})=(0,0)$.

\*

We define $\Theta^{*(k)}_{n,m}$ as the {\it set of
unequivalent labeled trees}, formed by following the rules
(1) to (8) given above, and with the further following constraints:
\\
(i) if $(n_{\ell_{0}},m_{\ell_{0}})$ denotes the momentum
flowing through the root line $\ell_{0}$ and
$(n_{\vvvv_{0}}',m_{\vvvv_{0}}')$ is the first mode label associated
to the node $\vvv_{0}$ which $\ell_{0}$ comes out from
(special vertex), then one has
$n=n_{\ell_{0}} + n_{\vvvv_{0}}'$ and $m=m_{\ell_{0}} + m_{\vvvv_{0}}'$;
\\
(ii) one has
$$ k = \left| V_{w}(\th) \right| + 
\sum_{\vvvv \in V_{v}^{1}(\th)} k_{\vvvv} ,
\Eq(3.6) $$
with $k$ called the {\it order} of the tree.

An example of tree is given in Figure 3.1, where only the labels $v,w$
of the nodes have been explicitly written.

\midinsert
\*
\insertplotttt{256pt}{151pt}{
\ins{30pt}{69pt}{$v$}
\ins{67pt}{52pt}{$w$}
\ins{67pt}{87pt}{$w$}
\ins{104pt}{33pt}{$w$}
\ins{104pt}{105pt}{$w$}
\ins{140pt}{15pt}{$v$}
\ins{140pt}{52pt}{$w$}
\ins{140pt}{105pt}{$v$}
\ins{140pt}{124pt}{$v$}
\ins{176pt}{15pt}{$w$}
\ins{176pt}{69pt}{$w$}
\ins{212pt}{33pt}{$v$}
\ins{212pt}{15pt}{$w$}
}
{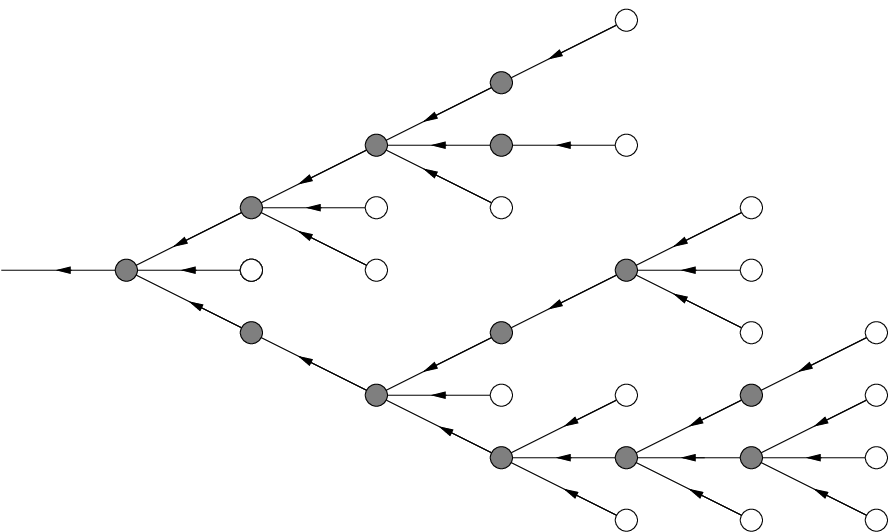}
\*
\line{\vtop{\line{\hskip1.5truecm\vbox{\advance\hsize by -3.1 truecm
\centerline{\nota
\figura(3.1)}
} \hfill} }}
\*
\endinsert

\*

\0{\bf Definition 3.}
{\it For all $\th\in \Theta^{*(k)}_{n,m}$, we call
$$ \Val(\th) = \Big( \prod_{\ell \in L(\th)} g_{\ell} \Big)
\Big( \prod_{\vvvv \in V(\th)}  \h_{\vvvv} \Big)
\Big( \prod_{\vvvv \in E(\th) } V_{\vvvv} \Big) ,
\Eq(3.8) $$
the {\rm value} of the tree $\th$.}

\*

Then the main result about the formal expansion of the
solution is provided by the following result. 

\*

\0{\bf Lemma 2.} {\it We can write
$$ u^{(k)}_{n,m} =
\sum_{\th\in \Theta^{*(k)}_{n,m}} \Val(\th) ,
\Eq(3.9) $$
and if the root line $\ell_{0}$ is a $v$-line the tree value
is a contribution to $v^{(k)}_{n,\pm n}$, while if
$\ell_{0}$ is a $w$-line the tree value is
a contribution to $w^{(k)}_{n,m}$.
The factors $C^{(k)}$ are defined as
$$ C^{(k)} = r_{0}^{-1}
{\mathop{\sum}_{\th\in \Theta^{*(k)}_{n,n}}}^{*} a_{0,-n} \Val(\th) ,
\Eq(3.10) $$
where the $*$ in the sum means the extra constraint
$s_{\vvvv_{0}}=3$ for the node immediately preceding the root
(which is the special vertex of the rooted tree).}

\*

\0{\it Proof.} 
The proof is done by induction in $k$.
Imagine to represent graphically $a_{0,n}$ as a (small) white bullet
with a line coming out from it, as in Figure 3.2a, and
$u^{(k)}_{n,m}$, $k\ge 1$, as a (big) black bullet
with a line coming out from it, as in Figure 3.2b.

\midinsert
\*
\insertplotttt{145pt}{21pt}{
\ins{20pt}{-5pt}{a}
\ins{110pt}{-5pt}{b}
}
{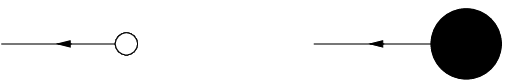}
\*
\line{\vtop{\line{\hskip1.5truecm\vbox{\advance\hsize by -3.1 truecm
\centerline{\nota
\figura(3.2)}
} \hfill} }}
\*
\endinsert

One should imagine that labels $k,n,m$ are associated
to the black bullet representing $u^{(k)}_{n,m}$,
while a white bullet representing $a_{0,n}$ carries the
labels $n,m=\pm n$.

For $k=1$ the proof of \equ(3.9) and \equ(3.10) is just a check
from the diagrammatic rules and the recursive
definitions \equ(2.25) and \equ(2.27), and it can be performed as follows.

Consider first the case $|n|\neq |m|$, so that $u^{(1)}_{n,m}=
w^{(1)}_{n,m}$. By taking into account only the badge
labels of the lines, by item (4) there is only one tree whose root line
is a $w$-line, and it has one node $\vvv_{0}$ (the special vertex
of the tree) with $s_{\vvvv_{0}}=3$,
hence three end-points $\vvv_{1}$, $\vvv_{2}$ and $\vvv_{3}$.
By applying the rules listed above one obtains, for $|n|\neq |m|$,
$$ w^{(1)}_{n,m}= {1\over -\o^2 n^2+\tilde\o^2_{m}}
\sum_{n_1+n_2+n_3=n \atop m_1+m_2+m_3=m} v^{(0)}_{n_1,m_1}
v^{(0)}_{n_2,m_2} v^{(0)}_{n_3,m_3} = 
\sum_{\th \in \Theta^{*(1)}_{n,m} } \Val(\th) ,
\Eq(3.11) $$
where the sum is over all trees $\th$ which can be obtained
from the tree appearing in Figure 3.3 by summing over
all labels which are not explicitly written.

\midinsert
\*
\insertplotttt{76pt}{79pt}{
\ins{30pt}{30pt}{$w$}
}
{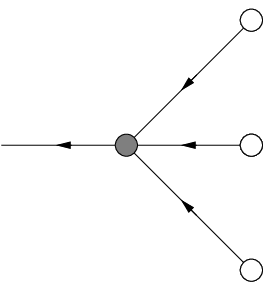}
\*
\line{\vtop{\line{\hskip1.5truecm\vbox{\advance\hsize by -3.1 truecm
\centerline{\nota
\figura(3.3)}
} \hfill} }}
\*
\endinsert

It is easy to realize that \equ(3.11)
corresponds to \equ(2.29) for $k=1$.
Each end-point $\vvv_{i}$ is graphically a white bullet with
first mode labels $(n_{i},m_{i})$ and second mode labels $(0,0)$,
and has associated an end-point factor
$(-1)^{1+\d_{n_{i},m_{i}}} a_{0,n_{i}}$ (see \equ(3.4) in item (7)).
The node $\vvv_{0}$ is represented as a (small) gray bullet,
with mode labels $(0,0)$ and $(0,0)$,
and the factor associated to it is $\h_{\vvvv_{0}}=\e$
(see \equ(3.3) in item (6)).
We associate to the line $\ell$ coming out from the node $\vvv_{0}$
a momentum $(n_{\ell},n_{\ell})$, with $n_{\ell}=n$, and a propagator
$g_{\ell}=1/(-\o^{2}n_{\ell}^{2}+\tilde\o_{m_{\ell}}^{2})$
(see \equ(3.1) in item (3)).

Now we consider the case $|n| =|m|$, so that $u^{(1)}_{n,m}=
\pm A^{(1)}_{n}$ (see \equ(2.14)).
By taking into account only the badge labels of the lines,
there are four trees contributing to $A^{(1)}_{n}$:
they are represented by the four trees in Figure 3.4
(the tree b and c are simply obtained from the tree a
by a different choice of the $w$-line entering the last node).

\midinsert
\*
\insertplotttt{407pt}{121pt}{
\ins{26pt}{55pt}{$v$}
\ins{54pt}{83pt}{$w$}
\ins{140pt}{55pt}{$v$}
\ins{166pt}{55pt}{$w$}
\ins{256pt}{55pt}{$v$}
\ins{286pt}{25pt}{$w$}
\ins{372pt}{55pt}{$v$}
\ins{40pt}{15pt}{a}
\ins{168pt}{15pt}{b}
\ins{275pt}{15pt}{c}
\ins{370pt}{15pt}{d}
}
{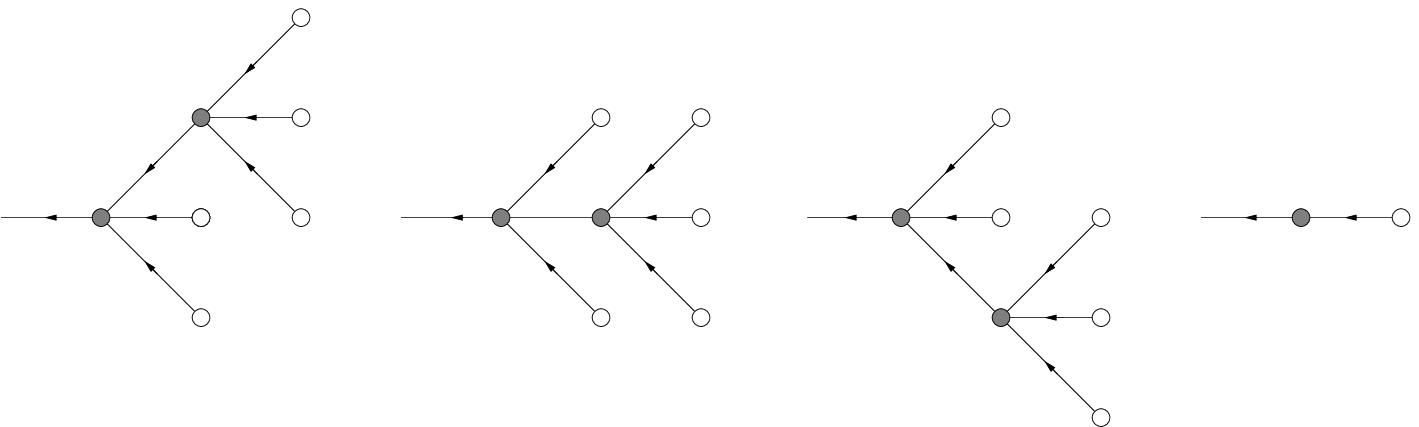}
\*
\line{\vtop{\line{\hskip1.5truecm\vbox{\advance\hsize by -3.1 truecm
\centerline{\nota
\figura(3.4)}
} \hfill} }}
\*
\endinsert

In the trees of Figures 3.4a, 3.4b and 3.4c the root line comes out from
a node $\vvv_{0}$ (the special vertex of the tree) with $s_{\vvvv_{0}}=3$,
and two of the entering lines come out from end-points:
then the other line has to be a $w$-line (by item (4)),
and \equ(3.6) requires that the subtree which has such a line as
root line is exactly the tree represented in Figure 3.2.
In the tree of Figure 4.4d the root line comes out from a node $\vvv_{0}$
with $s_{\vvvv_{0}}=1$, hence the line entering $\vvv_{0}$ is a $v$-line
coming out from an end-point (again see item (4)).

By defining $\Theta^{*(1)}_{n,n}$ as the set of all labeled
trees which can be obtained by assigning to the trees in Figure 3.4
the labels which are not explicitly written, one finds
$$ A^{(1)}_{n} = \sum_{\th \in \Theta^{*(1)}_{n,n} } \Val(\th) ,
\Eq(3.12) $$
which corresponds to the sum of two contributions. The first one
arises from the trees of Figures 3.4a, 3.4b and 3.4c, and it is given by
$$ 3 \sum_{n'\in\zzz} \LLL_{n,n'}
\sum_{\scriptstyle n'_1+n'_2+n'_3=n' \atop
\scriptstyle m'_1+m'_2+m'_3=n'}
v^{(0)}_{n'_1,m'_1} v^{(0)}_{n'_2,m'_2} w^{(1)}_{n'_3,m'_3} ,
\Eq(3.13) $$
where one has
$$ \eqalign{
\LLL_{n,n'} & =
B_{\mmmm} \O_{\mmmm}^{-2} D_{\mmmm}^{2}
{\mathop{\sum}_{\scriptstyle n_{1}+n_{2}=n-n' \atop
\scriptstyle n_{2}=-n'}}^{*}
s_{n_{1}} s_{n_{2}} +
B_{\mmmm} {\mathop{\sum}_{n_{1}+n_{2}=n-n'}}^{*}
{1 \over i^{2} (n_{2}+n')^{2}} cd_{n_{1}} cd_{n_{2}} \cr
& \qquad \qquad + 
B_{\mmmm} \O_{\mmmm}^{-1} D_{\mmmm}
{\mathop{\sum}_{n_{1}+n_{2}=n-n'}}^{*}
{1 \over i (n_{2}+n')} s_{n_{1}} cd_{n_{2}} \cr
& \qquad \qquad -
B_{\mmmm} \O_{\mmmm}^{-1} D_{\mmmm}
{\mathop{\sum}_{n_{1}+n_{2}=n-n'}}^{*}
{1 \over i (n_{2}+n')} cd_{n_{1}} s_{n_{2}} , \cr}
\Eq(3.14) $$
with the $*$ denoting the constraint $n_{2}+n'\neq 0$.
The first and second mode labels associated to the node $\vvv_{0}$ are,
respectively, $(m_{\vvvv_{0}}',n_{\vvvv_{0}}')=(n_{1},n_{1})$
and $(m_{\vvvv_{0}},n_{\vvvv_{0}})=(n_{2},n_{2})$, while
the momentum flowing through the root line is given
by $(n_{\ell},m_{\ell})$, with $|m_{\ell}|=|n_{\ell}|$ expressed
according to the definition \equ(3.5) in item (8):
the corresponding propagator is $(n_{\ell})^{\d_{\ell}}$ for
$n_{\ell}\neq 0$ and $1$ for $n_{\ell}=0$, as in \equ(3.1) in item (3).

The second contribution corresponds to the tree of Figure 3.4d,
and it is given by
$$ \sum_{n'\in\zzz} \LLL_{n,n'} C^{(1)} a_{0,n'} ,
\Eq(3.15) $$
with the same expression \equ(3.14) for $\LLL_{n,n'}$ and
$C^{(1)}$ still undetermined.
The mode labels of the node $\vvv_{0}$ and the momentum
of the root line are as before.

Then one immediately realize that the sum of \equ(3.13)
and \equ(3.15) corresponds to \equ(2.25) for $k=1$.

Finally that $C^{(1)}$ is given by \equ(3.9)
follows from \equ(2.10). This completes the check of the case $k=1$.

In general from \equ(2.29) one gets,
for $\th\in \Theta^{*(k)}_{n,m}$ contributing
to $w^{(k)}_{n,m}$, that the tree value $\Val(\th)$
is obtained by summing all contribution either of the form
$$ \eqalign{
& {1\over -\o^2 n^2 + \tilde\o_{m}^2} \, \e
\sum_{\scriptstyle n_{1}+n_{2}+n_{3}=n \atop
\scriptstyle m_{1}+m_{2}+m_{3}=m}
\sum_{k_{1}+k_{2}+k_{3}=k-1} \cr& \qquad 
\sum_{\th_{1} \in \Theta^{*(k_{1})}_{n_{1},m_{1}}}
\sum_{\th_{2} \in \Theta^{*(k_{2})}_{n_{2},m_{2}}}
\sum_{\th_{3} \in \Theta^{*(k_{3})}_{n_{3},m_{3}}}
\Val(\th_1)\,\Val(\th_2) \, \Val(\th_3) , \cr}
\Eq(3.16) $$
or of the form
$$ {1\over -\o^2 n^2 + \tilde\o_{m}^2} \sum_{c=a,b} \n_{m}^{(c)}
\sum_{\th_{1} \in \Theta^{*(k,1)}_{n,m^{(c)}}} \Val(\th_1) ,
\Eq(3.17) $$
with $m^{(a)}=m$ and $m^{(b)}=-m$; the corresponding
graphical representations are as in Figure 3.5.

\midinsert
\*
\insertplotttt{375pt}{183pt}{
\ins{25pt}{34pt}{$w$}
\ins{25pt}{135pt}{$w$}
\ins{126pt}{34pt}{$w$}
\ins{126pt}{135pt}{$w$}
\ins{228pt}{34pt}{$w$}
\ins{228pt}{135pt}{$w$}
\ins{328pt}{34pt}{$w$}
\ins{328pt}{135pt}{$w$}
\ins{30pt}{110pt}{a}
\ins{130pt}{110pt}{b}
\ins{230pt}{110pt}{c}
\ins{330pt}{110pt}{d}
\ins{30pt}{10pt}{e}
\ins{130pt}{10pt}{f}
\ins{230pt}{10pt}{g}
\ins{330pt}{10pt}{h}
}
{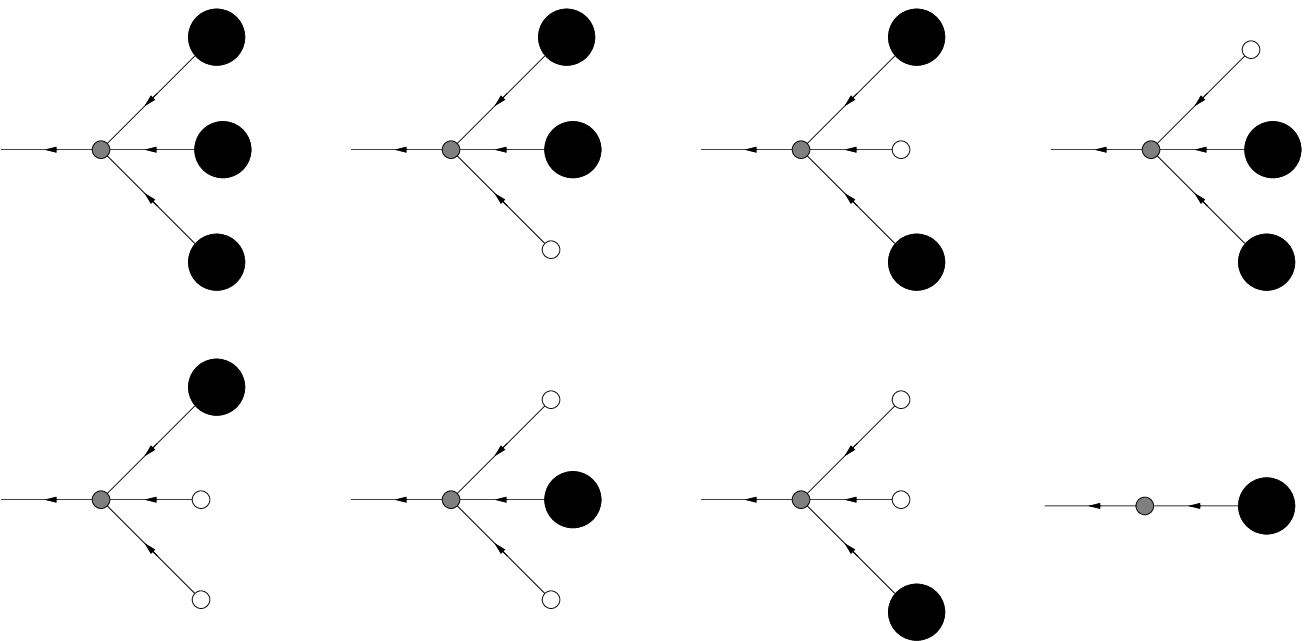}
\*
\line{\vtop{\line{\hskip1.5truecm\vbox{\advance\hsize by -3.1 truecm
\centerline{\nota
\figura(3.5)}} \hfill} }}
\*
\endinsert

Therefore, by simply applying the diagrammatic rules given above,
we see that by summing together the contribution \equ(3.16)
and \equ(3.17) we obtain \equ(3.9) for $|n| \neq |m|$.

A similar discussion applies to $A^{(k)}_{n}$,
and one finds that $A^{(k)}_{n}$ can be written as
sum of contribution either of the form
$$ \eqalign{
& \sum_{n'\in\zzz} \LLL_{n,n'}
\sum_{\scriptstyle n'_1+n'_2+n'_3=n' \atop
\scriptstyle m'_1+m'_2+m'_3=n'}
\sum_{k_{1}+k_{2}+k_{3}=k} \cr & \qquad 
\sum_{\th_{1} \in \Theta^{*(k_{1})}_{n_{1},m_{1}}}
\sum_{\th_{2} \in \Theta^{*(k_{2})}_{n_{2},m_{2}}}
\sum_{\th_{3} \in \Theta^{*(k_{3})}_{n_{3},m_{3}}}
\Val(\th_1)\,\Val(\th_2) \, \Val(\th_3) , \cr}
\Eq(3.18) $$
or of the form
$$ \sum_{n'\in\zzz} \LLL_{n,n'} C^{(k)} a_{0,n'} ,
\Eq(3.19) $$
with $C^{(k)}$ still undetermined. Both
\equ(3.18) and \equ(3.19) are of the form $\Val(\th)$,
for $\th\in\Theta^{*(k)}_{n,m}$.
A graphical representation is in Figure 3.6.

\midinsert
\*
\insertplotttt{429pt}{74pt}{
\ins{23pt}{30pt}{$v$}
\ins{112pt}{30pt}{$v$}
\ins{203pt}{30pt}{$v$}
\ins{294pt}{30pt}{$v$}
\ins{386pt}{30pt}{$v$}
\ins{25pt}{10pt}{a}
\ins{118pt}{10pt}{b}
\ins{209pt}{10pt}{c}
\ins{300pt}{10pt}{d}
\ins{390pt}{10pt}{e}
}
{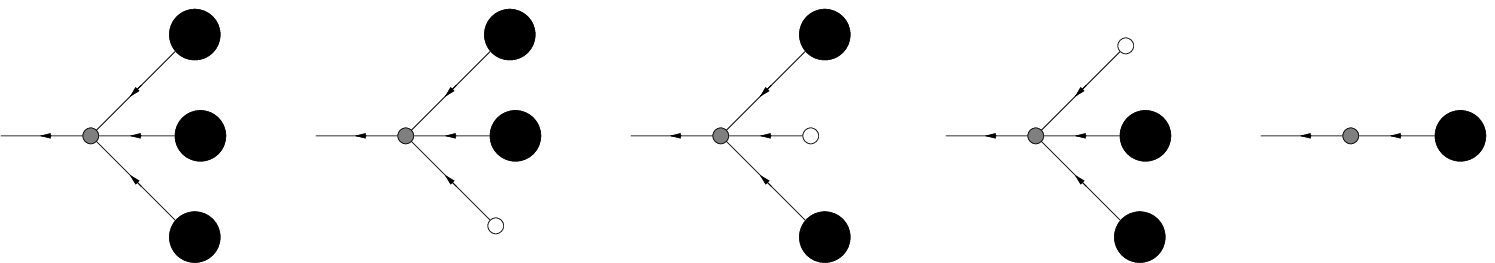}
\*
\line{\vtop{\line{\hskip1.5truecm\vbox{\advance\hsize by -3.1 truecm
\centerline{\nota
\figura(3.6)}} \hfill} }}
\*
\endinsert

Analogously to the case $k=1$ the coefficients $C^{(k)}$
are found to be expressed by \equ(3.10).
Then the lemma is proved. \qed

\*

\0{\bf Lemma 3.}
{\it For any rooted tree $\th$ one has
$|V_{v}^{3}(\th)| \le 2|V_{w}^{3}(\th)|+2|V_{v}^{1}(\th)|$
and $|E(\th)|\le 2(|V_{v}^{3}(\th)|+|V_{w}^{3}(\th)|)+1$.}

\*

\0{\it Proof.} First of all note that $|V_{w}^{3}(\th)|=0$
requires $|V_{v}^{1}(\th)| \ge 1$, so that one has
$|V_{w}^{3}(\th)|+|V_{v}^{1}(\th)|\ge 1$ for all trees $\th$.

We prove by induction on the number $N$ of nodes the bound
$$ \left| V_{v}^{3}(\th) \right| \le \cases{ 
2|V_{w}^{3}(\th)|+2|V_{v}^{1}(\th)|- 1, &
if the root line is a $v$-line , \cr
2|V_{w}^{3}(\th)|+2|V_{v}^{1}(\th)| - 2 &
if the root line is a $w$-line , \cr}
\Eq(3.20) $$
which will immediately imply the first assertion.

For $N=1$ the bound is trivially satisfied, as Figures 3.3 and 3.4 show.

Then assume that \equ(3.20) holds for the trees with $N'$ nodes,
for all $N'<N$, and consider a tree $\th$ with $V(\th)=N$.

If the special vertex $\vvv_{0}$ of $\th$ is not in $V_{v}^{3}(\th)$
(hence it is in $V_{w}(\th)$) the bound \equ(3.20) follows
trivially by the inductive hypothesis.

If $\vvv_{0}\in V_{v}^{3}(\th)$ then we can write
$$ |V_{v}^{3}(\th)| = 1 + \sum_{i=1}^{s} |V_{v}^{3}(\th_{i})| ,
\Eq(3.21) $$
where $\th_{1},\ldots,\th_{s}$ are the subtrees whose
root line is one of the lines entering $\vvv_{0}$.
One must have $s\ge1$, as $s=0$ would correspond to have
all the entering lines of $\vvv_{0}$ coming out from end-points,
hence to have $N=1$.

If $s\ge 2$ one has from \equ(3.21) and from the inductive hypothesis
$$ |V_{v}^{3}(\th)| \le 1 + \sum_{i=1}^{s}
\left( 2 |V_{w}^{3}(\th_{i})| + 2 |V_{v}^{1}(\th_{i})| - 1 \right)
\le 1 + 2|V_{w}^{3}(\th)| + 2|V_{v}^{1}(\th)| - 2 ,
\Eq(3.22) $$
and the bound \equ(3.20) follows.

If $s=1$ then the root line of $\th_{1}$ has to be a $w$-line
by item (4), so that one has
$$ |V_{v}^{3}(\th)| \le 1 + \left( 2 |V_{w}^{3}(\th_{1})|
+ 2 |V_{v}^{1}(\th)| - 2 \right)
\Eq(3.23) $$
which again yields \equ(3.20).

Finally the second assertion follows from the
standard (trivial) property of trees
$$ \sum_{\vvvv\in V(\th)} \left( s_{\vvvv}-1 \right) = |E(\th)| - 1 ,
\Eq(3.24) $$
and the observation that in our case one has $s_{\vvvv}\le 3$. \qed

\*\*
\section(4,Tree expansion: the multiscale decomposition)

\0We assume the Diophantine conditions \equ(2.31). We introduce a
multiscale decomposition of the propagators of the $w$-lines.
Let $\chi(x)$ be a $C^\io$ non-increasing function such that
$\chi(x)=0$ if $|x|\ge C_{0}$ and $\chi(x)=1$ if $|x|\le 2C_{0}$
($C_0$ is the same constant appearing in \equ(2.31)),
and let $\chi_h(x)=\chi(2^{h}x)-\chi(2^{h+1}x)$ for
$h \ge 0$, and $\chi_{-1}(x)=1-\chi(x)$; such functions
realize a smooth partition of the unity as
$$ 1=\chi_{-1}(x)+\sum_{h=0}^{\io} \chi_h(x)
= \sum_{h=-1}^{\io} \chi_h(x) .
\Eq(4.1) $$
If $\chi_{h}(x) \neq 0$ for $h\ge 0$ one has
$2^{-h-1}C_{0}\le|x|\le2^{-h+1}C_{0}$, while if $\chi_{-1}(x)\neq 0$
one has $|x| \ge C_{0}$.

We write the propagator of any $w$-line
as sum of propagators on single scales in the following way:
$$ g (\o n,m) =
\sum_{h=-1}^{\io}
{ \chi_h( | \o n|-\tilde\o_{m})
\over - \o^{2} n^2+\tilde\o_{m}^2 }
= \sum_{h=-1}^{\io} g^{(h)}(\o n,m) .
\Eq(4.2)$$
Note that we can bound $|g^{(h)}(\o n,m)|\le 2^{-h+1}C_{0}$.

This means that we can attach to each $w$-line $\ell$ in $L(\th)$
a scale label $h_{\ell}\ge -1$, which is the scale of the propagator
which is associated to $\ell$. We can denote with
$\Theta^{(k)}_{n,m}$ the set of trees
which differ from the previous ones simply because
the lines carry also the scale labels.
The set $\Theta^{(k)}_{n,m}$ is defined according to the
rules (1) to (8) of Section \secc(3), by changing item (3) into
the following one.

\*

\0(3$^{\prime}$) To each line $\ell$ coming out from nodes of $w$-type
we associate a {\it scale label} $h_{\ell}\ge -1$.
For notational convenience we associate a scale label
$h=-1$ to the lines coming out from the nodes
of $v$-type and to the lines coming out from the end-points.
To each line $\ell$ we associate a {\it propagator}
$$ g^{(h_{\ell})}_{\ell} \= g^{(h_{\ell})} (\o n_{\ell},m_{\ell}) =
\cases{
\displaystyle
{ \chi_{h_{\ell}} ( | \o n_{\ell}|-\tilde\o_{m_{\ell}})
\over - \o^{2} n_{\ell}^2 +\tilde\o_{m_{\ell}}^2} , &
if $\g_{\ell}=w$, \cr
\displaystyle {1 \over (i n_{\ell})^{\d_{\ell}} } , &
if $\g_{\ell}=v$, $n_{\ell}\neq 0$ , \cr
& \cr
\displaystyle 1 , &
if $\g_{\ell}=v$, $n_{\ell} = 0$ , \cr}
\Eq(4.3) $$
with momentum $(n_{\ell},m_{\ell})$.

\*

\0{\bf Definition 4.}
{\it For all $\th\in \Theta^{(k)}_{n,m}$, we define
$$ \Val(\th) = 
\Big( \prod_{\ell \in L(\th)} g_{\ell}^{(h_{\ell})} \Big)
\Big( \prod_{\vvvv \in V(\th)} \h_{{\vvvv}} \Big)
\Big( \prod_{\vvvv \in E(\th) } V_{\vvvv} \Big) ,
\Eq(4.4) $$
the {\rm value of the tree} $\th$.}

\*

Then \equ(3.9) and \equ(3.10) are replaced, respectively, with
$$ u^{(k)}_{n,m} =
\sum_{\th\in\Th^{(k)}_{ n,m}} \Val(\th) ,
\Eq(4.5) $$
and
$$ C^{(k)} = r_{0}^{-1}
{\mathop{\sum}_{\th\in \Th^{(k)}_{n,n}}}^{*}
a_{0,-n} \Val(\th) ,
\Eq(4.6) $$
with the new definition for the tree value $\Val(\th)$
and with $*$ meaning the same constraint as in \equ(3.10).

\*

\0{\bf Definition 5.} {\it
A {\rm cluster} $T$ is a connected set of nodes which are linked
by a continuous path of lines with the same scale label $h_{T}$ or
a lower one and which are maximal; we shall say that the cluster
has scale $h_{T}$. We shall denote with $V(T)$ and $E(T)$
the set of nodes and the set of end-points, respectively, which are
contained inside the cluster $T$, and with $L(T)$
the set of lines connecting them. As for trees
we call $V_{v}(T)$ and $V_{w}(T)$ the sets
of nodes $\vvv\in V(T)$ which are of $v$-type and
of $w$-type respectively.}

\*

We define the {\rm order} $k_{T}$ of a cluster $T$
as the order of a tree (see item (ii) before Definition 3),
with the sums restricted to the nodes internal to the cluster.

An inclusion relation is established between clusters,
in such a way that the innermost clusters are the clusters
with lowest scale, and so on. Each cluster $T$ can have an arbitrary
number of lines entering it (incoming lines), but only one or zero line
coming from it (outcoming line); we shall denote the latter
(when it exists) with $\ell_{T}^{1}$.
We shall call {\it external lines} of the cluster $T$
the lines which either enter or come out from $T$, and
we shall denote by $h_{T}^{(e)}$
the minimum among the scales of the external lines of $T$.
Define also
$$ \eqalign{
K(\th) & = \sum_{\vvvv \in V(\th)}
\left( |n_{\vvvv}'|+|n_{\vvvv}| \right) +
\sum_{\vvvv \in E(\th) } |n_{\vvvv}'| , \cr
K(T) & = \sum_{\vvvv \in V(T)}
\left( |n_{\vvvv}'|+|n_{\vvvv}| \right) +
\sum_{\vvvv \in E(\th) } |n_{\vvvv}'| , \cr}
\Eq(4.7) $$
where we recall that one has $(n_{\vvvv}',m_{\vvvv}')=
(n_{\vvvv},m_{\vvvv})=0$ if $\vvv\in V(\th)$ is of $w$-type.

\*

If a cluster has only one entering line $\ell_{T}^{2}$
and $(n,m)$ is the momentum of such a line,
for any line $\ell\in L(T)$ one can write
$ (n_{\ell},m_{\ell}) = (n^{0}_{\ell},m^{0}_{\ell}) + \h_{\ell}(n,m)$,
where $\h_{\ell}=1$ if the line $\ell$ is along the path connecting
the external lines of $T$ and $\h_{\ell}=0$ otherwise.

\*

\0{\bf Definition 6.} {\it A cluster $T$ with only one incoming
line $\ell_{T}^{2}$ such that one has
$$  n_{\ell_{T}^{1}} =  n_{\ell_{T}^{2}} \quad \hbox{ and } \quad
m_{\ell_{T}^{1}}= \pm m_{\ell_{T}^{2}}
\Eq(4.8) $$
will be called a {\rm self-energy graph} or a {\rm resonance}.
In such a case we shall call a {\rm resonant line}
the line $\ell_{T}^{1}$, and we shall refer to its momentum
as the momentum of the self-energy graph.}

\*

Examples of self-energy graphs $T$ with $k_{T}=1$ are represented
in Figure 4.1. The lines crossing the encircling bubbles
are the external lines, and in the figure are assumed to be on
scales higher than the lines internal to the bubbles. There
are 9 self-energy graphs with $k_{T}=1$: they are all obtained
by the three drawn in Figure 4.1, simply by considering
all possible unequivalent trees.

\midinsert
\*
\insertplotttt{221pt}{113pt}{
\ins{24pt}{64pt}{$w$}
\ins{140pt}{64pt}{$w$}
\ins{168pt}{36pt}{$v$}
\ins{24pt}{12pt}{a}
\ins{142pt}{12pt}{b}
}
{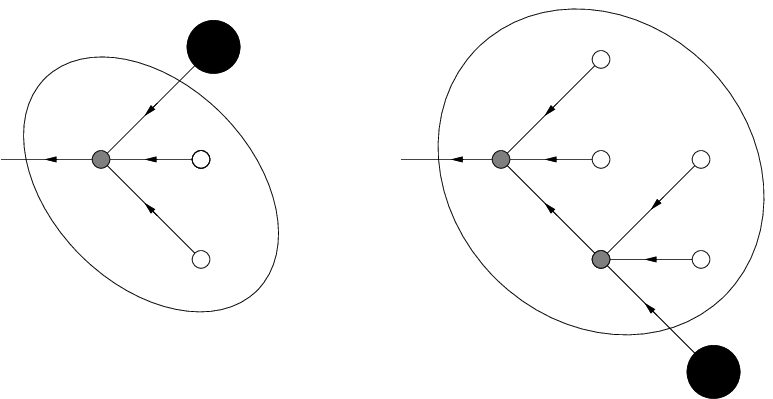}
\*
\line{\vtop{\line{\hskip1.5truecm\vbox{\advance\hsize by -3.1 truecm
\centerline{\nota
\figura(4.1)}} \hfill} }}
\*
\endinsert

\*

\0{\bf Definition 7.}
{\it The {\rm value of the self-energy graph $T$} with momentum $(n,m)$
associated to the line $\ell_{T}^{2}$ is defined as
$$\VV_{T}^{h}(\o n,m)=
\Big( \prod_{\ell \in T} g^{(h_{\ell})}_{\ell} \Big)
\Big( \prod_{\vvv \in V(T)} \h_{\vvvv} \Big)
\Big( \prod_{\vvv \in E(T) } V_{\vvvv} \Big) .
\Eq(4.9) $$
where $h=h^{(e)}_T$ is the minimum between the scales of the
two external lines of $T$ (they can differ at most by a unit), 
and one has
$$ \eqalign{
n(T) & \= \sum_{\vvvv \in V(T)} \left( n_{\vvvv}' + n_{\vvvv} \right) +
\sum_{\vvvv\in E(T) } n_{\vvvv}' = 0 , \cr
m(T) & \= \sum_{\vvvv \in V(T)} \left( m_{\vvvv}' + m_{\vvvv} \right) +
\sum_{\vvvv\in E(T) } m_{\vvvv}' \in \{0,2m\} , \cr}
\Eq(4.10) $$
by definition of self-energy graph; one says that $T$ is
a resonance of type $c=a$ when $m(T)=0$ and a resonance
of type $c=b$ when $m(T)=2m$.}

\*

\0{\bf Definition 8.}
{\it Given a tree $\th$, we shall denote by
$N_{h}(\th)$ the number of lines with scale $h$, and by
$C_{h}(\th)$ the number of clusters with scale $h$.}

\*

Then the product of propagators appearing in \equ(4.4)
can be bounded as
$$ \Big| \prod_{\ell \in L(\th)} g_{\ell}^{(h_{\ell})} \Big| \le
\Big( \prod_{h=0}^{\io} 2^{h N_{h}(\th)} \Big)
\Big( \prod_{\scriptstyle \ell\in L(\th) \atop \scriptstyle \g_{\ell} =w}
{\chi_{h_{\ell}}(|\o n_{\ell}|-
\tilde \o_{m_{\ell}}) \over |\o n_{\ell}|+\tilde\o_{m_{\ell}} } \Big)
\Big( \prod_{\scriptstyle \ell\in L(\th) \atop
\scriptstyle \g_{\ell} =v, \, n_{\ell}\neq 0}
{1 \over |\o n_{\ell}|} \Big) ,
\Eq(4.11) $$
and this will be used later.

\*

\0{\bf Lemma 4.} {\it Assume $0<C_{0}<1/2$ and that there is a
constant $C_{1}$ such that one has $| \tilde\o_{m}-|m| |\le C_{1}\e/|m|$.
If $\e$ is small enough for any tree $\th \in \Th_{n,m}^{(k)}$
and for any line $\ell$ on a scale $h_{\ell} \ge 0$ one has 
$\min\{m_{\ell},n_{\ell}\} \ge 1/2\e$.}

\*

\0{\it Proof.} If a line $\ell$ with momentum $(n,m)$ is on scale
$h \ge 0$ then one has
$$ \eqalign{
{1 \over 2} > C_{0} \ge \left| |\o n| - \tilde\o_{m} \right| & \ge
\left| \left( \sqrt{1-\e} - 1 \right) |n| + \left( |n| - |m| \right)
- C_{1} \e/|m| \right| \cr
& \ge \left| \left| {\e \, |n|  \over 1 + \sqrt{1-\e}} -
\left(|n|-|m|\right) \right| - C_{1}\e/|m| \right| , \cr}
\Eq(4.12) $$
with $|n|\neq|m|$, hence $|n-m|\ge 1$, so that $|n|\ge 1/2\e$.
Moreover one has $||\o n|-\tilde \o_{m}| \le 1/2$ and
$\tilde\o_{m}-|m|=O(\e)$, and one obtains also $|m|>1/2\e$. \qed

\*

\0{\bf Lemma 5.} {\it
Define $h_{0}$ such that $2^{h_{0}} \le 16 C_{0}/\e < 2^{h_{0}+1}$,
and assume that there is a constant $C_{1}$ such that
one has $| \tilde\o_{m}-|m| |\le C_{1}\e/|m|$.
If $\e$ is small enough for any tree $\th \in \Th_{n,m}^{(k)}$
and for all $h\ge h_{0}$ one has
$$ N_h(\th) \le  4 K(\th) 2^{(2-h)/\t} -
C_{h}(\th) + S_{h}(\th) + M_{h}^{\n}(\th) , 
\Eq(4.13) $$
where $K(\th)$ is defined in \equ(4.7),
while $S_{h}(\th)$ is the number of self-energy graphs $T$
in $\th$ with $h_{T}^{(e)}=h$
and $M_{h}^{\n}(\th)$ is the number of $\n$-vertices in $\th$
such that the maximum scale of the two external lines is $h$.}

\*

\0{\it Proof.}
We prove inductively the bound
$$ N^{*}_{h}(\th) \le \max\{ 0 , 2 K(\th) 2^{(2-h)/\t} - 1 \} ,
\Eq(4.14) $$
where $N_{h}^{*}(\th)$ is the number of non-resonant lines in $L(\th)$
on scale $h'\ge h$.

First of all note that for a tree $\th$ to have a line on scale
$h$ the condition $K(\th) > 2^{(h-1)/\t}$ is necessary,
by the first Diophantine conditions in \equ(2.31).
This means that one can have $N_{h}^{*}(\th) \ge 1$ only
if $K=K(\th)$ is such that $K > k_{0}\=2^{(h-1)/\t}$:
therefore for values $K \le k_{0}$ the bound \equ(4.14) is satisfied.

If $K=K(\th)>k_{0}$, we assume that
the bound holds for all trees $\th'$ with $K(\th')<K$.
Define $E_{h}=2^{-1}(2^{(2-h)/\t})^{-1}$: so we have
to prove that $N^{*}_{h}(\th)\le \max\{0,K(\th) E_{h}^{-1}-1\}$.

Call $\ell$ the root line of $\th$ and
$\ell_{1},\ldots,\ell_{m}$ the $m\ge 0$ lines on scale $\ge h$
which are the closest to $\ell$ (i.e. such that no other line
along the paths connecting the lines $\ell_{1},\ldots,\ell_{m}$
to the root line is on scale $\ge h$).

If the root line $\ell$ of $\th$ is either on scale $< h$
or on scale $\ge h$ and resonant, then
$$ N_{h}^{*}(\th) = \sum_{i=1}^{m} N_{h}^{*}(\th_{i}) ,
\Eq(4.15) $$
where $\th_{i}$ is the subtree with $\ell_{i}$ as root line,
hence the bound follows by the inductive hypothesis.

If the root line $\ell$ has scale $\ge h$ and is non-resonant,
then $\ell_{1},\ldots,\ell_{m}$ are the entering line of a cluster $T$.

By denoting again with $\th_{i}$ the subtree
having $\ell_{i}$ as root line, one has
$$ N_{h}^{*}(\th) = 1 + \sum_{i=1}^{m} N_{h}^{*}(\th_{i}) ,
\Eq(4.16) $$
so that the bound becomes trivial if either $m=0$ or $m\ge 2$.

If $m=1$ then one has a cluster $T$ with two external lines
$\ell$ and $\ell_{1}$, which are both with scales $\ge h$; then
$$ \left| |\o n_{\ell}| - \tilde\o_{m_{\ell}} \right|
\le 2^{-h+1} C_{0} , \qquad
\left| |\o n_{\ell_{1}}| - \tilde\o_{m_{\ell_{1}}} \right|
\le 2^{-h+1} C_{0} ,
\Eq(4.17) $$
and recall that $T$ is not a self-energy graph.

Note that the validity of both inequalities in \equ(4.17)
for $h \ge h_{0}$ imply that one has $|n_{\ell}-n_{\ell_{1}}|\neq
|m_{\ell}\pm m_{\ell_{1}}|$, as we are going to show.

By Lemma 4 we know that one has $\min\{m_{\ell},n_{\ell}\} \ge 1/2\e$.
Then from \equ(4.17) we have,
for some $\h_{\ell},\h_{\ell_{1}}\in\{\pm 1\}$,
$$ 2^{-h+2} C_{0} \ge \big| \o
( n_{\ell}-n_{\ell_{1}}) + \h_{\ell} \tilde\o_{m_{\ell}}+
\h_{\ell_{1}} \tilde\o_{m_{\ell_{1}}} \big|,
\Eq(4.18) $$
so that if one had $|n_{\ell}-n_{\ell_{1}}|
= |m_{\ell}\pm m_{\ell_{1}}|$ we would obtain for $\e$ small enough
$$ 2^{-h+2} C_{0}  \ge {\e \over 1 + \sqrt{1-\e}} \left|
n_{\ell} - n_{\ell_{1}} \right| - {C_{1} \e \over |m_{\ell}|} -
{C_{1} \e \over |m_{\ell_{1}}|} \ge {\e \over 2} - 4 C_{1} \e^{2} >
{\e \over 4} ,
\Eq(4.19) $$
which is contradictory with $h \le h_{0}$; hence one has
$|n_{\ell}-n_{\ell_{1}}| \neq |m_{\ell}\pm m_{\ell_{1}}|$.

Then, by \equ(4.17) and for $|n_{\ell}-n_{\ell_{1}}|
\neq |m_{\ell}\pm m_{\ell_{1}}|$,
one has, for suitable $\h_{\ell},\h_{\ell_{1}}\in\{+,-\}$,
$$ 2^{-h+2} C_{0} \ge \big| \o
( n_{\ell}-n_{\ell_{1}}) + \h_{\ell} \tilde\o_{m_{\ell}}+
\h_{\ell_{1}} \tilde\o_{m_{\ell_{1}}} \big| \ge
C_{0} |n_{\ell}-n_{\ell_{1}}|^{-\t} ,
\Eq(4.20) $$
where the second Diophantine conditions in \equ(2.31) have been used.
Hence $K(\th)-K(\th_{1})> E_{h}$,
which, inserted into \equ(4.16) with $m=1$, gives,
by using the inductive hypothesis,
$$ \eqalign{
N_{h}^{*}(\th) & = 1 + N_{h}^{*}(\th_{1})
\le 1 + K(\th_{1}) E_{h}^{-1} - 1 \cr
& \le 1 + \Big( K(\th) - E_{h}\Big) E_{h}^{-1} - 1
\le K(\th) E_{h}^{-1} - 1 , \cr}
\Eq(4.21) $$
hence the bound is proved also
if the root line is on scale $\ge h$.

In the same way one proves that, if we denote with $C_{h}(\th)$
the number of clusters on scale $h$, one has
$$ C_{h}(\th) \le \max\{ 0 , 2K(\th) 2^{(2-h)/\t} - 1 \} ;
\Eq(4.22) $$
see \cita{GM} for details. \qed

\*

Note that the argument above is very close to \cita{GM}:
this is due to the fact that the external lines of any self-energy graph
$T$ are both $w$-lines, so that the only effect of the presence
of the $v$-lines and of the nodes of $v$-type is in the
contribution to $K(T)$.

\*

The following lemma deals with the lines on scale $h <h_{0}$.

\*

\0{\bf Lemma 6.} {\it Let $h_{0}$ be defined as in Lemma 2
and $C_{0}<1/2$, and assume that there is a constant $C_{1}$
such that one has $|\tilde\o_{m}-|m| |\le C_{1}\e$.
If $\e$ is small enough for $h<h_{0}$ one has $|g^{(h)}_{\ell}|\le 32$.}

\*

\0{\it Proof.}
Either if $h\neq h_{\ell}$ of $h=h_{\ell}=-1$ the bound is trivial.
If $h=h_{\ell} \ge 0$ one has
$$ g^{(h)}_{\ell} = { \chi_{h}(|\o n_{\ell}|-\tilde\o_{m_{\ell}})
\over -|\o n_{\ell}|+\tilde\o_{m} } \,
{1 \over |\o n_{\ell}|+\tilde\o_{m} } ,
\Eq(4.23) $$
where $|\o n_{\ell}|+\tilde\o_{m} \ge 1/2\e$ by Lemma 4.
Then one has
$$ { 1 \over |\o n_{\ell}|+\tilde\o_{m} }\le 2 \e ,
\Eq(4.24) $$
which, inserted in \equ(4.23), gives $|g^{(h)}_{\ell}| \le
2^{h+2} \e/C_{0} \le 32$, so that the lemma is proved. \qed

\*\*
\section(5,The renormalized expansion)

\0It is an immediate consequence of Lemma 5 and Lemma 6
that all the trees $\th$ with no self-energy graphs or
$\n$-vertices admit a bound $O(C^k\e^k)$,
where $C$ is a constant. However the generic tree
$\th$ with 
$S_h(\th)\not=0$ admits a much worse bound, namely $O(C^k\e^k k!^\a)$,
for some constant $\a$, and the presence of factorials
prevent us to prove the convergence of the series; 
in KAM theory this is called {\it accumulation
of small divisors}.
It is convenient then to consider another expansion
for $u_{n,m}$, which is essentially a resummation
of the one introduced in Sections \secc(3) and \secc(4).

We define the set $\Th^{(k)\RR}_{n,m}$ 
of {\it renormalized trees}, which are defined
as $\Th^{(k)}_{n,m}$  
except that the following rules are added.

\*

\0(9) To each self-energy graph (with $|m|\ge 1$) the 
$\RR=\openone-\LL$ operation is applied, where $\LL$ acts  on
the self-energy graphs in the following way,
for $h\ge 0$ and $|m|\ge 1$,
$$ \LL \VV_T^{h}(\o n,m) = 
\VV_T^{h}({\rm sgn}(n)\,\tilde\o_{m},m) ,
\Eq(5.1) $$
$\RR$ is called {\it regularization operator}; its action
simply means that each self-energy graph  
$\VV_T^{h}(\o n,m)$ must be replaced by $\RR
\VV_T^{h}(\o n,m)$.

\*

\0(10) To the nodes $\vvv$ of $w$-type with $s_{\vvvv}=1$
(which we still call $\n$-vertices) and with
$h\ge 0$ the minimal scale among the lines entering or exiting $\vvv$,
we associate a factor $2^{-h} \n^{(c)}_{h,m}$, $c=a,b$, 
where $(n,m)$ and $(n,\pm m)$, with $|m|\ge 1$, 
are the momenta of the lines
and $a$ corresponds to the sign $+$ and $b$ to the sign $-$ in $\pm m$.

\*

\0(11) The set $\{h_\ell\}$ of the scales  
associated to the lines $\ell\in L(\th)$ must
satisfy the following constraint
(which we call {\it compatibility}):
fixed $(n_\ell,m_\ell)$ for any  $\ell\in L(\th)$ and replaced
$\RR$ with $\openone$ at each self-energy graph, one must have
$\chi_{h_{\ell}}(|\o n_{\ell}|-\tilde\o_{m_{\ell}}) \neq 0$.

\*

\0(12) The factors $C^{(k_{\vvvv})}$
in \equ(3.2) are replaced with $\G$, to be considered a parameter.

\*

The set $\Th^{(k)\RR}_{n,m}$ is defined as $\Th^{(k)}_{n,m}$ with
the new rules and with the constrain that the order $k$ is given by
$k=|V_{w}(\th)|+|V_{v}^{1}(\th)|$.

\*

We consider the following expansion
$$ \tilde u_{n,m} = \sum_{k=1}^{\io}\m^k
\sum_{\th \in \Th^{(k)\RR}_{n,m}} \Val(\th) ,
\Eq(5.4) $$
where, for $|m|\ge 1$ and $h\ge 0$,  $\n^{(c)}_{h,m}$ is given by
$$ 2^{-h} \m \n^{(c)}_{h,m}=\m \n^{(c)}_{m} + {1 \over 2}
\sum_{\s=\pm} \sum_{T \in \TT^{(c)}_{< h}} \m^{k_{T}}
\VV_T^{h}(\s\tilde \o_{m},m) ,
\Eq(5.3) $$
with $c=a,b$, and $\TT^{(c)}_{< h}$ denoting the set of
self-energy graphs $T$ of type $c$
(see item (6) in Section \secc(3)) with $h_{T}<h$,
and $\G$ is determined by the self-consistence equation
$$ \G = r_{0}^{-1} 
\sum_{k=1}^{\io} \m^{k-1}
{\mathop{\sum}_{\th \in \Th^{(k)\RR}_{n,n}}}^{*}
a_{0,-n} \Val(\th) ,
\Eq(5.3a) $$
with $*$ denoting the same constraint as in \equ(3.10).

We shall set also $\n_{m}^{(c)}=\n_{-1,m}^{(c)}$.
Note that $\VV_T^{h}(\s\tilde \o_{m},m)$
is indipendent from $\s$.

Calling $L_0(\th), V_0(\th), E_0(\th)$
the set of lines, node and end-points 
not contained in any self-energy graph, and $S_{0}(\th)$
the {\it maximal self-energy graphs}, \ie the self-energy
graphs which are not contained in any self-energy graphs,
we can write $\Val(\th)$ in \equ(5.4) as 
$$ \Val(\th) = 
\Big( \prod_{\ell \in L_0(\th)} g_{\ell}^{(h_{\ell})} \Big)
\Big( \prod_{\vvv \in V_0(\th)} \h_{{\vvvv}} \Big)
\Big( \prod_{\vvv \in E_0(\th) } V_{\vvvv}\Big) 
\Big(\prod_{T\in S_0(\th)} \RR \VV^{h_T^{e}}_T(\o n_{\ell_T},m_{\ell_T})
\Big) ,
\Eq(5.4aa) $$
and by definition
$$ \RR \VV^{h_T^{e}}_T(\o n_{\ell_T},m_{\ell_T})=
 \VV^{h_T^{e}}_T(\o n_{\ell_T},m_{\ell_T})-
\VV_T^{h_T^{e}}({\rm sgn}(n_{\ell_T})\,\tilde\o_{m_{\ell_T}},m_{\ell_T}) ,
\Eq(5.1aa) $$
and $\VV^{h_T^{e}}_T(\o n_{\ell_T},m_{\ell_T})$
is given by
$$\VV^{h_T^{e}}_T(\o n_{\ell_T},m_{\ell_T})=
\Big( \prod_{\ell \in L_0(T)} g^{(h_{\ell})}_{\ell} \Big)
\Big( \prod_{\vvv \in V_0(T)} \h_{\vvvv} \Big)
\Big( \prod_{\vvv \in E_0(T) } V_{\vvvv} \Big)
\Big(\prod_{T'\in S_0(T)} \RR \VV^{h_{T'}^{e}}_{T'}
(\o n_{\ell_{T'}},m_{\ell_{T'}}) \Big) .
\Eq(5.9aa) $$

First (Lemma 7) we will show that the expansion \equ(5.4)
is well defined, for $\n_{h,m}, \G=O(\e)$;
then (Lemma 8 and 9) we show that under the same conditions also
the r.h.s. of \equ(5.3) is well defined; moreover (Lemma 10)
we prove by using \equ(5.3) that it is indeed 
possible to choose $\n_{m}^{(c)}$
such that $\n_{h,m}=O(\e)$ for any $h$; then (Lemma 11 and 12)
we show that \equ(5.4) admit a solution $\G=O(\e)$;
finally (Lemma 13) we show that
indeed \equ(5.4) solves the last of \equ(1.13) and \equ(2.2);
this completes the proof of Proposition 1. In the next
section we will solve the implicit function
problem of \equ(2.1) so completing the proof of Theorem 1.

We start from the following lemma stating that, if the 
$\n_{h,m}^{(c)}$ and $\G$ functions are bounded, then the expansion 
\equ(5.4) is well defined.

\*

\0{\bf Lemma 7.} {\it Assume that there exist a constant $C$ such that
one has $|\G|\le C\e $ and $|\n_{h,m}^{(c)}| \le C\e$, with $c=a,b$,
for all $|m|\ge 1$ and all $h\ge 0$.
Then for all $\m_0>0$ there exists $\e_0>0$
such that for all $|\m|\le\m_0$ and
for all $0<\e< \e_{0}$ and for all $(n,m)\in \ZZZ^{2}$ one has
$$ \left| \tilde u_{n,m} \right| \le D_{0} \e \m \,
e^{-\k(|n|+|m|)/4} ,
\Eq(5.6) $$
where $D_{0}$ is a positive constant.
Moreover $u_{n,m}$ is analytic in $\m$ and in 
the parameters $\n_{m',h'}^{(c)}$, with $c=a,b$ and $|m'| \ge 1$.}

\*

\0{\it Proof.} In order to take into account the $\RR$
operation we write \equ(5.1aa) as
$$\RR \VV^{h_T^{e}}_T(\o n_{\ell_T},m_{\ell_T})=
\left(\o n_{\ell_T}-\tilde\o_{m_{\ell_T}} \right)
\int_{0}^{1} {\rm d}t \dpr
\VV^{h_T^{e}}_T(\o n_{\ell_T}+t(\o n_{\ell_T}-\tilde\o_{m_{\ell_T}}) ,
m_{\ell_T}) ,
\Eq(5.1aaa) $$
where $\dpr$ denotes the derivative with respect to the argument
$\o n_{\ell_T}+t(\o n_{\ell_T}-\tilde\o_{m_{\ell_T}})$.

By \equ(5.9aa) we see that the derivatives can be applied
either on the propagators in $L_0(T)$, or 
on the $\RR \VV^{h_{T'}^{e}}_{T'}$.
In the first case there is an extra factor 
$2^{-h^{(e)}_T+h_{T}}$ with respect to the bound
\equ(4.11): $2^{-h^{(e)}_T}$ is obtained from
$\o n_{\ell_T}-\tilde\o_{m_{\ell_T}}$
while $\partial g^{(h_T)}$ is bounded proportionally to $2^{2 h_T}$;
in the second case note that $\dpr_t \RR \VV^{h_{T'}^{e}}_{T'}=
\dpr_t \VV^{h_{T'}^{e}}_{T'}$
as $\LL\VV_{T'}^{h^{(e)}_{T'}}$ is independent of
$t$; if the derivative acts on the propagator
of a line $\ell\in L(T)$, we get a gain factor
$$ 2^{-h^{(e)}_T+h_{T'}} \le 2^{-h^{(e)}_T+h_{T}}
2^{-h^{(e)}_{T'}+h_{T'}} ,
\Eq(5.7) $$
as $h^{(e)}_{T'}\le h_{T}$. We can iterate
this procedure until all the $\RR$
operations are applied on propagators;
at the end (i) the propagators are derived at most one time;
(ii) the number of terms so generated is $\le k$;
(iii) to each self-energy graph $T$ a factor
$2^{-h^{(e)}_T+h_{T}}$ is associated.

Assuming that $|\n^{(c)}_{h,m}|\le C\e$ and $|\G|\le C\e$,
for any $\th$ one obtains,
for suitable constants $\overline{D}$
$$ \eqalign{
& \left| \Val(\th) \right| \le
\e^{|V_w(\th)|+|V_v^{(1)}(\th)|}
\overline{D}^{|V(\th)|} \cr
& \qquad \qquad
\Big( \prod_{h=h_0}^{\io} \exp \Big[ h \log 2 \Big( 4 K(\th)
2^{-(h-2)/\t} - C_{h}(\th) + S_{h}(\th) +
M_{h}^{\n}(\th) \Big) \Big] \Big) \cr
& \qquad \qquad
\Big( \prod_{T \in \SS(\th)\atop h^{(e)}_T \ge h_0 } 
2^{-h^{(e)}_T+h_T} \Big)
\Big( \prod_{h=h_{0}}^{\io} 2^{-h M_{h}^{\n}(\th)} \Big)
\cr & \qquad \qquad
\Big( \prod_{\vvv\in V(\th) \cup E(\th)}
e^{-\k(|n_{\vvvv}|+|n_{\vvvv}'|)} \Big)
\Big( \prod_{\vvv\in V(\th) \cup E(\th)}
e^{-\k(|m_{\vvvv}|+|m_{\vvvv}'|)} \Big) , \cr}
\Eq(5.8) $$
where the second line 
is a bound for $\prod_{h\ge h_0} 2^{h N_h(\th)}$
and we have used that by item (12) $N_h(\th)$ can be
bounded through Lemma 5), and Lemma 4 has been used
for the lines on scales $h< h_{0}$;
moreover $ \prod_{h=h_{0}}^\io 2^{-h M_\n^h(\th)}$
takes into account the factors $2^{-h}$ arising from
the running coupling constants $\n^{(c)}_{h,m}$
and the action of $\RR$ produces, as discussed above, 
the factor $\prod_{T \in \SS(\th)} 2^{-h^{(e)}_T+h_T}$.
Then one has
$$ \eqalign{
& \Big( \prod_{h=h_0}^{\io} 2^{h S_{h}(\th)} \Big)
\Big( \prod_{T \in \SS(\th) } 2^{-h^{(e)}_T} \Big) = 1 , \cr
& \Big( \prod_{h=h_0}^{\io} 2^{-h C_{h}(\th)} \Big)
\Big(\prod_{T \in \SS(\th)}2^{h_T} \Big)  \le 1 . \cr}
\Eq(5.11) $$
We have to sum the values of all trees,
so we have to worry about the sum of the labels.
Recall that a labeled tree is obtained from an unlabeled tree
by assigning all the labels to the points and the lines:
so the sum over all possible labeled trees can be written
as sum over all unlabeled trees and of labels.
For a fixed unlabeled tree $\th$ with a given number of nodes,
say $N$, we can assign first the mode labels
$\{(n_{\vvvv}',m_{\vvvv}'),
(n_{\vvvv},m_{\vvvv})\}_{v\in V(\th)\cup E(\th)}$,
and we sum over all the other labels, which
gives $4^{|V_{v}(\th)|}$ (for the labels $j_{\vvvv}$)
times $2^{|L(\th)|}$ (for the scale labels):
then all the other labels are uniquely fixed.
Then we can perform the sum over the mode labels by
using the exponential decay arising from the node factors
\equ(3.2) and end-point factors \equ(3.4).
Finally we have to sum over the unlabeled trees, and this
gives a factor $4^{N}$ \cita{HP}. By Lemma 3,
one has $|V(\th)|= |V_w(\th)|+|V_{v}^{(1)}(\th)|+|V_{v}^{(3)}(\th)|\le
3(|V_{w}^{(3)}(\th)|+|V_v^{(1)}(\th)|)$, hence $N \le 3k$,
so that $\sum_{\th\in\Th^{(k)\RR}_{n,m}} |\Val(\th)|\le D^k\e^{k}$,
for some positive constant $D$.

Therefore, for fixed $(n,m)$ one has
$$ \sum_{k=1}^{\io} \sum_{\th\in\Th^{(k)}_{n,m}}\m^k
\left| \Val(\th) \right| \le
D_{0} \m\e \, e^{-\k(|n|+|m|)/4} ,
\Eq(5.13) $$
for some positive constant $D_{0}$,
so that \equ(5.6) is proved. \qed

\*

From \equ(5.3) we know that
the quantities $\n_{h,m}^{(c)}$, for $h\ge0$ and $|m|\ge 1$,
verify the recursive relations
$$ \m \n^{(c)}_{h+1,m} = 2 \m \n^{(c)}_{h,m} + \b^{(c)}_{h,m}
(\tilde\o,\e,\{\n_{h',m'}^{(c')}\}) ,
\Eq(5.14) $$
where, by defining $\TT^{(c)}_h$ as the set of self-energy graphs
in $\TT^{(c)}_{<h+1}$ which are on scale $h$, the {\it beta function}
$$ \b^{(c)}_{h,m} \=
\b^{(c)}_{h,m}(\tilde\o,\e,\{\n_{h',m'}^{(c')}\}) = 2^{h+1}
{1 \over 2} \sum_{\s=\pm} \sum_{T\in\TT_h^{(c)}} \m^{k_{T}}
\VV_T^{h+1}(\s\tilde\o_{m},m) ,
\Eq(5.15) $$
depends only on the scales $h'\le h$.

In order to obtain a bound on the beta function, hence
on the running coupling constants, we need to bound
$\VV_{T}^{h+1}(\pm \tilde\o_{m},m)$ for $T\in\TT_h^{(c)}$.
We define $\wt \Theta^{(k)\RR}_{n,m}$
as the set $\Theta^{(k)\RR}_{n,m}$ 
introduced before,
but by changing item (7) into the following one:

\*

\0(7$^{\prime}$) We divide the set $\wt E(\th)$
of end-points into two sets $E(\th)$ and $E_{0}(\th)$.
To each end-point $\vvv \in E(\th)$ we associate a first mode
label $(n_{\vvvv}',m_{\vvvv}')$, with $|m_{\vvvv}'|=|n_{\vvvv}'|$,
a second mode label $(0,0)$
and an {\it end-point factor} $ V_{\vvvv} =
(-1)^{1+\d_{n_{\vvvv}',m_{\vvvv}'}} a_{0,n_{\vvvv}'}$,
while $E_{0}(\th)$ is either the empty set or
a single end-point $\vvvv_{0}$, and, in the latter case,
to the end-point $\vvv \in E_{0}(\th)$ 
we associate a first mode
label $(\lis\o_{m_{\vvvv}},n_{\vvvv})$,
where $\lis\o_{m_{\vvvv}}=\tilde\o_{m_{\vvvv}}/\o$,
a second mode label $(0,0)$
and an {\it end-point factor} $V_{\vvvv} =1$.

\*

Then we have the following generalization of Lemma 4.

\*

\0{\bf Lemma 8.} {\it
If $\e$ is small enough for any tree $\th \in \wt\Th_{n,m}^{(k)\RR}$
one has
$$ N_{h}(\th) \le 4 K(\th) 2^{(2-h)/\t} -
C_{h}(\th) + S_{h}(\th) + M_{h}^{\n}(\th) ,
\Eq(5.20) $$
where the notations are as in lemma 4.}

\*

\0{\it Proof.}
Lemma 4 holds for $E_0(\th)=0$;
we mimic the proof of Lemma 4 proving that
$$ N^{*}_{h}(\th) \le \max\{ 0 , 2 K(\th) 2^{(2-h)/\t} \} ,
\Eq(5.22) $$
for all trees $\th$ with $E_{0}(\th) \neq \emptyset$,
again by induction on $K(\th)$.

For any line $\ell \in L(\th)$ set $\h_{\ell}=1$ if the line
is along the path connecting $\vvvv_{0}$ to the root and
$\h_{\ell}=0$ otherwise, and write
$$ n_{\ell} = n_{\ell}^{0} + \h_{\ell} \lis\o_{m} ,
\qquad m_{\ell} = m_{\ell}^{0} + \h_{\ell} m ,
\Eq(5.23) $$
which implicitly defines $n_{\ell}^{0}$ and $m_{\ell}^{0}$.

Define $k_{0}=2^{(h-1)/\t}$. One has $N_{h}^{*}(\th)=0$ for
$K(\th)<k_{0}$, because if a line $\bar\ell\in L(\th)$ is indeed on scale
$h$ then $|\o n_{\bar\ell}-\tilde \o_{m_{\bar\ell}}|<C_{0}2^{1-h}$,
so that \equ(5.23) and the Diophantine conditions imply
$$ K(\th) \ge \left| n_{\bar\ell}^{0} \right|
> 2^{(h-1)/\t} \= k_{0} .
\Eq(5.24) $$
Then, for $K \ge k_{0}$, we assume that the bound \equ(5.22)
holds for all $K(\th)=K'<K$,
and we show that it follows also for $K(\th)=K$.

If the root line $\ell$ of $\th$ is either on scale $<h$
or on scale $\ge h$ and resonant,
the bound \equ(5.22) follows immediately from the bound \equ(4.13)
and from the inductive hypothesis.

The same occurs if the root line is on scale $\ge h$ and non-resonant,
and, by calling $\ell_{1},\ldots,\ell_{m}$ the lines on scale
$\ge h$ which are the closest to $\ell$, one has $m \ge 2$:
in fact in such a case at least $m-1$ among the subtrees
$\th_{1},\ldots,\th_{m}$ having $\ell_{1},\ldots,\ell_{m}$,
respectively, as root lines have $E(\th_{i})=\emptyset$, so that
we can write, by \equ(4.13) and by the inductive hypothesis,
$$ N_{h}^{*}(\th) = 1 + \sum_{i=1}^{m} N_{h}^{*}(\th_{i})
\le 1 + E_{h}^{-1} \sum_{i=1}^{m} K(\th_{i}) - (m-1)
\le E_{h}   K(\th) ,
\Eq(5.24a) $$
so that \equ(5.22) follows.

If $m=0$ then $N_{h}^{*}(\th)=1$ and $K(\th) 2^{(2-h)/\t} \ge 1$
because one must have $K(\th) \ge k_{0}$.

So the only non-trivial case is when
one has $m=1$. If this happens $\ell_{1}$ is, by construction,
the root line of a tree $\th_{1}$ such that
$K(\th) =  K(T) +  K(\th_{1})$, where $T$ is the cluster
which has $\ell$ and $\ell_{1}$ as external lines and 
$K(T)$, defined in \equ(4.7), satisfies the bound
$K(T) \ge |n_{\ell_{1}} - n_{\ell}|$.

Moreover, if $E_{0}(\th_{1})\neq \emptyset$, one has
$$ \eqalign{
& \left| |\o n_{\ell}^{0}+
\tilde\o_{m}| - \tilde\o_{m_{\ell}} \right|
\le 2^{-h+1} C_{0} , \cr
& \left| |\o n_{\ell_{1}}^{0} +
\tilde\o_{m}| - \tilde\o_{m_{\ell_{1}}} \right|
\le 2^{-h+1} C_{0} , \cr}
\Eq(5.26a) $$
so that, for suitable $\h_{\ell},\h_{\ell_{1}}\in \{-,+\}$, we obtain
$$ 2^{-h+2} C_{0} \ge \big| \o
( n^{0}_{\ell}- n^{0}_{\ell_{1}}) +\h_{\ell} \tilde\o_{m_{\ell}}+
\h_{\ell_{1}} \tilde\o_{m_{\ell_{1}}} \big|
\ge C_{0} |n^{0}_{\ell}-n^{0}_{\ell_{1}}|^{-\t} \=
C_{0} |n_{\ell}-n_{\ell_{1}}|^{-\t} ,
\Eq(5.27) $$
by the second Diophantine conditions in \equ(2.31),
as the quantities $\tilde \o_{m}$
appearing in \equ(5.26a) cancel out.
Therefore one obtains by the inductive hypothesis
$$ N_{h}^{*}(\th) \le 1 + K(\th_{1}) E_{h}^{-1}
\le 1 + K(\th) E_{h}^{-1} -
K(T) E_{h}^{-1} \le  K(\th) E_{h}^{-1} ,
\Eq(5.28) $$
hence the first bound in \equ(5.22) is proved.

If $E_{0}(\th_{1})= \emptyset$, one has 
$$ N_{h}^{*}(\th) \le 1 +  K(\th_{1}) E_{h}^{-1} - 1
\le 1 + K(\th) E_{h}^{-1} - 1 \le  K(\th) E_{h}^{-1} ,
\Eq(5.28a) $$
and \equ(5.22) follows also in such a case. \qed

\*

The following bound for $\VV_T^{h+1}(\pm\tilde\o_{m},m)$,
$h \ge h_{0}$, can then be obtained.

\*

\0{\bf Lemma 9.} {\it 
Assume that there exist a constant $C$ such that
one has $|\G|\le C\e $ and $|\n_{h,m}^{(c)}| \le C\e$, with $c=a,b$,
for all $|m|\ge 1$ and all $h\ge 0$.
Then if $\e$ is small enough
for all $h\ge 0$ and for all $T\in\TT_{h}^{(c)}$ one has
$$|\VV_T^{h+1}(\pm\tilde\o_{m},m)|\le
B^{|V(T)|} e^{-\k 2^{(h-1)/\t}/4} e^{-\k K(T)/4}
\e^{|V_v^{(1)}(T)|+|V_w(T)|} ,
\Eq(5.17) $$
where $B$ is a constant and $K(T)$ is defined in \equ(4.7).}

\*

\0{\it Proof.} 
By using lemma 7 
we obtain for all $T\in\TT_{h}^{(c)}$ and
assuming $h\ge h_0$ we get the bound 
$$ \eqalignno{
& \left| \VV_{T}^{h+1}(\pm\tilde\o_{m},m) \right|
\le \lis B^{|V(T)|} \e^{|V_v^{(1)}(T)|+|V_w(T)|} \cr
& \qquad \qquad \qquad
\prod_{h'=h_0}^{h} \exp \Big[ 4 K(T) \log 2 h' 2^{(2-h')/\t} -
C_{h'}(T) + S_{h'}(T) + M_{h'}^{\n}(T) \Big] \cr
& \qquad \qquad \qquad
\Big( \prod_{T'\subset T\atop h_T^{(\e)}\ge h_0} 2^{-h^{(e)}_T+h_T} \Big)
\Big( \prod_{h'=h_0}^{h} 2^{-h' M_{h'}^{\n}(T)} \Big) e^{-\k|K(T)|/2} ,
& \eq(5.36a)  \cr} $$
where $\lis B$ is a suitable constant.
If $h< h_0$ the bound trivializes as the r.h.s. reduces simply
to  $C^{|V(T)|} |\e|^{|V_v^{(1)}(T)|+|V_w(T)|} 
e^{-{\k\over 2}|K(T)|}$.
The main difference with 
respect to Lemma 6 is that, given
a self-energy graph $T\in\TT_{h}^{(c)}$, there is at least a
line $\ell\in L(T)$ on scale $h_{\ell}=h$ and with
propagator
$$ {1 \over - \o^{2} (n^{0}_{\ell} + \h_{\ell}
\lis \o_{m} )^{2} +  \tilde \o^{2}_{m^{0}_{\ell} + \h_{\ell} m} }  ,
\Eq(5.37) $$
where $\h_{\ell}=1$ if the line $\ell$ belongs to the path of
lines connecting the entering line (carrying a momentum $(n,m)$)
of $T$ with the line coming out of $T$, and $\h_{\ell}=0$ otherwise.
Then one has by the Mel'nikov conditions
$$ C_{0} |n^{0}_{\ell}|^{-\t}\le \left| \o n^{0}_{\ell}
+ \h_{\ell} \tilde \o_m \pm \tilde \o_{m^{0}_{\ell} +\h_{\ell} m}  \right|
\le C_{0} 2^{-h+1} ,
\Eq(5.38) $$
so that $|n^{0}_{\ell}|\ge 2^{(h-1)/\t}$.
On the other hand one has $|n^{0}_{\ell}| \le K(T)$,
hence $K(T)\ge 2^{(h-1)/\t}$; so we get the bound \equ(5.17).\qed

\*

It is an immediate consequence of the above Lemma
that for all $\m_0>0$ there exists $\e_0>0$
such that for all $|\m|\le\m_0$ and $0<\e<\e_0$ one has
$|\b^{(c)}_{h,m}|\le B_1\e|\m|$, with $B_1$
a suitable constant.

We have then proved convergence {\it assuming that}
the parameters $\n_{h,m}$ and $\G$ are bounded; we have
to show that this is actually the case, if the $\n_m^{(c)}$
in \equ(2.32) are chosen in a proper way.

We start proving that it is possible to choose
$\n^{(c)}=\{\n^{(c)}_{m}\}_{|m| \ge 1}$ such that,
for a suitable positive constant $C$, one has
$|\n^{(c)}_{h,m}|\le C\e$ for all $h\ge 0$ and for all $|m|\ge 1$.

For any sequence $a\=\{a_{m}\}_{|m|\ge 1}$ we introduce the norm
$$ \left\| a \right\|_{\io} =\sup_{|m|\ge 1} |a_{m}| .
\Eq(5.20a) $$
Then we have the following result.

\*

\0{\bf Lemma 10.}
{\it Assume that there exists a constant $C$ such that $|\G|\le C\e$.
Then for all $\m_0>0$ there exists $\e_0>$
such that, for all $|\m|\le\m_0$ and
for all $0< \e <\e_{0}$
there is
a family of intervals $I^{(\bar h)}_{c,m}$, $\bar h\ge 0$, $|m|\ge 1$,
$c=a,b$, such that $I^{(\bar h+1)}_{c,m}\subset I^{(\bar h)}_{c,m}$, 
$|I^{\bar h}_{c,m}|\le 2\e (\sqrt{2})^{-(\bar h+1)}$ and,
if $\n^{(c)}_{m}\in I^{(\bar h)}_{c,m}$, then
$$ \| \n^{(c)}_{h} \|_{\io}\le D \e, \qquad \bar h\ge h\ge 0 ,
\Eq(5.21) $$
for some positive constant $D$. Finally one has
$\n_{h,-m}^{(c)}=\n_{h,m}^{(c)}$, $c=a,b$,
for all $\bar h \ge h\ge 0$ and for all $|m|\ge 1$.
Therefore one has $\|\n_{h}\|_{\io}\le C\e$ for all $h \ge 0$,
for some positive constant $D$;
in particular $|\n_{m}|\le D\e$ for all $m\ge 1$.}

\*

\0{\it Proof.} 
\0The proof is done by induction on $\bar h$.
Let us define $J^{(h)}_{c,m}=[-\e,\e]$ and
call $J^{(h)}=\times_{|m|\ge 1,c=a,b} J^{(h)}_{c,m}$
and $I^{(h)}=\times_{|m|\ge 1,c=a,b} I^{(h)}_{c,m}$.

We suppose that there exists $I^{(\bar h)}$ such that, 
if $\n$ spans $I^{(\bar h)}$ then $\n_{\bar h}$
spans $J^{(\bar h)}$ and $|\n^{(c)}_{h,m}|\le D\e$ for 
$\bar h\ge h\ge 0$; we want to show that the same holds for $\bar h+1$.
Let us call $\tilde J^{(\bar h+1)}$ the interval spanned by
$\{\n^{(c)}_{\bar h+1,m}\}_{|m|\ge 1,c=a,b}$
when $\{\n^{(c)}_{m}\}_{|m|\ge 1,c=a,b}$ span $I^{(\bar h)}$.
For any $\{\n^{(c)}_{m}\}_{|m|\ge 1,c=a,b}\in I^{(\bar h)}$ one has 
$\{\n_{\bar h+1,m}\}_{|m|\ge 1,c=a,b} \in [-2\e-D\e^2, 2\e+D\e^{2}]$,
where the bound \equ(5.17) has been used.
This means that $J^{(\bar h+1)}$
is strictly contained in $\tilde J^{(\bar h+1)}$.
On the other hand it is obvious that there
is a one-to-one correspondence between 
$\{\n^{(c)}_{m}\}_{|m|>1,c=a,b}$
and the sequence $\{\n^{(c)}_{h,m}\}_{|m|\ge 1,c=a,b}$,
$\bar h+1\ge h\ge 0$. Hence there is a set
$I^{(\bar h+1)}\subset I^{(\bar h)}$
such that, if $\{\n^{(c)}_m\}_{|m|\ge 1,c=a,b}$
spans $I^{(\bar h+1)}$, then $\{\n^{(c)}_{\bar h+1,m}\}_{|m|\ge 1,c=a,b}$
spans the interval $J^{(\bar h)}$ and, for $\e$ small enough, 
$|\n_{h}|_{\io}\le C\e$ for $\bar h+1\ge h\ge 0$. 

The previous computations also
show that the inductive hypothesis is verified also for $\bar h=0$
so that we have proved that there exists a decreasing sets
of intervals $I^{(\bar h)}$ such that if

$\{\n^{(c)}_{m}\}_{|m|>1,c=a,b} \in I^{(\bar h)}$
then the sequence $\{\n^{(c)}_{h,m}\}_{|m|\ge 1,c=a,b}$
is well defined for $h\le\bar h$ and it verifies
$|\n_{h,m}^{(c)}|\le C|\e|$. In order to prove the bound on the size of
$I^{(\bar h)}_{c,m}$ let us denote by $\{\n^{(c)}_{h,m}\}_{|m|\ge 1,c=a,b}$
and $\{\n^{\prime(c)}_{h,m}\}_{|m|\ge 1,c=a,b}$,
$0\le h\le \bar h$, the sequences corresponding to
$\{\n^{(c)}_{m}\}_{|m|\ge 1,c=a,b}$ and

$\{ \n^{\prime(c)}_{m} \}_{|m|\ge 1,c=a,b}$ in 
$I^{(\bar h)}$. We have
$$ \m\n^{(c)}_{h+1,m} - \m\n^{\prime(c)}_{h+1,m}=
2 \left( \m\n^{(c)}_{h,m}-\m\n^{\prime(c)}_{h,m} \right) +
\b^{(c)}_{h,m} - \b^{\prime(c)}_{h,m} ,
\Eq(5.64a)$$
where $\b^{(c)}_{h,m}$ and $\b^{\prime(c)}_{h,m}$ are
shorthands for the beta functions. Then, as
$|\n_k-\n'_k|_{\io} \le |\n_h-\n'_h|_{\io}$ for all $k\le h$, we have
$$ | \n_h-\n'_h |_{\io} \le {1 \over 2} |\n_{h+1}-\n'_{h+1} |_{\io}
+ D \e^{2} |\n_h-\n'_h |_{\io} .
\Eq(5.65) $$
Hence if $\e$ is small enough then one has
$$ \|\n-\n'\|_{\io} \le (\sqrt{2})^{-(\bar h+1)}
\|\n_{\bar h}-\n'_{\bar h}\|_{\io} .
\Eq(5.66a) $$
Since, by definition, if $\n$ spans $I^{(\bar h)}$,
then $\n_{\bar h}$ spans the interval $J^{(\bar h)}$,
of size $2 |\e|$, the size of $I^{(\bar h)}$ is
bounded by $2|\e| (\sqrt{2})^{(-\bar h-1)}$.

Finally note that  one can choose
$\n^{(c)}_{m}=\n^{(c)}_{-m}$ and then
$\n^{(c)}_{h,m}=\n^{(c)}_{h,-m}$ for any $|m|\ge 1$ and any 
$\bar h\ge h\ge 0$; this follows from the fact
that the function $\b_{k,m}^{(c)}$ in \equ(5.15) is even
under the exchange $m\to-m$; it depends on $m$ through $\tilde \o_m$
(which is an even function of $m$), through the end-points
$v\in E(\th)$ (which are odd under the exchange $m\to-m$; but their 
number must be even) and finally through 
$\nu^{(q-1)}_{k,m}$ which are assumed inductively to be even.\qed
\*
It will be useful to explicitly
construct the $\n^{(c)}_{h,m}$ by a contraction method.
By iterating \equ(5.14) we find, for $|m|\ge 1$
$$ \m\n^{(c)}_{h,m}=2^{h+1} \left( \m\n^{(c)}_m+\sum_{k=-1}^{h-1}
2^{-k-2} \b^{(c)}_{k,m}(\tilde\o,\e,\{ \n^{(c')}_{k',m'} \} ) \right) ,
\Eq(5.22a)$$
where $\b^{(c)}_{k,m}(\tilde\o,\e,\{ \n^{(c')}_{k',m'} \} )$ depends on
$\n_{k',m'}$ with $k'\le k-1$. If we put $h=\bar h$ 
in \equ(5.22a) we get
$$ \m\n_{m}^{(c)}=-\sum_{k=-1}^{\bar h-1}
2^{-k-2} \b^{(c)}_{k,m}(\tilde\o,\e,\{ \n^{(c')}_{k',m'} \} )
+2^{-\bar h-1}\m\n_{\bar h,m}^{(c)}
\Eq(5.23a) $$
and, combining \equ(5.22a) with \equ(5.23a), we find,
for $\bar h> h\ge 0$,
$$ \m\n^{(c)}_{h,m}=-2^{h+1} \left( \sum_{k=h}^{\bar h-1}
2^{-k-2} \b^{(c)}_{k,m}(\tilde\o,\e, \{ \n^{(c')}_{k',m'} \})
\right) + 2^{h-\bar h} \m\n_{\bar h,m}^{(c)} .
\Eq(5.24b) $$
The sequences $\{\n^{(c)}_{h,m}\}_{|m|>1}$, $\bar h> h\ge h_{0}$,
parametrized by $\{\n^{(c)}_{\bar h,m}\}_{|m| \ge 2}$
such that $\|\n^{(c)}_{\bar h}\|_\io\le C\e$,
can be obtained as the limit as $q\to\io$
of the sequences $\{\n^{(c)(q)}_{h,m}\}$, $q\ge 0$,
defined recursively as
$$ \eqalign{
\m\n_{h,m}^{(c)(0)} & = 0 , \cr
\m\n_{h,m}^{(c)(q)} & = - 2^{h+1}\left(\sum_{k=h}^{\bar h-1}
2^{-k-2} \b^{(c)}_{k,m}(\tilde\o,\e,
\{\n^{(c')(q-1)}_{k',m'} \})\right)+2^{h-\bar h} \m\n_{\bar h,m}^{(c)} . \cr}
\Eq(5.25) $$
In fact, it is easy to show inductively  that, if
$\e$ is small enough, $\|\n_h^{(q)}\|_{\io}\le C\e$, 
so that \equ(5.25) is meaningful, and
$$ \max_{0\le h\le \bar h} \|\n_h^{(q)}-\n_h^{(q-1)} \|_{\io}
\le (C\e)^{q} .
\Eq(5.26) $$
For $q=1$ this is true as $\n_h^{(c)(0)}=0$;
for $q>1$ it follows by the fact that
$\b^{(c)}_k(\tilde\o,\e, \{\n^{(c')(q-1)}_{k',m'} \})-
\b^{(c)}_k(\tilde\o,\e, \{\n^{(c')(q-2)}_{k',m'} \})$
can be written as a sum of terms
in which there are at least one $\n$-vertex,
with a difference $\n_{h'}^{(c')(q-1)}-\n_{h'}^{(c')(q-2)}$,
with $h'\ge k$, in place of the corresponding
$\nu_{h'}^{(c')}$, and one node carrying an $\e$.
Then $\n_h^{(q)}$ converges as $q\to\io$,
for $\bar h<h\le 1$, to a limit $\n_h$, satisfying
the bound $\|\n_h\|_{\io}\le C\e$. 
Since the solution is unique, it must coincide with one in
Lemma 10.

\*

We have then constructed a sequence of $\n_{h,m}^{(c)}$
solving \equ(5.24b) for any $\bar h>1$ and any $\n_{\bar h,m}^{(c)}$;
we shall call $\n_{h,m}^{(c)}(\G)$ the solution of \equ(5.24b)
with $\bar h=\io$ and $\n^{(c)}_{\io,m}=0$,
to stress the dependence on $\G$.

We will prove the following Lemma.

\*

\0{\bf Lemma 11.} {\it Under the 
the same conditions of Lemma 10
it holds that for any $h\ge 0$
$$ \|\nu_h(\G^1)-\nu_h(\G^2)\|_{\io}\le D \e |\G^1-\G^2| ,
\Eq(5.100)$$
for a suitable constant $D$.}

\*

\0{\it Proof.} Calling 
$\n_h^{(q)}(\G)$ the l.h.s. of \equ(5.25)
with $\bar h=\io$ and $\n_{\io,m}=0$, we can 
show by induction on $q$ that
$$ \|\nu_h^{(q)}(\G^1)-\nu_h^{(q)}(\G^2)\|_{\io}\le D \e |\G^1-\G^2| ,
\Eq(5.101) $$
We find convenient to write explicitly the dependence
of the function $\b^{(c)}_{h,m}$ from the parameter 
$\G$, so that we rewrite
$\b^{(c)}_{k,m}(\tilde\o,\e,
\{\n^{(c')(q-1)}_{k',m'} \}$
in the r.h.s. of \equ(5.25) as 
$\b^{(c)}_{k,m}(\tilde\o,\e,\G,
\{\n^{(c')(q-1)}_{k',m'}(\G) \}$.
Then from \equ(5.25) we get
$$\eqalign{
& \m\n_{m}^{(c)(q)}(\G^1)-\m\n_{m}^{(c)(q)}(\G^2) 
 = \cr & \sum_{k=h}^{\io}
2^{h-k-1} 
[\b^{(c)}_{k,m}(\tilde\o,\e,\G^1,
\{\n^{(c')(q-1)}_{k',m'}(\G^1) \})-
\b^{(c)}_{k,m}(\tilde\o,\e,\G^2,
\{\n^{(c')(q-1)}_{k',m'}(\G^2) \})]\cr}
\Eq(5.25aa) $$
When $q=1$ we have that $\b^{(c)}_{k,m}(\tilde\o,\e,\G^1,
\{\n^{(c')(q-1)}_{k',m'}(\G^1) \})-
\b^{(c)}_{k,m}(\tilde\o,\e,\G^2,
\{\n^{(c')(q-1)}_{k',m'}(\G^2) \})$
is given by a sum of self energy graphs with
one node $\vvv$ with a factor $\h_{\vvvv}$ 
with $\G$ replaced by $\G^1-\G^2$; as there is at least 
a vertex $\vvv$ of w-type  by the definition of the self energy graphs
we obtain
$$ \|\nu_h^{(1)}(\G^1)-\nu_h^{(1)}(\G^2)\|_{\io}\le
\left( D_{1} \e+\tilde D_1\e^2 \right) |\G^1-\G^2| ,
\Eq(5.102) $$
for positive constants $D_{1}<2 D$ and $\tilde D_2$, where 
$D_{1} \e|\G^1-\G^2|$ is a bound for the self-energy first order contribution.

For $q>1$ we can write the difference
in \equ(5.25aa) as 

$$\eqalign{
& \left( \b^{(c)}_{k,m}(\tilde\o,\e,\G^1,
\{\n^{(c')(q-1)}_{k',m'}(\G^1) \})-
\b^{(c)}_{k,m}(\tilde\o,\e,\G^2,
\{\n^{(c)(q-1)}_{k',m'}(\G^1) \}) \right) \cr
& + \qquad \left( \b^{(c)}_{k,m}(\tilde\o,\e,\G^2,
\{\n^{(c)(q-1)}_{k',m'}(\G^1) \})-
\b^{(c)}_{k,m}(\tilde\o,\e,\G^2,
\{\n^{(c)(q-1)}_{k',m'}(\G^2) \}) \right) . \cr}
\Eq(5.677)$$
The first factor is given by a sum over self-energy graphs
with one node $\vvv$ with a factor $\h_{\vvvv}$ 
with $\G$ replaced by $\G^1-\G^2$; the other difference
is given by a sum over self energy graphs with a $\n$-vertex
to which is associated a factor
$\n^{(c)(q-1)}_{k',m'}(\G^1)-\n^{(c)(q-1)}_{k',m'}(\G^2)$;
hence we find
$$ \|\nu_h^{(q)}(\G^1)-\nu_h^{(q)}(\G^2)\|_{\io} \le
\left( D_{1} \e + D_3\e^2 \right)
|\G^1-\G^2| + \e D_2 \sup_{h \ge 0}
\|\n_h^{(q-1)}(\G^1)- \n_h^{(q-1)}(\G^2) \|_\io ,
\Eq(5.103) $$
where $D_1\e |\G^1-\G^2|$ is a bound for
the first order contribution coming from the first line in \equ(5.677),
while the last summand in \equ(5.103)
is a bound from the terms from the last line of \equ(5.677).  
Then \equ(5.101) follows with $D=2D_{1}$, for $\e$ small enough. \qed

\*

By using Lemma 11 we can show that the self consistence equation
for $\G$ \equ(5.3a) has a unique solution $\G=O(\e)$.

\*

\0{\bf Lemma 12.} {\it 
For all $\m_0>0$ there exists $\e_0>0$
such that, for all $|\m| \le \m_0$ and
for all $0<\e <\e_{0}$, given $\n_{m}^{(c)}(\G)$
chosen as in Lemma 9
(with $\bar h=\io$ and $\n_{\io,m}^{(c)}=0$)
it holds that \equ(5.3a) has a solution $|\G|\le C\e$
where $C$ is a suitable constant.}

\*

\0{\it Proof.} 
The solution of \equ(5.3a) can be
obtained as the limit as $q\to\io$ of the sequence $\G^{(q)}$,
$q\ge 0$, defined recursively as
$$ \eqalign{
\G^{(0)} & = 0 , \cr
\G^{(q)} & = r_{0}^{-1} 
\sum_{k=1}^{\io} {\mathop{\sum}_{\th \in \Th^{(k)(q-1))\RR}_{n,n}}}^{*}
a_{0,-n} \Val(\th)  , \cr}
\Eq(5.256) $$
where we define $\Th^{(k)(q)\RR}_{n,m}$ as the set of trees
identical to  $\Th^{(k)\RR}_{n,m}$
except that $\G$ in $\h_{\vvvv}$
is replaced by $\G^{(q)}$ and each $\n_{h,m}$ is replaced by 
$\n_{h,m}(\G^{(q)})$, for all $h\ge 0, |m|\ge 1$.
\equ(5.256) is a contraction
defined on the set $|\G|\le C\e$, for $\e$ small.
In fact if $|\G^{(q-1)}|\le C\e$, then by \equ(5.256)
$|\G^{(q)}|\le C_1\e+C_2 C\e^2$, where we have used
that the first order contribution to the r.h.s.
of \equ(5.256) is $\G$-indipendent (see Section \secc(3)),
and $C_1\e<\e C /$ is a bound for it;
hence for $\e$ small enough \equ(5.256)
send the interval $|\G|\le C\e$ to itself. 

Moreover we can show inductively  that
$$\|\G^{(q)}-\G^{(q-1)} \|_{\io}
\le (C\e)^{q} .
\Eq(5.266) $$
For $q=1$ this is true;
for $q>1$ $\G^{(q)}-\G^{(q-1)}$
can be written as sum of trees
in which a)or to a node $\vvv'$
is associated a factor proportional to $\G^{(q-1)}-\G^{(q-2)}$;
or b) to a $\nu$ vertex is associated $\nu_{h',m'}(\G^{(q-1)})-
\nu_{h',m'}(\G^{(q-2)})$ for some $h',m'$. 
In the first case we note that the constraint in
the sum in the r.h.s. of \equ(5.256)
implies that $s_{v_0}=3$ for the special vertex of $\th$;
hence, item (4) in Section \secc(3), says that $\vvv'\not= v_0$
so that such terms are bounded by $D_1 \e|\G^{(1)}-\G^{(2)}|$
(a term order $O(\G^{1}-\G^2)$ should have three
$v$ lines entering $v_0$ and two of them coming from
end points, which is impossible). 
In the second case we use 
\equ(5.100), and we bound such terms by $D_2 \e|\G^{(1)}-\G^{(2)}|$ .
Hence by induction \equ(5.266) is found, if $C\ge D_1/4,
D_2/4, C_1/4$). \qed

\*

We have finally to prove that $\tilde u_{n,m}$ 
solves the last of \equ(1.13) and \equ(2.2).

\*

\0{\bf Lemma 13.} {\it 
For all $\m_0>0$ there exists $\e_0>$
such that, for all $|\m|\le\m_0$ and
for all $0<\e < \e_{0}$, given $\n_{m}^{(c)}(\G)$
chosen as in Lemma 9 and $\G$ chosen as in lemma 12
then $\tilde u_{n,m}$ solves the last of \equ(1.13) and \equ(2.2).}

\*

\0{\it Proof.} Let us consider first the case
in which $|n|\not=|m|$ and we call $\Th_{n,m}^{\RR}=\bigcup_{k}
\Th^{(k)\RR}_{n,m}$; assume also
(what of course is not restrictive) that $n,m$
is such that $\chi_{h_0}(|\o n|-\tilde\o_m)+ 
\chi_{h_0+1}(|\o n|-\tilde\o_m)=1$. 
We call $\Th_{n,m,}^{\RR}$
the set of trees $\th\in \Th_{n,m}^{\RR}$
with root line at scale $h_0$, so that  
$$ \tilde u_{n,m}=\sum_{\th\in \Th_{n,m}^{\RR}}\Val(\th)=
\sum_{\th\in \Th_{n,m,h_0}^{\RR}}\Val(\th)+
\sum_{\th\in \Th_{n,m,h_0+1}^{\RR}}\Val(\th) ,
\Eq(5.cvv)$$
and we write 
$\Th_{n,m,h_0}^{\RR}=\Th_{n,m,h_0}^{\a,\RR}\bigcup
\Th_{n,m,h_0}^{\b,\RR}$, where $\Th_{n,m,h_0}^{\a,\RR}$
are the trees with $s_{\vvvv_0}=1$, while $\Th_{n,m,\bar h}^{\b,\RR}$
are the trees with $s_{\vvvv_0}=3$ and $\vvv_0$ is the special
vertex (see Definition 2). Then
$$ \eqalign{
& \sum_{\th\in \Th_{n,m,h_0}^{\a,\RR}}\Val(\th)
= \sum_{c=a,b} \Big( g^{(\bar h)}(n,m) 2^{-h_0}\n_{h_0,m}^{(c)} 
\sum_{\th \in \Th^{\RR}_{n,m_c,h_0}}
\Val(\th) \cr
& \qquad + g^{(\bar h)}(n,m) 2^{-\bar h}\n_{h_0,m}^{(c)}
\sum_{\th \in \Th^{\RR}_{n,m_c,h_0+1}}
\Val(\th) \cr
& \qquad + g^{(h_0+1)}(n,m) 2^{-h_0}\n^{(c)}_{h_0,m}
\sum_{\th \in \Th^{\RR}_{n,m,h_0}}
\Val(\th) \cr
& \qquad + g^{(h_0+1)}(n,m) 2^{-h_0-1}\n^{(c)}_{h_0+1,m}
\sum_{\th \in \Th^{\RR}_{n,m_c,h_0+1}}
\Val(\th) \Big)  , \cr}
\Eq(5.1aaaa)
$$
where $m_c$ is such that $m_a=m$ and $m_b=-m$.
On the other hand we can write 
$\Th_{n,m,h_0 }^{\b,\RR}=\Th_{n,m,h_0}^{\b1,\RR}\bigcup
\Th_{n,m,h_0}^{\b2,\RR}$, where 
$\Th_{n,m,h_0}^{\b1,\RR}$ are the trees 
such that the root line $\ell_0$
is the external line of a self-energy graph,
and $\Th_{n,m,h_0}^{\b2,\RR}$ is the complement.
Then we can write
$$\eqalign{
&\sum_{\th\in\Th_{n,m,h_0}^{\b1,\RR}}
\Val(\th)=\sum_{c=a,b} 
g^{(h_0)}(n,m) \sum_{T \in \tilde\TT^{(c)}_{< h_0}}
\RR\VV_T^{\bar h}(\o n,m) 
\sum_{\th \in \Th^{\RR}_{n,m_c,h_0}}
\Val(\th)\cr
& \qquad +g^{(h_0)}(n,m) \sum_{T \in \tilde\TT^{(c)}_{< h_0}}
\RR\VV_T^{h_0}(\o n,m) 
\sum_{\th \in \Th^{\RR}_{n,m_c,h_0+1}}
\Val(\th) \cr 
& \qquad + g^{(h_0+1)}(n,m) \sum_{T \in 
\tilde\TT^{(c)}_{< h_0}}
\RR\VV_T^{h_0}(\o n,m)
\sum_{\th \in \Th^{\RR}_{n,m_c,h_0}}
\Val(\th) \cr
& \qquad + g^{(h_0+1)}(n,m) \sum_{T \in \tilde\TT^{(c)}_{< h_0+1}}
\RR\VV_T^{h_0+1}(\o n,m)
\sum_{\th \in \Th^{\RR}_{n,m_c,h_0+1}}
\Val(\th) , \cr}
\Eq(5.bbb) $$
where $\tilde\TT^{(c)}_{<h_0}$
is the set of self-energy graphs such that if $\vvv$
is the vertex to which the external line $\ell_0$
is attached, than $s_{\vvvv}=3$. Note that if $T$ belongs
to the complementary of $\tilde\TT^{(c)}_{< \bar h}$,
then $\LL\VV_T^h=0$ by the compact support properties of $g^{(h)}$. 

Summing \equ(5.1aa) and \equ(5.bbb)
and using \equ(5.1), \equ(5.3) we get 
$$\eqalign{
& \sum_{\th\in\Th_{n,m,h_0}^{\b1,\RR}}
\Val(\th)=\sum_{c=a,b} 
g^{(h_0)}(n,m) \sum_{T \in \TT^{(c)}_{< h_0}}
\VV_T^{h_0}(\o n,m) 
\sum_{\th \in \Th^{\RR}_{n,m_c,h_0}}
\Val(\th)\cr
& \qquad +g^{(h_0)}(n,m) \sum_{T \in \TT^{(c)}_{< h_0}}
\VV_T^{h_0}(\o n,m) 
\sum_{\th \in \Th^{\RR}_{n,m,h_0+1}}
\Val(\th) \cr
& \qquad + g^{(h_0+1)}(n,m) \sum_{T \in \TT^{(c)}_{< h_0}}
\VV_T^{h_0}(\o n,m)
\sum_{\th \in \Th^{\RR}_{n,m_c,h_0}}
\Val(\th) \cr
& \qquad + g^{(h_0+1)}(n,m) \sum_{T \in \TT^{(c)}_{< h_0+1}}
\VV_T^{h_0+1}(\o n,m)
\sum_{\th \in \Th^{\RR}_{n,m_c,h_0+1}}
\Val(\th) \cr
& \qquad +\n^{(c)}_m \left( g^{(h_0)}(n,m)+g^{(h_0+1)}(n,m) \right)
\Big( \sum_{\th \in \Th^{\RR}_{n,m_c,h_0}}
\Val(\th)+\sum_{\th \in \Th^{\RR}_{n,m_c,h_0+1}}
\Val(\th) \Big) . \cr}
\Eq(5.bbbb) $$
The last line is equal to 
$$ {1 \over -\o^2 n^2 + \tilde\o_{m}^2}
[\n_{m}^{(a)} \tilde u_{n,m}+ \n_{m}^{(b)} \tilde u_{n,-m}] ,
\Eq(5.ss)$$
while adding the first three lines in \equ(5.bbb)
to $\sum_{\th\in\Th_{n,m,h_0}^{\b1,\RR}}
\Val(\th)$ we get
$$\e\sum_{n_1+n_2+n_3=n\atop m_1+m_2+m_3=m}
\tilde u_{n_1,m_1}\tilde u_{n_2,m_2}\tilde u_{n_3,m_3} ,
\Eq(5.sss)$$
from which we get that 
$\sum_{\th\in \Th_{n,m}^{\RR}}\Val(\th)$, for $\m=1$,
is a formal solution of \equ(2.2).
A similar result holds for $|n|=|m|$. \qed

\*\*
\section(6,Construction of the perturbed frequencies)

\0In the following it will be convenient to set
$\tilde\o=\{\o_{m}\}_{|m|\ge 2}$.
By the analysis of the previous sections we have 
found the counterterms $\{\n_m(\tilde\o,\e)\}_{|m|\ge 2}$
as functions of $\e$ and $\tilde\o$.
We have now to invert the relations
$$ \tilde\o^2_m-\n_{m}(\tilde\o,\e) = m^{2} ,
\Eq(6.1)$$
in order to prove Proposition 2.

We shall show that there exists a sequence of sets
$\{\EE^{(p)}\}_{p=0}^{\io}$ in $[0,\e_{0}]$, such that
$\EE^{(p+1)} \subset \EE^{(p)}$, and a sequence of functions
$\{\tilde\o^{(p)}(\e)\}_{p=0}^{\io}$, with each
$\tilde\o^{(p)} \= \tilde\o^{(p)}(\e)$ defined
for $\e\in\EE^{(p)}$, such that for all $\e\in\EE$, with
$$ \EE = \bigcap_{p=0}^{\io} \EE^{(p)} =
\lim_{p\to\io} \EE^{(p)} ,
\Eq(6.2) $$
there exists the limit
$$ \tilde\o^{(\io)}(\e) = \lim_{p\to\io} \tilde\o^{(p)}(\e) ,
\Eq(6.3) $$
and it solves \equ(6.1).

To fulfill the program above we shall define since the beginning
$\o=\o_{\e}=\sqrt{1-\e}$, and we shall follow an
iterative scheme by setting, for $|m|\ge 1$,
$$ \eqalign{
& \tilde \o_{m}^{(0)2} = m^2 , \cr
& \tilde \o_{m}^{(p)2} \= \tilde \o_m^{(p)2}(\e) =
m^{2} + \n_{m}(\tilde\o^{(p-1)},\e) , \qquad p\ge 1 , \cr}
\Eq(6.4) $$
and reducing recursively the set of admissible values of $\e$.

We start by imposing on $\e$ the Diophantine conditions
$$ | \o n \pm m |
\ge 2 C_{0} | n|^{-\t_{0}}
\qquad \forall  n\in\ZZZ \setminus \{ 0 \} \hbox{ and }
\forall m \in \ZZZ \setminus \{0\} \hbox{ such that } |m| \neq |n| ,
\Eq(6.5) $$
where $C_{0}$ and $\t_{0}$ are two positive constants.

This will imply some restrictions on the admissible
values of $\e$, as the following result shows.

\*

\0{\bf Lemma 14.}
{\it For all $0 < C_{0} \le 1/2$ there exist $\e_{0}>0$
and $\g_{0},\d_{0}>0$ such that the set $\EE^{(0)}$ of
values $\e\in[0,\e_{0}]$ for which \equ(6.5) are satisfied has
Lebesgue measure $|\EE^{(0)}|\ge \e_{0}(1-\g_{0}\e_{0}^{\d_{0}})$
provided that one has $\t_{0}>1$.}

\*

\0{\it Proof.}
For $(n,m)$ such that $|\o n \pm m|\ge 2C_0$ 
the Diophantine conditions in \equ(6.5) are trivially satisfied.
We consider then $(n,m)$ such that $|\o n \pm m|<2C_0$
and we write, if $0 < C_0 \le 1/2$
$$ 1 > 2 C_0 > |\o n\pm m| \ge
\left| {\e n\over 1 + \sqrt{1-\e}} - n \pm m \right| ,
\Eq(6.6)$$
and as $|n\pm m|\ge 1$ one gets $|n| \ge 1/2\e \ge 1/2\e_{0}$.
Moreover, for fixed $n$, the set $\MM$ of $m$'s such that
$|\o n \pm m|<1$ contains at most $2+\e_{0}|n|$ values. By writing
$$ f(\e(t)) = n \sqrt{1-\e(t)} \pm m = t {2C_{0} \over |n|^{\t_{0}}} ,
\qquad t\in[0,1] ,
\Eq(6.7) $$
and, calling $\III^{(0)}$ the set of $\e$ such that  $| \o n \pm m |
< 2C_{0} | n|^{-\t_{0}}$ is verified for some $(n,m)$,
one finds for the Lebesgue measure of $\III^{(0)}$
$$ {\rm meas}(\III^{(0)}) = \int_{\III^{(0)}} {\rm d}\e =
\sum_{|n|\ge 1/2\e_{0}} \sum_{|m| \in \MM}
\int_{-1}^{1} {\rm d}t \left|
{{\rm d}\e(t) \over {\rm d}t} \right| .
\Eq(6.8) $$
We have from \equ(6.7)
$$ {{\rm d} f \over {\rm d}t} = {{\rm d} f \over {\rm d} \e}
{{\rm d} \e \over {\rm d} t} = { 2 C_{0} \over |n|^{\t_0} } ,
\Eq(6.9) $$
so that, noting that one has $|\dpr f/\dpr \e| \ge |n|/4$, 
we have to exclude a set of measure
$$ \sum_{|n| \ge N} { 8 C_{0} \over |n|^{\t_{0}+1}}
\left( 2 + \e_{0} |n| \right) \le
{\rm const.} \, \e_{0}^{\t_{0}} ,
\Eq(6.10) $$
and one has to impose $\t_{0}=1+\d_{0}$, with $\d_{0}>0$. \qed

\*

For $p\ge 1$ the sets $\EE^{(p)}$ will be defined recursively as
$$ \eqalign{
\EE^{(p)} & = \Big\{ \e \in \EE^{(p-1)} \;:\;
| \o n \pm \tilde\o_{m}^{(p)} | >
C_{0} | n|^{-\t} \qquad \forall |m| \neq |n|, \cr
& \qquad \qquad
| \o n \pm ( \tilde\o_{m}^{(p)} \pm
\tilde \o_{m'}^{(p)} ) | > C_{0} | n|^{-\t} \qquad
|n| \neq |m \pm m'| \Big\} ,
\qquad p \ge 1 , \cr}
\Eq(6.11) $$
for $\t>\t_{0}$ to be fixed.

In Appendix \secc(A4) we prove the following result.

\*

\0{\bf Lemma 15.} {\it For all $p\ge 1$ one has
$$ \left\| \tilde\o^{(p)}(\e)-\tilde\o^{(p-1)}(\e) \right\|_{\io}
\le C \e_{0}^{p} \qquad \forall \e\in \EE^{(p)} ,
\Eq(6.12) $$
for some constant $C$.}

\*

Therefore we can conclude that there exist a sequence
$\{\tilde\o^{(p)}(\e)\}_{p=0}^{\io}$
converging to $\tilde\o^{(\io)}(\e)$ for $\e\in\EE$.
We have now to show that the set $\EE$ has positive (large) measure.

It is convenient to introduce a set of
variables $\m(\tilde\o,\e)$ such that
$$ \tilde\o_{m} + \m_{m}(\tilde\o,\e) = \o_{m} \= |m| ;
\Eq(6.13) $$
the variables $\m(\tilde\o,\e)$ and the counterterms are
trivially related by
$$ - \n_{m}(\tilde\o,\e) = \m_{m}^{2}(\tilde\o,\e) +
2\tilde\o_{m}\m_{m}(\tilde\o,\e) .
\Eq(6.14) $$

One can write $\tilde\o^{(p)}_{m}=\o_{m}-\m_{m}(\tilde\o^{(p-1)},\e)$,
according to \equ(6.4). We shall impose the Diophantine conditions
$$ \eqalignno{
& | \o n \pm ( \o_{m} - \m_{m}(\tilde\o^{(p-1)},\e))| > C_{0} | n|^{-\t} ,
& \eq(6.15) \cr
& | \o n \pm ( ( \o_{m} \pm \o_{m'}) -
(\m_{m}(\tilde\o^{(p-1)},\e) \pm \m_{m'}(\tilde\o^{(p-1)},\e) )) | >
C_{0} | n|^{-\t} . \cr} $$
Suppose that for $\e\in\EE^{(p-1)}$ the functions
$\m_{m}(\tilde\o^{(p-1)},\e)$ are well defined; then define
$\III^{(p)}=\III^{(p)}_{1}\cup \III^{(p)}_{2}\cup\III^{(p)}_{3}$,
where $\III^{(p)}_{1}$ is the set of values $\e \in \EE^{(p-1)}$
verifying the conditions
$$ \left| \o n \pm \left( \o_{m} - \m_{m}(\tilde\o^{(p-1)},\e)
\right) \right| \le C_{0} \left| n \right|^{-\t} ,
\Eq(6.16) $$
$\III^{(p)}_{2}$ is the set of values $\e$ verifying the conditions
$$ \left| \o n \pm \left( ( \o_{m} - \o_{m'}) -
(\m_{m}(\tilde\o^{(p-1)},\e) - \m_{m'}(\tilde\o^{(p-1)},\e) )
\right) \right| \le C_{0} | n|^{-\t} ,
\Eq(6.17) $$
and $\III^{(p)}_{3}$ is the set of values $\e$ verifying the conditions
$$ \left| \o n \pm \left( ( \o_{m} + \o_{m'}) -
(\m_{m}(\tilde\o^{(p-1)},\e) + \m_{m'}(\tilde\o^{(p-1)},\e) ) \right)
\right| \le C_{0} | n|^{-\t} .
\Eq(6.18) $$
For future convenience we shall call, for $i=1,2,3$, $\III^{(p)}_{i}(n)$
the subsets of $\III^{(p)}_{i}$ which verify the Diophantine conditions
\equ(6.16), \equ(6.17) and \equ(6.18), respectively, for fixed $n$.

We want to bound the measure of the set $\III^{(p)}$.
First we need to know a little better the dependence on $\e$ and
$\tilde \o$ of the counterterms: this is provided by the following result.

\*

\0{\bf Lemma 16.} {\it For all $p\ge 1$ and for
all $\e\in\EE^{(p)}$ there exists a positive constant $C$ such that
$$ \eqalign{
& \left| \n_{m}(\tilde\o^{(p)},\e) \right| \le C \e , \qquad
\left| \dpr_{\tilde\o^{(p)}_{m'}} \n_{m}(\tilde\o^{(p)},\e) \right|
\le C \e , \qquad m\ge 1 , \cr
& \left| \m_{m}(\tilde\o^{(p)},\e) \right| \le {C \e \over m} , \qquad
\left| \dpr_{\tilde\o^{(p)}_{m'}} \m_{m}(\tilde\o^{(p)},\e) \right|
\le {C \e \over m} , \qquad m\ge 1 , \cr
& \left| \dpr_{\e} \m^{(p)}_{m}(\o,\e) \right| \le C ,
\qquad m\ge 1 , \cr}
\Eq(6.19) $$
where the derivatives are in the sense of Whitney {\rm \cita{Wh}}.}

\*

\0{\it Proof.}
The bound for $\n_{m}(\tilde\o^{(p)},\e)$ is obvious by construction.

In order to prove the bound for
$\dpr_{\o^{(p)}_{m'}} \n_{m}(\tilde\o^{(p)},\e)$ note that one has
$$ \left| g^{(h_{\ell_{j}})}_{\ell_{j}}(\tilde \o') -
g^{(h_{\ell_{j}})}_{\ell_{j}}(\tilde \o)- (\tilde\o-\tilde\o') \,
\dpr_{\tilde\o} g^{(h_{\ell_{j}})}_{\ell_{j}}(\tilde \o)\right|
\le C \, n_{\ell_{j}}^{2} 2^{-3h_{\ell_{j}}}
\left\| \tilde \o' - \tilde \o \right\|^2_{\io} ,
\Eq(6.20) $$
from the compact support properties of the propagator.

Let us consider the quantity $\nu(\tilde\o',\e)-
\nu(\tilde\o,\e) - (\o-\o') \, \dpr_{\tilde\o} \nu(\tilde\o,\e)$,
where $\dpr_{\tilde\o} \nu(\tilde\o,\e)$ denotes the 
derivative in the sense of Whitney,
and note that it can be expressed as a sum over
trees each one containing a line with propagator
$g^{(h_{\ell_{j}})}_{\ell_{j}}(\tilde \o') -
g^{(h_{\ell_{j}})}_{\ell_{j}}(\tilde \o)- (\tilde\o-\tilde\o') \,
\dpr_{\tilde\o} g^{(h_{\ell_{j}})}_{\ell_{j}}(\tilde \o)$,
by proceeding as in the proof of lemma 9 of \cita{GM}. Then we find
$$ \left\| \nu(\tilde\o',\e)-\nu(\tilde\o,\e) - (\tilde\o-\tilde\o')
\, \dpr_{\tilde\o} \n(\tilde\o,\e) \right\|_{\io} \le
C \e \left\| \tilde \o' - \tilde \o \right\|^2_{\io} ,
\Eq(6.21) $$
and the second bound in \equ(6.19) follows.

The bounds for $\m_{m}(\tilde\o^{(p)},\e)$ and
$\dpr_{\o^{(p)}_{m'}} \m_{m}(\tilde\o^{(p)},\e)$ simply follow
from \equ(6.14) which gives
$$ \m_{m} (\tilde\o^{(p)},\e) =
{ \n_{m}(\tilde\o^{(p)},\e) \over |m| } \,
{- 1 \over 1 + \sqrt{1 - \n_{m}(\tilde\o^{(p)},\e)/m} } .
\Eq(6.22) $$

In order to prove the last bound in \equ(6.19) we prove
by induction the bound
$$ \left| \dpr_{\e} \tilde\o^{(p)}_{m}(\o,\e) \right| \le C ,
\qquad m\ge 1 ,
\Eq(6.23) $$
by assuming that it holds for $\tilde \o^{(p-1)}$;
then from \equ(6.4) we have
$$ \eqalignno{
& 2\tilde\o^{(p)}\dpr_{\e} \tilde\o_m^{(p)}(\e) = -
\dpr_{\e} \n_m(\tilde\o^{(p-1)}(\e),\e)
& \eq(6.24) \cr
& \qquad = - \sum_{m'\in\zzz} \dpr_{\tilde\o^{(p-1)}_{m'}}
\n_m(\tilde\o^{(p-1)}(\e),\e) \, \dpr_{\e} \tilde\o_{m'}^{(p-1)}(\e)
- \left. \dpr_{\h} \n_m(\tilde\o_m^{(p-1)}(\e),\h)
\right|_{\h=\e}, \cr} $$
as it is easy to realize by noting that $\m_{m}$ can
depend on $\tilde \o_{m'}$ when $|m'|>|m|$ only if
the sum of the absolute values of the mode labels $m_{v}$
is greater than $|m'-m|$, while we can bound the sum of the
contributions with $|m'|<|m|$ by a constant times $\e$,
simply by using the second line in \equ(6.19).

Hence from the inductive hypothesis and the proved bounds
in \equ(6.19), we obtain 
$$ \left\| \tilde\o^{(p)}(\e')-\tilde\o^{(p)}(\e)-
(\e-\e') \, \dpr_\e \tilde\o^{(p)}(\e) \right\|_{\io}
\le C \left| \e' - \e \right| ,
\Eq(6.25) $$
so that also the bound \equ(6.23) and hence the last
bound in \equ(6.16) follow. \qed

\*

Now we can bound the measure of the set we have to exclude:
this will conclude the proof of Proposition 2.

We start with the estimate of the measure of the set $\III^{(p)}_{1}$.
When \equ(6.16) is satisfied one must have
$$ \eqalign{
& C_{2} |m| \le | \o_{m} - \m_{m}(\tilde\o^{(p-1)},\e) |
\le \o|n| + C_{0} | n|^{-\t} \le C_{2}'| n| , \cr
& C_{1} |m| \ge | \o_{m} - \m_{m}(\tilde\o^{(p-1)},\e) |
\ge \o|n| - C_{0} | n|^{-\t} \ge C_{1}'| n| , \cr}
\Eq(6.26) $$
which implies
$$ M_{1} |n| \le |m| \le M_{2} |n| ,
\qquad M_{1} = {C_{1} \over C_{1}'} ,
\qquad M_{2} = {C_{2}' \over C_{2}} ,
\Eq(6.27) $$
and it is easy to see that, for fixed $n$, the set $\MM_{0}(n)$ of $m$'s
such that \equ(6.16) are satisfied contains at most $2+\e_{0}|n|$ values.

Furthermore (by using also \equ(6.5)) from \equ(6.16)
one obtains also
$$ \eqalign{
2 C_{0} | n|^{-\t_{0}} & \le | \o |n| - \o_{m} |
\le | \o |n| - \o_{m} + \m_{m}(\tilde \o^{(p-1)},\e) | +
| \m_{m}(\tilde \o^{(p-1)},\e) | \cr
& \le C_{0} | n|^{-\t} + C\e_{0}/|m| , \cr}
\Eq(6.28) $$
which implies, together with \equ(6.36), for $\t \ge \t_{0}$,
$$ |n| \ge N_{0} \= \left( { C_{0} M_{1} \over C \e_{0} }
\right)^{1/(\t_{0}-1)} .
\Eq(6.29) $$
Let us write $\o(\e)=\o_{\e}=\sqrt{1-\e}$ and
consider the function $\m(\tilde\o^{(p-1)},\e)$:
we can define a map $t\to \e(t)$ such that
$$ f(\e(t)) \= \o(\e(t)) |n| - \o_{m} + \m_{m}(\tilde\o^{(p-1)},\e(t)) =
t {C_{0} \over | n|^{\t} } , \qquad t \in [-1,1] ,
\Eq(6.30)$$
describes the interval defined by \equ(6.16); then the Lebesgue measure
of $\III_{1}^{(p)}$ is
$$ {\rm meas}(\III^{(p)}_{1}) = \int_{\III^{(p)}_1} {\rm d}\e =
\sum_{|n|\ge N_{0}} \sum_{m \in \MM_{0}(n)}
\int_{-1}^{1} {\rm d}t \left|
{{\rm d}\e(t) \over {\rm d}t} \right| .
\Eq(6.31) $$
We have from \equ(6.30)
$$ {{\rm d} f \over {\rm d}t} = {{\rm d} f \over {\rm d} \e}
{{\rm d} \e \over {\rm d} t} = { C_{0} \over |n|^{\t} } ,
\Eq(6.32) $$
hence
$$ {\rm meas}(\III^{(p)}_{1}) =
\sum_{|n|\ge N_{0}} \sum_{m \in \MM_{0}(n)}
{ C_{0} \over | n|^{\t}} \int_{-1}^{1} {\rm d} t
\left| {{\rm d} f(\e(t))\over {\rm d} \e(t)} \right|^{-1} .
\Eq(6.33) $$
In order to perform the derivative in \equ(6.30) we write
$$ {{\rm d} \m_{m}\over {\rm d}\e} = \dpr_\e \m_{m} +
\sum_{|m'|\ge 1} \dpr_{\o^{(p-1)}_{m'}}\m_{m}
\dpr_\e \tilde\o_{m'}^{(p-1)} ,
\Eq(6.34) $$
where $\dpr_\e \m_{m}$ and $\dpr_\e \tilde\o_{m'}^{(p-1)}$
are bounded through lemma 11. Moreover one has
$$ \left| \dpr_{\tilde\o^{(p-1)}_{m'}} \m_{m} \right| \le
C \e\,  e^{-\k|m'-m|/2} , \qquad |m'| > |m| ,
\Eq(6.35)$$
and the sum over $m'$ can be dealt with as in \equ(6.24).

At the end we get that the sum in \equ(6.34)
is $O(\e)$, and from \equ(6.29) and \equ(6.30) we obtain,
if $\e_{0}$ is small enough,
$$ \left| {\dpr f(\e(t)) \over \dpr \e(t)} \right|
\ge { \left| n \right| \over 4} , 
\Eq(6.36) $$
so that one has
$$ \eqalign{
\int_{\III^{(p)}_1} {\rm d}\e
& \le \hbox{const.} \sum_{|n|\ge N_{0}}
{ C_{0} \over | n|^{\t+1}}
\left( 2 + \e_{0} |n| \right) \cr
& \le {\rm const.} \, \e_{0}
\left( \e_{0}^{(\t-1)/(\t_{0}-1)} + \e_{0}^{(\t-\t_{0}+1)/(\t_{0}-1)}
\right) . \cr}
\Eq(6.37) $$
Therefore for $\t \ge \max\{1,\t_{0}-1\}$ the Lebesgue measure of the set
$\III^{(p)}_{1}$ is bounded by $\e_{0}^{1+\d_{1}}$, with $\d_{1}>0$.

Now we discuss how to bound the measure of the set
$\III^{(p)}_{2}$. We start by noting that from \equ(6.30) we obtain,
if $m,\ell>0$,
$$ \left| \tilde\o^{(p)}_{m+\ell} - \tilde\o^{(p)}_{m} -
\ell \right| \le {C \e \over m} + {C \e \over m + \ell} ,
\Eq(6.38) $$
for all $p\ge 1$.

By the parity properties of
$\tilde\o_{m}$ without loss of generality we can confine
ourselves to the case $n>0$, $m'>m \ge 2$,
and $|\o n-(\tilde\o_{m'}^{(p)}-\tilde\o_{m}^{(p)})|<1$.
Then the discussion proceeds as follows.

When the conditions \equ(6.17) are satisfied, one has
$$ \eqalign{
2 C_{0} | n|^{-\t_{0}} & \le
| \o n - (\o_{m'}-\o_{m}) | \cr
& \le | \o n -
(\o_{m'} - \m^{(p-1)}_{m'}(\tilde\o^{(p-1)},\e) ) \cr
& \qquad \qquad + (\o_{m} - \m^{(p-1)}_{m}(\tilde\o^{(p-1)},\e) ) |
\cr & \qquad \qquad + | \m_{m'}(\tilde \o^{(p-1)},\e) -
\m_{m}(\tilde \o^{(p-1)},\e) | \cr
& \le C_{0} | n|^{-\t} + {C \e \over m} + {C \e \over m'}
\le C_{0}|n|^{-\t} + {2C\e_{0} \over m} , \cr}
\Eq(6.39) $$
which implies for $\t \ge \t_{0}$
$$ | n| \ge N_{1} \= \left( { C_{0} \over 2 C \e_{0} }
\right)^{1/\t_{0}} .
\Eq(6.40) $$

We can bound the Lebesgue measure
of the set $\III^{(p)}_{2}$ by distinguishing,
for fixed $(n,\ell)$, with $\ell=m'-m>0$,
the values $m\le m_{0}$ and $m>m_{0}$,
where $m_{0}$ is determined by the request that
one has for $m>m_{0}$
$$ {2 C \e_{0} \over m} \le {C_{0} \over |n|^{\t_{0}}} ,
\Eq(6.41) $$
which gives
$$ m_{0} = m_{0}(n) = \left( {2 C |n|^{\t_{0}} \e_{0} \over C_{0}} \right) .
\Eq(6.42) $$
Therefore for $m>m_{0}$ and $\t \ge \t_{0}$
one has, from \equ(6.5), \equ(6.38) and \equ(6.41),
$$ \left| \o n -
(\tilde \o_{m'}^{(p)} - \tilde \o_{m}^{(p)} ) \right|
\ge \left| \o n - \ell \right|
- {2C\e \over m } \ge {2C_{0} \over |n|^{\t_{0}}} - 
{C_{0} \over |n|^{\t_{0}} } \ge {C_{0} \over |n|^{\t_{0}} }
\ge {C_{0} \over |n|^{\t} } ,
\Eq(6.43) $$
so that one has to exclude no further value from
$\EE^{(p-1)}$, provided one takes $\t \ge \t_{0}$.

For $m<m_{0}$ define $L_{0}$ such that
$$ C_{3} \ell \le \left| \tilde \o_{m+\ell}^{(p)} -
\tilde \o_{m}^{(p)} \right| < \o n + 1 < C_{3}' |n| , \qquad
L_{0} = { C_{3}' \over C_{3} } ,
\Eq(6.44) $$
where \equ(6.26) has been used. Again, for fixed $n$,
the set $\LL_{0}(n)$ of $\ell$'s such that \equ(6.26) is satisfied
with $m'-m=\ell$ contains at most $2+\e_{0}|n|$ values.

Therefore, by reasoning as in obtaining \equ(6.37),
one finds that for $m<m_{0}$ one has to exclude from $\EE^{(p-1)}$
a set of measure bounded by a constant times
$$ \sum_{|n| \ge N_{1}} \sum_{\ell \in \LL_{0}(n)} 
\sum_{m < m_{0}(n)} {C_{0} \over |n|^{\t+1}} \le
{\rm const.} \, \e_{0} \left(
\e_{0}^{(\t-\t_{0})/\t_{0}} + 
\e_{0}^{1+ (\t-\t_{0}-1)/\t_{0}} \right),
\Eq(6.45) $$
provided that one has $\t>\t_{0}+1>2$ the Lebesgue measure of the set
$\III^{(p)}_{2}$ is bounded by $\e_{0}^{1+\d_{2}}$, with $\d_{2}>0$.

Finally we study the measure of the set $\III^{(p)}_{3}$.
If $n>0$, $|m'|>|m|$ and $|(\tilde\o_{m}^{(p)}+\tilde\o_{m'})-
\o n| < 1$ (which again is the only case we can confine ourselves
to study), then one has to sum over $|n|\le N_{1}$,
with $N_{1}$ given by \equ(6.40). For such values of $n$ one has
$$ \eqalign{
\left| \o n - (\tilde \o_{m'}^{(p)} + \tilde \o_{m}^{(p)} ) \right|
& \ge \left| \o n - ( |m|+|m'|) \right| - {2C\e \over |m| } \cr
& \ge {2C_{0} \over |n|^{\t_{0}}} - 
{C_{0} \over |n|^{\t_{0}} } \ge {C_{0} \over |n|^{\t_{0}} }
\ge {C_{0} \over |n|^{\t} } , \cr}
\Eq(6.46) $$
as soon as $|m|>m_{0}$, with $m_{0}$ given by \equ(6.42).
Therefore we have to take into account only the values
of $m$ such that $|m|<m_{0}$, and we can also note that $|m'|$
is uniquely determined by the values of $n$ and $m$.
Then one can proceed as in the previous case and in the end one
excludes a further subset of $\EE^{(p-1)}$
whose Lebesgue measure is bounded by a constant times
$$ \sum_{|n| \ge N_{1}}
\sum_{m < m_{0}(n)} {C_{0} \over |n|^{\t+1}} \le
{\rm const.} \, \e_{0}^{1+ (\t-\t_{0})/\t_{0}} ,
\Eq(6.47) $$
so that the Lebesgue measure of the set
$\III^{(p)}_{3}$ is bounded by $\e_{0}^{1+\d_{3}}$, with $\d_{3}>0$,
provided that one takes $\t>\t_{0}$.

By summing together the bounds for $\III^{(p)}_{1}$,
for $\III^{(p)}_{2}$ and for $\III^{(p)}_{3}$, then the bound
$$ {\rm meas}(\III^{(p)}) \le b \, \e_{0}^{\d+1}
\Eq(6.48) $$
with $\d>0$, follows for all $p\ge 1$, if $\t$ is chosen to be
$\t>\t_{0}+1>2$.

We can conclude the proof of Proposition 2 through the
following result, which shows that the bound \equ(6.42)
essentially extends to the union of all $\III^{(p)}$
(at the cost of taking a larger constant $B$ instead of $b$).

\*

\0{\bf Lemma 17.} {\it Define $\III^{(p)}$ as the set of values
in $\EE^{(p)}$ verifying \equ(6.26) to \equ(6.28) for $\t>1$.
Then one has, for two suitable positive constants $B$ and $\d$,
$$ {\rm meas}\left( \cup_{p=0}^{\io} \III^{(p)} \right)
\le B \e_{0}^{\d+1} ,
\Eq(6.49) $$
where ${\rm meas}$ denotes the Lebesgue measure.}

\*

\0{\it Proof.}
First of all we check that, if we call $\e^{(p)}_{j}(n)$
the centers of the intervals $\III^{(p)}_{j}(n)$, with $j=1,2,3$,
then one has
$$ \left| \e^{(p+1)}_{j}(n) - \e^{(p)}_{j}(n) \right|
\le D \e_{0}^{p}  ,
\Eq(6.50) $$
for a suitable constant $D$.

The center $\e^{(p)}_{1}(n)$ is defined by the condition
$$ \o n \pm \left( \o_{m} - \m_{m}(\tilde\o^{(p-1)}(\e^{(p)}_{1}(n)),
\e^{(p)}_{1}( n)) \right) = 0 , 
\Eq(6.51) $$
where Whitney extensions are considered outside $\EE^{(p-1)}$;
then, by subtracting \equ(6.51) from the equivalent
expression for $p+1$, we have
$$ \left(
\m_{m}(\tilde\o^{(p)}(\e^{(p+1)}_{1}( n)),\e^{(p+1)}_{1}( n)) -
\m_{m}(\tilde\o^{(p-1)}(\e^{(p)}_{1}( n)),\e^{(p)}_{1}( n))
\right) = 0 .
\Eq(6.52) $$
In \equ(6.48) one has
$$ \eqalign{
& \m_{m}(\tilde\o'(\e'),\e') - \m_{m}(\tilde\o(\e),\e) =
\m_{m}(\tilde\o'(\e'),\e') - \m_{m}(\tilde\o(\e'),\e') \cr
& \qquad +
\m_{m}(\tilde\o(\e'),\e') - \m_{m}(\tilde\o(\e),\e')
+ (\m_{m}(\tilde\o(\e),\e') -
\m_{m}(\tilde\o(\e),\e)) , \cr}
\Eq(6.53) $$
and, from Lemma 13,
$$ \eqalign{
& \left| \m_{m}(\tilde \o'(\e'),\e') - \m_{m}(\tilde \o(\e'),\e') \right|
\le C \e \left\|\o'-\o\right\|_{\io} , \cr
& \left| \m_{m}(\tilde \o(\e'),\e') - \m_{m}(\tilde \o(\e),\e') \right|
\le C \e  \left| \e'-\e \right| , \cr}
\Eq(6.54) $$
so that we get, by lemma 10,
$$ \left| \e^{(p+1)}_{1}( n) - \e^{(p)}_{1}( n) \right|
\le C \e_{0}^{p} ,
\Eq(6.55) $$
for a suitable positive constant $C$.
This proves the bound \equ(6.50) for $j=1$.
Analogously one can consider the cases $j=2$ and $j=3$,
and a similar result is found.

By
\equ(6.10), \equ(6.51), \equ(6.54) and \equ(6.12) it follows
for $p>p_0$
$$|\e^{(p)}_{j}( n) - \e^{(p_{0})}_{j}( n)|
\le C\sum_{k=p_0}^\io \e_0^{k}=C {\e_0^{p_0}\over 1-\e_0}
\Eq(6.56) $$
%
%
so that one can ensure that $ |\e^{(p)}_{j}( n) - \e^{(p_{0})}_{j}( n)|
\le C_{0} | n|^{-\t}$ for $p\ge p_0$ by choosing
$$ p_{0} = p_{0}( n,j) \le \hbox{const.} \log | n| .
\Eq(6.57) $$
For all $p\le p_{0}$ define $\JJ^{(p)}_{j}( n)$ as the set
of values $\e$ such that \equ(6.16), \equ(6.17) and \equ(6.18)
are satisfied with $C_{0}$ replaced with $2C_{0}$.
By the definition of $p_{0}$ all the intervals $\III^{(p)}_{j}( n)$
fall inside the union of the intervals
$\JJ^{(0)}_{j}( n),\ldots,\JJ^{(p_{0})}_{j}( n)$
as soon as $p>p_{0}$. This means that, by \equ(6.31) to \equ(6.37) 
$$ \eqalign{
{\rm meas}\left( \cup_{p=0}^{\io} \III_1^{(p)} \right) &
\le {\rm meas}\left( \cup_{p=0}^{p_0} \JJ_1^{(p)} \right)
\le \sum_{|n|\ge N_{0}}
{\rm meas} (\III^{(p)}_{1}(n)) \cr
& \le \hbox{const.}
\sum_{|n|\ge N_{0}} \sum_{p=0}^{p_{0}(n,1)}
{2 C_{0} \over | n|^{\t+1} }
\le {\rm const.} \sum_{n=N_{0}}^{\io} n^{-(\t+1)} \log n \le
B \e_{0}^{1+\d} , \cr}
\Eq(6.58) $$
with suitable $B$ and $\d$, in order to take into account
the logarithmic corrections due to \equ(6.55).
Analogously one obtains the bounds
$ {\rm meas}(\III_{2}) \le B \e_{0}^{1+\d}$ and
${\rm meas}(\III_{3}) \le B \e_{0}^{1+\d}$
for the the Lebesgue measures of the sets $\III_{2}$ 
and $\III_{3}$
(possibly redefining the constant $B$).
This completes the proof of the bound \equ(6.48). \qed

\*\*
\section(7,The case 1$\minore\tt\leq$2)

\0Proposition 1 was proved assuming that $\tilde\o$ verifies
the first and the second Diophantine condition \equ(2.31) with $\t>2$.
Her we wnat to prove that it is possible to obtain a result similar
to Proposition 1 assuming only the first Diophantine
condition and $1<\t\le 2$, and that
also in such a case the set of allowed values of $\e$ 
have large Lebesgue measure, so that a result analogous
to Proposition 2 holds.

\*

The proof of the analogue of Proposition 1 is an immediate adaptation
of the analysis in Sectiones \secc(4) and \secc(5).
First we consider a slight different multiscale
decomposition of the propagator; instead of \equ(4.2)
we write
$$ {1\over - \o^{2} n_{\ell}^2 +\tilde\o_{m_{\ell}}^2}=
{\chi( \sqrt{6|\o^{2} n_{\ell}^2 -\tilde\o_{m_{\ell}}^2|})
\over - \o^{2} n_{\ell}^2 +\tilde\o_{m_{\ell}}^2} +
{\chi_{-1}( \sqrt{6|\o^{2} n_{\ell}^2 -\tilde\o_{m_{\ell}}^2|})
\over - \o^{2} n_{\ell}^2 +\tilde\o_{m_{\ell}}^2} ,
\Eq(7.1aa) $$
and the denominator in the first addend of the r.h.s. 
of \equ(7.1aa) is smaller than $C_{0}/6$; as in Section \secc(4)
we assume $C_{0}\le 1/2$. If $|n|\ne |m|$ and
$|\o^2 n^2 -\tilde\o_{m}^2| < C_{0}/6$, then,
by reasoning as in \equ(4.12), we obtain
$$|n|> |m|> {3\over 4} |n| , \qquad
\min\{ |m| , |n| \} > {1 \over 2\e} ,
\Eq(7.fon)$$
and
$$ |\o |n|- \tilde\o_m|< {C_{0} \over 6 (\o |n|+ \tilde\o_{m})} <
{C_{0}\over 6|n|} .
\Eq(7.10e)$$

We can decompose the first summand in \equ(7.1aa), obtaining
$$\chi( \sqrt{6|\o^{2} n_{\ell}^2 -\tilde\o_{m_{\ell}}^2|})
\sum_{h=-1}^{\io}
{ \chi_h( | \o n|-\tilde\o_{m})
\over - \o^{2} n^2+\tilde\o_{m}^2 }
\equiv \sum_{h=-1}^{\io} g^{(h)}(\o n,m) ,
\Eq(7.2aa)$$
and the scales from $-1$ to $h_{0}$, with $h_{0}$ given in the statement
of Lemma 5, can be bounded as in Lemma 6.

If $N_{h}(\th)$ is the number of lines on scale $h\ge h_{0}$,
and assuming only the first Melnikov condition in \equ(2.31),
the following reult holds.

\*

\0{\bf Lemma 18.} {\it
Assume that there is a constant $C_{1}$ such that
one has $| \tilde\o_{m}-|m| |\le C_{1}\e/|m|$ and $\tilde\o$
verifies the first Melnikov condition in \equ(2.31) with $\t>1$.
If $\e$ is small enough for any tree $\th \in \Th_{n,m}^{(k)}$ and
for all $h\ge h_{0}$ one has
$$ \bar N_h(\th) \le  4 K(\th) 2^{(2-h)/\t^2} -
C_{h}(\th) + M_{h}^{\n}(\th) + S_{h}(\th) , 
\Eq(7.11) $$
where $M_{h}^{\n}(\th)$ and $S_{h}(\th)$ are defined as
the number of $\n$-vertices in $\th$ such that the maximum scale
of the two external lines is $h$ and, respectively, the number of
self-energy graphs in $\th$ with $h_{T}^{e}=h$.}

\*

\0{\it Proof.}
We prove inductively the bound
$$ N^{*}_{h}(\th) \le \max\{ 0 , 2 K(\th) 2^{(1-h)/\t^2} - 1 \} ,
\Eq(7.12) $$
where $N_{h}^{*}(\th)$ is the number of non-resonant lines in $L(\th)$
on scale $h'\ge h$.
We proceed exactly as in the proof of Lemma 5.
First of all note that for a tree $\th$ to have a line on scale
$h$ the condition $K(\th) > 2^{(h-1)/\t}> 2^{(h-1)/\t^2} $ is necessary,
by the first Diophantine conditions in \equ(2.31).
This means that one can have $N_{h}^{*}(\th) \ge 1$ only
if $K=K(\th)$ is such that $K > k_{0}\=2^{(h-1)/\t}$:
therefore for values $K \le k_{0}$ the bound \equ(7.12) is satisfied.
If $K=K(\th)>k_{0}$, we assume that
the bound holds for all trees $\th'$ with $K(\th')<K$.
Define $E_{h}=2^{-1}(2^{(1-h)/\t^2})^{-1}$, we have
to prove that $N^{*}_{h}(\th)\le \max\{0,K(\th) E_{h}^{-1}-1\}$.

The dangerous case is if $m=1$ then one has a cluster $T$
with two external lines $\ell$ and $\ell_{1}$,
which are both with scales $\ge h$; then
$$ \left| |\o n_{\ell}| - \tilde\o_{m_{\ell}} \right|
\le 2^{-h+1} C_{0} , \qquad
\left| |\o n_{\ell_{1}}| - \tilde\o_{m_{\ell_{1}}} \right|
\le 2^{-h+1} C_{0} ,
\Eq(7.15a) $$
and recall that $T$ is not a self-energy graph,
so that $n_{\ell}\not=n_{\ell_1}$.
Note that the validity of both inequalities in \equ(7.15a)
for $h \ge h_{0}$ imply, by Lemma 5, that one has 
$$|n_{\ell}-n_{\ell_{1}}|\neq
|m_{\ell}\pm m_{\ell_{1}}| .
\Eq(7.10c)$$
Moreover from \equ(7.10e), \equ(7.10c) and \equ(7.fon) one obtains
$$ {C_0\over |n_{\ell}-n_{\ell_{1}}|^{\t}}\le  
|\o|n_{\ell}-n_{\ell_1}| -| m_{\ell_1}\pm m_{\ell_1}|| <
{C_{0} \over 3 \min \{ |n_{\ell}| , |n_{\ell_{1}}| \} } ,
\Eq(7.100) $$ 
which implies
$$ \left| n_{\ell}-n_{\ell_1} \right| \ge 3^{1/\t}
\min \{ |n_{\ell}| , |n_{\ell_{1}}| \}^{1/\t} , 
\Eq(7.10d)$$

Finally if $C_{0}|n|^{-\t} \le | |\o n| - \tilde\o_{m} | <
C_{0} 2^{-h+1 } $ we have $|n| \ge 2^{(h-1)/\t}$, so that
from \equ(7.10c) we obtain $|n_{\ell}-n_{\ell_{1}}| \ge
3^{1/\t} 2^{(h-1)/\t^2}$.
Then $K(\th)-K(\th_{1})> E_{h}$, which gives \equ(7.12),
by using the inductive hypothesis. \qed

\*

The analysis in Section \secc(5) can be repeated with the following
modifications. The {\it renormalized trees} are defined as in 
Section \secc(5), but the rule (9) is replaced with

\*

\0(9$^{\prime}$) To each self-energy graph $T$ the $\LL'$ operation
is applied, where  
$$\LL'\VV_T^{h}(\o n,m)= \VV_T^{h}(\o n,m) ,
\Eq(7.1aaaaa)$$
if $T$ is such that its two external lines are attached
to the same node $\vvv_{0}$ of $T$ (see as an example
the first graph of Figure 4.1), and $\LL'\VV_T^{h}(\o n,m)=0$ otherwise.

\*

With this definition we have the expansion \equ(5.4) to \equ(5.3a),
with $\n_{h,m}$ replaced with $\n_{h}$ (as $\LL'\VV_T^{h}(\o n,m)$
is independent of $n,m$).

We have that the analogue of Lemma 7 still holds also with $\LL',\RR'$
replacing $\LL,\RR$; indeed by construction
the dependence on $(n,m)$ in $\RR' \VV_T^{h_T}(\o n,m)$ 
is due to the propagators of lines $\ell$ along the path
connecting the external lines of the self-energy graph.
The propagators of such lines have the form
$$ { \chi_h(|\o n_{\ell}|-\tilde \o_{m_{\ell}})
\over - \left( |\o n_{\ell}| +\tilde\o_{m_{\ell}} \right)
 \left( |\o n_{\ell}| - \tilde\o_{m_{\ell}} \right)} ,
\Eq(7.71a) $$
and the second factor in the denominator is
bounded proportionally to $2^{-h_{\ell}}$, 
while the first is bounded proportionally
to $|m_{\ell}^{0}+m_{\ell_{T}^{e}}|$.
Then  $\RR' \VV_T^{h_T}(\o n,m)$ is bounded by a constant 
times $e^{-\k |m_{\ell_{0}}|}/|m_{\ell}^{0}+m_{\ell_{T}^{e}}|$,
hence by a constant times $1/m_{\ell_{T}^{e}}$.

We can also bound the propagator associated to one of the
external lines of $T$ with $|n_{\ell^{e}_{T}}|^{\t-1}$, by using
the first Diophantine condition, so that
then bound \equ(5.8) is replaced, using also Lemma 18, by
$$ \eqalign{
& \left| \Val(\th) \right| \le \e^{k}
\overline{D}^{k}
\Big( \prod_{h=h_0}^{\io} \exp \Big[ h \log 2 \Big( 4 K(\th)
2^{-(h-2)/\t^{2}} - C_{h}(\th) +
M_{h}^{\n}(\th) \Big) \Big] \Big) \cr
& \qquad \qquad
\Big( \prod_{T \in \SS(\th)}
{|n_{\ell^e_T}^{\t-1}|\over|m_{\ell^e_T}|}
\Big)
\Big( \prod_{h=h_{0}}^{\io} 2^{-h M_{h}^{\n}(\th)} \Big)
\cr & \qquad \qquad
\Big( \prod_{\vvv\in V(\th) \cup E(\th)} e^{-\k(|n_{v}|+|n_{v}'|)} \Big)
\Big( \prod_{\vvv\in V(\th) \cup E(\th)} e^{-\k(|m_{v}|+|m_{v}'|)} \Big)
, \cr}
\Eq(6.8a) $$
where $|n_{\ell^e_T}^{\t-1}|/|m_{\ell^{e}_{T}}|$ is less than a contant
provided one takes $\t\le 2$. Hence \equ(5.6) follows.
There is no need of Lemma 8, while the analogues of Lemma 9 and Lemma 10
can be proved with some obvious changes.
For instance in \equ(5.38) one has $\h_{\ell}=0$, which shows that
the first Mel'nikov is required, and $\n_{m}$ is replaced with $\n$,
i.e. one needs only one counterterm.

Finally also the analogue of Proposition 2 can be proved by
reasoning as in Section \secc(6), simply by observing that
in order to impose the first Mel'nikov conditions \equ(6.16)
one can take $\t_{0}=\t$; note indeed that the condition $\t>2$
was made necessary to obtain the second Mel'nikov conditions
\equ(6.17) and \equ(6.18).

\*\*
\section(8,Generalizations of the results)

\0So far we have considered only the case $\f(u)=u^{3}$ in \equ(1.1).
Now we consider the case in which the function $\f(u)$ in \equ(1.1)
is replaced with any odd analytic function
$$ \f(u)= \F u^{3} + O(u^{5}) , \qquad \F \neq 0 .
\Eq(8.1) $$
Define $\o=\sqrt{1-\l\e}$ as in Section \secc(1). Then by choosing
$\l=\s$ we obtain, instead of \equ(1.13),
$$ \eqalign{
\hbox{Q} \qquad & \cases{
n^{2} a_{n} = \left[ (v+w)^{3} + O(\e) \right]_{n,n} , & \cr
- n^{2} a_{n} = \left[ (v+w)^{3} + O(\e) \right]_{n,-n} , & \cr} \cr
\hbox{P} \qquad & \left( - \o^2 n^2 + m^2 \right) w_{n,m} = \e
\left[ (v+w)^{3} + O(\e) \right]_{n,m} , \qquad |m| \neq |n| , \cr}
\Eq(8.2) $$
where $O(\e)$ denotes analytic functions in $u$ and $\e$
of order at least one in $\e$. Then we can introduce an auxiliary
parameter $\m$, by replacing $\e$ with $\e\m$ in \equ(8.2)
(recall that at the end one has to set $\m=1$).
Then we can proceed as in the previous sections.
The equation for $a_{0}$ is the same as before, and only the
diagrammatic rules for the coefficients $u^{(k)}_{n,m}$
have to be slightly modified.
The nodes can have any odd number of entering lines,
that is $s_{\vvvv}=1,3,5,\ldots$.
For nodes $\vvv$ of $w$-type with $s_{\vvvv} \ge 3$ the node factor
is given by $\h_{\vvvv}=\e^{(s_{\vvvv}-1)/2}$, which has to be added
to the list of node factors \equ(3.3), while for nodes
$\vvv$ of $v$-type one has to add to the list in \equ(3.2)
other (obvious) contributions, arising from the fact that
the function $f^{(k)}_{n}$ in \equ(2.12)
has to be replaced with
$$ f^{(k)}_{n} = \sum_{k'=1}^{k} f^{(k,k')}_{n} ,
\Eq(8.3) $$
where the function $f^{(k,1)}_{n}$ is given by \equ(2.12),
while for $k\ge 2$ and $1< k' \le k$ one has
$$ f^{(k,k')}_{n} =
- \sum_{\scriptstyle k_{1}+\ldots+k_{2k'+1}=k-k'}
\sum_{\scriptstyle n_{1}+\ldots+n_{2k'+1} = n \atop
\scriptstyle m_{1}+\ldots+m_{2k'+1} = m}
u^{(k_{1})} _{n_{1},m_{1}} \ldots u^{(k_{2k'+1})} _{n_{2k'+1},m_{2k'+1}} ,
\Eq(8.4) $$
where the symbols have to be interpreted according to \equ(2.13)
and \equ(2.14); note that $s=2k'+1$ is the number of lines
entering the node in the corresponding graphical representation.

In the same way one has to modify \equ(2.29) by replacing
the last term in the r.h.s. with
$$ [\f(v+w)]^{(k)}_{n,m} = \sum_{k'=1}^{k}
[\f(v+w)]^{(k,k')}_{n,m} ,
\Eq(8.5) $$
where $[\f(v+w)]^{(k,1)}_{n,m}$ is given by the old \equ(2.30),
while all the other terms are given by expressions analogous
to \equ(8.4).

The discussion then proceeds exactly as in the previous cases.
Of course we have to use that, by writing
$$ \f(u) = \sum_{k=1}^{\io} \F_{k} u^{2k+1} ,
\Eq(8.6) $$
with $\F_{1}=\F$, the constants $\F_{k}$ can be bounded
for all $k \ge 1$ by a constant to the power $k$
(which follows from the analyticity assumption).

By taking into account the new diagrammatic rules,
one change in the proper way the definition of the tree value,
and a result analogous to Lemma 2 is easily obtained.
The second statement of Lemma 3 has to be changed into
$|E(\th)|\le \sum_{\vvvv\in V(\th)} (s_{\vvvv}-1)+1$
(see \equ(3.24)), while the bound on $|V_{v}^{3}(\th)|$
still holds. For $V_{v}^{s}(\th)$, with $s\ge 5$,
one has to use the factors $\e^{(s-1)/2}$.

Nothing changes in the following sections, except that
the bound \equ(5.11) has to be suitably modified
in order to take into account the presence of the new kinds of
nodes (that is the nodes with branching number more than three).

At the end we obtain the proof of Theorem 1 for general
odd nonlinearities starting from the third order.
Until now we are still confining ourselves to the case $j=1$.

If we choose $j>1$ we have perform a preliminary rescaling
$u(t,x) \to j u (jt,jx) $, and we write down the equation
for $U(t,x)=u(jt,jx)$. If $\f(u)=Fu^{3}$ we see immediately
that the function $U$ solves the same equation as before,
so that the same conditions on $\e$ has to be imposed
in order to find a solution. In the general case
\equ(8.2) holds with $j^{2}\e$ replacing $\e$.
This completes the proof of Theorem 1 in all cases.

\*\*
\appendix(A1,The solution of {\equ(1.5)})

\0The odd $2\pi$-periodic solutions of \equ(1.5) can be found in the
following way \cita{LS}, \cita{PBC}.
First we consider \equ(1.5), with $\langle a_{0}^{2} \rangle$
replaced by a parameter $c_{0}$, and we see that
$$ a_{0}(\x) = V \, \sn (\O\x,\mm) ,
\Eqa(A1.1) $$
where $\sn(\O\x,\mm)$ is the sine-amplitude function
with modulus $\sqrt{m}$ \cita{GrR}, \cita{AS}, is an odd solution of
$$ \ddot a_{0} = - 3 c_{0} a_{0} - a_{0}^{3} ,
\Eqa(A1.2) $$
if the following relations are verified by $V,\O,c_0,\mm$
$$ V = \sqrt{ -2\mm} \O , \qquad \qquad
{V^{2} \over 6 c_{0} + V^{2}} = - \mm ,
\Eqa(A1.3) $$
Of course both $V$ and $\mm$ can be written as function of $c_0$
and $\O$ as
$$ \mm = {{3 c_0\over\O^2}}-1 , \qquad \qquad
V = \sqrt{2-{6 c_0\over\O^2}} \, \O .
\Eqa(A1.4) $$
In particular one finds
$$ \dpr_\O V = \sqrt{-2 \mm} + {1\over \sqrt{-2\mm}} \, {6 c_0\over\O^2}
= {\sqrt{-2 \mm}\over \mm} \left( \mm -{3 c_0\over\O^2} \right) =
{\sqrt{2\over -\mm}} ,
\Eqa(A1.5)$$
so that one has
$$ {\O{\dpr_{\O}}V \over V} = {1\over -\mm } .
\Eqa(A1.6) $$

If we impose also that the solution \equ(A1.1) is $2\p$-periodic (and we
recall that $4 K(\mm)$ is the natural period of the sine-amplitude
$\sn(\x,\mm)$, \cita{AS}), we obtain $\O = \O_{\mmmm} = 2K(\mm)/\p$,
and $V$ is fixed to the value $V_{\mmmm}=\sqrt{-2\mm}\O_{\mmmm}$.

Finally imposing that $c_{0}$ equals the average of $a_{0}^{2}$
fixes $\mm$ to be the solution of \cita{PBC}
$$ E(\mm) = K(\mm) \, {7+\mm \over 6} ,
\Eqa(A1.7) $$
where
$$ K(\mm) = \int_{0}^{\p/2} {\rm d}\th \, {1 \over
\sqrt{1 - \mm \sin^{2}\th} } , \qquad
E(\mm) = \int_{0}^{\p/2} {\rm d}\th \, 
\sqrt{ 1 - \mm \sin^{2} \th}
\Eqa(A1.8) $$
are, respectively, the complete elliptic integral of first kind and
the complete elliptic integral of second kind \cita{AS}, and we have
used that the average of $\sn^{2}(\O\x,\mm)$ is $(K(\mm)-E(\mm))/K(\mm)$.
This gives $\mm\approx -0.2554$.

One can find $2\p/j$-periodic solutions by noticing that
equations \equ(1.5) are invariant by the simmetry $a_{0}(\x)
\rightarrow \a \, a_{0}(\a \x)$, so that a complete set of odd
$2\p$-periodic solutions is provided by $a_0(\x,j)= j \, a_0(j\x)$.

\*\*
\appendix(A2, Proof of Lemma 1)

\0The solution of \equ(2.19) is found by variation of constants.
We write $\langle a_{0}^{2} \rangle = c_{0}$, and consider $c_{0}$
as a parameter. Let $a_{0}(\x)$ be the solution \equ(A1.1)
of equation \equ(A1.2), and define the Wronskian matrix
$$ W(\x) = \left( \matrix{
w_{11}(\x) & w_{12}(\x) \cr
w_{21}(\x) & w_{22}(\x) \cr} \right) ,
\Eqa(A2.1) $$
which solves the linearized equation 
$$ \dot W = M(\x) W , \qquad
M(\x) = \left( \matrix{
0 & 1 \cr -3a_{0}^{2}(\x)-3 c_{0} & 0 \cr} \right) ,
\Eqa(A2.2) $$
and is such that $W(0)=\openone$ and $\det W(\x)=1$ $\forall\x$.
We need two independent solutions of the the linearized
equation. We can take one as $\dot a_0$, the other as
$\partial_\O a_0$ (we recall that $\mm$ and $V$ are function of
$\O$ and $c_0$ through \equ(A1.3)); then one has
$$ \eqalignno{
w_{11}(\x) & = {1\over \O V} \dot a_{0}(\x) =
\cn(\O\x,\mm)\,\dn(\O\x,\mm) , \cr
w_{21}(\x) & = \dot w_{11}(\x) =
-\sn(\O\x,\mm) \left( \dn^{2}(\O\x,m) +
\mm \, \cn^{2}(\O\x,\mm) \right) , \cr
w_{12}(\x) & = {1 \over V (1+D_{\mmmm})} \dpr_{\O} a_{0}(\x)
& \eqa(A2.3) \cr
& = B_{\mmmm} \left( \x \, \cn(\O\x,\mm) \, \dn(\O\x,\mm) +
\O^{-1} D_{\mmmm} \sn(\O\x,\mm) \right) , \cr w_{22}(\x)
& = \dot w_{12}(\x) \cr
& = \cn(\O\x,\mm) \, \dn(\O\x,\mm) - \O B_{\mmmm} \x \, \sn(\O\x,\mm)
\left( \dn^{2}(\O\x,\mm) + \mm \, \cn^{2}(\O\x,\mm) \right) , \cr} $$
where we have used \equ(A1.3) to define the dimensionless constants
$$ D_{\mmmm} = {\O \dpr_{\O} V \over V } = {1\over -\mm} , \qquad\qquad
B_{\mmmm} = {1 \over \left( 1 + D_{\mmmm} \right) } = {- \mm \over 1 - \mm} .
\Eqa(A2.4) $$

As we are interested in the case $c_{0}=\langle a_{0}^{2} \rangle$
we set $\O=\O_{\mmmm}=2K(\mm)/\p$ (in order to have the period equal to
$2\p$) and we fix $\mm$ as in Appendix \secc(A1).

Then, by defining $X=(y,\dot y)$ and $F=(0,h)$,
we can write the solution of \equ(2.19) as the first component of
$$ X(\x) = W(\x)\,\lis X +
W(\x) \int_{0}^{\x} {\rm d}\x' \, W^{-1}(\x') \, F(\x') ,
\Eqa(A2.5) $$
where $\lis X = (0,\dot y(0))$ denote the corrections
to the initial conditions $(a_{0}(0),\dot a_{0}(0))$,
and $y(0)=0$ as  we are looking for an odd solution. 

Shorten $c(\x)\=\cn(\O_{\mmmm}\x,\mm)$,
$s(\x)\=\sn(\O_{\mmmm}\x,\mm)$, and $d(\x)\=\dn(\O_{\mmmm}\x,\mm)$,
and define $cd(\x)=\cn(\O_{\mmmm}\x,\mm)\,\dn(\O_{\mmmm}\x,\mm)$.
One can write the first component of \equ(A2.5) as
$$ \eqalign{
y(\x) & = w_{12}(\x) \, \dot  y(0) + \int_{0}^{\x} {\rm d} \x'
\left( w_{12}(\x)\,w_{11}(\x') -
w_{11}(\x)\,w_{12}(\x') \right) h(\x') \cr
& =   B_{\mmmm} \left( \x \, cd(\x) +
\O_{\mmmm}^{-1}D_{\mmmm} s(\x) \right)\dot  y(0) \cr
& \qquad + B_{\mmmm} \left( cd(\x) \int_{0}^{\x} {\rm d}\x'
\int_{0}^{\x'} {\rm d}\x'' \, cd(\x'')\, h (\x'') \right. \cr
& \qquad \left. + \O_{\mmmm}^{-1} D_{\mmmm}
\left( s(\x) \int_{0}^{\x} {\rm d}\x' \, cd(\x')\,h(\x') -
cd(\x) \int_{0}^{\x} {\rm d}\x'\, s(\x') \, h(\x')
\right) \right) , \cr}
\Eqa(A2.6) $$
as we have explictly written
$$ \eqalignno{
& w_{12}(\x)\,w_{11}(\x')-w_{11}(\x)\,w_{12}(\x')
& \eqa(A2.7) \cr
& \quad = B_{\mmmm} \left( \x\, cd(\x)\,cd(\x')
+ \O_{\mmmm}^{-1}D_{\mmmm} s(\x)\,cd(\x') - cd(\x)\,\x'\, cd(\x') -
\O_{\mmmm}^{-1}D_{\mmmm} cd(\x)\,s(\x') \right) , \cr} $$
and integrated by parts
$$ \eqalign{
& \int_{0}^{\x} {\rm d}\x' \left(
\x \, cd(\x) \, cd(\x') - cd(\x)\,\x'\,cd(\x') \right) h(\x') \cr
& \qquad \qquad \qquad \qquad =
cd(\x) \int_{0}^{\x} {\rm d}\x' \int_{0}^{\x'}
{\rm d}\x'' cd(\x'')\,h(\x'') . \cr}
\Eqa(A2.8) $$
By using that if $\PP[F]=0$ then $\PP[\II[F]]=0$ 
and that $\II$ switches parity we can rewrite \equ(A2.6) as
$$ \eqalign{
y(\x) & =  B_m \Big( \x \, cd(\x)  \Big(\dot y(0)-
\II[ cd\, h](0) - \O_{\mmmm}^{-1}D_{\mmmm}\la s\, h\ra\Big) + 
\O_{\mmmm}^{-1}D_{\mmmm} s(\x) \Big( \dot y(0)-\II[cd\,h](0)\Big) \cr
& \qquad + \O_{\mmmm}^{-1}D_m\big (s(\x) \, \II[cd\, h](\x) - cd(\x) \,
\II[\PP[s \, h]](\x) \big) + cd(\x) \, \II[\II[cd\,h]](\x) \Big) . \cr}
\Eqa(A2.9) $$
This is an odd $2\p$-periodic analytic function provided that we choose
$$ \dot y(0) - \II[cd \, h](0) -
\O_{\mmmm}^{-1}D_{\mmmm}\la s \, h \ra = 0 ,
\Eqa(A2.10) $$ 
which fixes the parameter $\dot y(0)$; hence \equ(2.20) is found

\*\*
\appendix(A3, Proof of {\equ(2.25)})

\0We have to compute $\langle a_0 \LLL[a_0] \rangle$, which is
given by (recall that one has $a_0(\x)=V_{\mmmm}\,s(\x)$
and $cd(\x)=\O_{\mmmm}^{-1} \dot s(\x)$)
$$ \la a_0 \LLL[a_0]\ra = B_{\mmmm} \O_{\mmmm}^{-2} V_{\mmmm}^{2}
\left( D_{\mmmm}^{2} \la s^{2} \ra^{2}  +
D_{\mmmm} \la s^{2} \, \II [s \dot s ] \ra -
D_{\mmmm} \la s \dot s \,\II [\PP[ s^{2}] ] \ra + 
\la s\dot s \II[\II [s \dot s] ] \ra \right) .
\Eqa(A3.1)$$
Using that one has $\II[s \dot s] = (s^2- \langle s^2 \rangle)/2$,
we obtain
$$ \la s^2 \II [s \dot s ] \ra = {1\over 2}
\left( \la s^4\ra-\la s^2\ra^2 \right) ;
\Eqa(A3.2) $$
moreover, integrating by parts, we find
$$ \eqalign{
& \la s \dot s  \,\II[\PP[ s^2]] \ra = -
{1\over 2} \la s^4\ra+{1\over 2}\la s^2\ra^2, \cr
& \la s\dot s \II[\II[s \dot s]] \ra = -
{1\over 4} \la s^4 \ra+{1\over 4}\la s^2\ra^2 , \cr}
\Eqa(A3.3) $$
so that we finally get 
$$ \la a_0 \LLL[a_0] \ra=  {1 \over 2} V_{\mmmm}^{2}
\O_{\mmmm}^{-2} B_{\mmmm}
 \left( \left( 2D_{\mmmm} - {1\over2} \right)
\langle s^{4} \rangle  + \left( 2D_{\mmmm} \left( D_{\mmmm}-1 \right)
+ {1\over 2} \right) \langle s^{2} \rangle^{2} \right) ,
\Eqa(A3.4) $$
which is strictly positive by \equ(A2.4) and by the choice of $\mm$
according to Appendix \secc(A1); then \equ(2.25) follows.

\*\*
\appendix(A4, Proof of Lemma 15)

\0We shall prove inductively on $p$ the bounds \equ(6.12).
From \equ(6.4) we have
$$ | \tilde\o^{(p)}_{m}(\e) - \tilde\o^{(p-1)}_{m}(\e) |
\le C | \n_{m} (\tilde \o^{(p-1)}(\e),\e) -
\n_{m} (\tilde \o^{(p-2)}(\e),\e) | ,
\Eqa(A4.1) $$
as we can bound $|\tilde\o^{(p)}_{m}(\e)+\tilde\o^{(p-1)}_{m}(\e)|
\ge 1$ for $\e \in \EE^{(p)}$.

We set, for $|m| \ge 1$
$$ \D \n_{h,m} \= \n_{h,m} (\tilde \o^{(p-1)}(\e),\e) -
\n_{h,m} (\tilde \o^{(p-2)}(\e),\e)
= \lim_{q \to \io} \D \n^{(q)}_{h,m} ,
\Eqa(A4.2) $$
where we have used the notations \equ(5.25) to define
$$ \D \n^{(q)}_{h,m} = \n^{(q)}_{h,m} (\tilde\o^{(p-1)}(\e),\e) -
\n^{(q)}_{h,m} (\tilde \o^{(p-2)}(\e),\e) .
\Eqa(A4.3) $$
We want to prove inductively on $q$ the bound
$$ \left| \D \n^{(q)}_{h,m} \right| \le C \e 
\| \tilde \o^{(p-1)}(\e) - \tilde \o^{(p-2)}(\e) \|_{\io} ,
\Eqa(A4.4) $$
for some constant $C$, uniformly in $q$, $h$ and $m$.

For $q=0$ the bound \equ(A4.4) is trivially satisfied.
Then assume that \equ(A4.4) hold for all $q'<q$.

For simplicity we set $\tilde \o=\tilde \o^{(p-1)}(\e)$ and
$\tilde \o'=\tilde \o^{(p-2)}(\e)$.
We can write, from \equ(6.25), for $|m|\ge 1$
and for $h\ge 0$,
$$ \eqalign{
& \D \n^{(q)}_{h,m} = - \sum_{k = h}^{\bar h - 1} 2^{-k-2} 
\Big( \b^{(q)}_{k,m}(\tilde\o,\e,\{\n^{(q-1)}_{k'}(\tilde\o,\e)\})
\cr & \qquad \qquad -
\b^{(q)}_{k,m}(\tilde \o',\e,\{\n^{(q-1)}_{k'}(\tilde\o',\e)\})
\Big) , \cr} 
\Eqa(A4.5) $$
where we recall that $\b^{(q)}_{k,m}
(\tilde\o,\e,\{\n^{(q-1)}_{k'}(\tilde\o,\e)\})$
depend only on $\n^{(q-1)}_{k'}(\tilde\o,\e)$ with $k'\le k-1$, and
we can set
$$ \eqalign{
& \b^{(q)}_{k,m}(\tilde\o,\e,\{\n^{(q-1)}_{k'}(\tilde\o,\e)\}) \cr
& \qquad \qquad
= \b^{(a)(q)}_{k,m}(\tilde\o,\e,\{\n^{(q-1)}_{k'}(\tilde\o,\e)\}) -
\b^{(b)(q)}_{k,m}(\tilde\o,\e,\{\n^{(q-1)}_{k'}(\tilde\o,\e)\}) , \cr}
\Eqa(A4.6) $$
according to the settings in \equ(2.3).
Then we can split the differences in \equ(A4.6) into
$$ \eqalign{
& \b^{(q)}_{k,m}(\tilde\o,\e,\{\n^{(q-1)}_{k'}(\tilde\o,\e)\}) -
\b^{(q)}_{k,m}(\tilde\o',\e,\{\n^{(q-1)}_{k'}(\tilde\o',\e)\}) \cr
& \qquad =
\left( \b^{(q)}_{k,m}(\tilde\o,\e,\{\n^{(q-1)}_{k'}(\tilde\o,\e)\}) -
\b^{(q)}_{k,m}(\tilde\o',\e,\{\n^{(q-1)}_{k'}(\tilde\o,\e)\}) \right) \cr
& \qquad \qquad +
\left( \b^{(q)}_{k,m}(\tilde\o',\e,\{\n^{(q-1)}_{k'}(\tilde\o,\e)\}) -
\b^{(q)}_{k,m}(\tilde\o',\e,\{\n^{(q-1)}_{k'}(\tilde\o',\e)\})
\right) , \cr}
\Eqa(A4.7) $$
and we bound separately the two terms.

The second term can be expressed as sum of trees
$\th$ which differ from the previously considered ones as,
among the nodes $v$ of $w$-type with only one entering line,
there are some with $\n^{(c_{\vvvv})(q-1)}_{k_{\vvvv}}(\tilde\o,\e)$,
some with $\n^{(c_{\vvvv})(q-1)}_{k_{\vvvv}}(\tilde\o',\e)$ and one with
$\n^{(c_{\vvvv})(q-1)}_{k_{\vvvv}}(\tilde\o,\e)-
\n^{(c_{\vvvv})(q-1)}_{k_{\vvvv}}(\tilde\o',\e)$. Then we can bound
$$ \eqalign{
& \left| \b^{(q)}_{k,m}(\tilde\o',\{\n^{(q-1)}_{k'}(\tilde\o,\e)\}) -
\b^{(q)}_{k,m}(\tilde\o',\{\n^{(q-1)}_{k'}(\tilde\o',\e)\}) \right| \cr
& \qquad \qquad 
\le D_{1} \e\, \sup_{h'\ge 0} \sup_{|m'| \ge 2}
\D \n^{(q-1)}_{h',m'}
\le D_{1} C \e^{2} \left\| \tilde\o' - \tilde\o \right\|_{\io} , \cr}
\Eqa(A4.8) $$
by the inductive hypothesis.

We are left with the first term in \equ(A4.8).
We can reason as in \cita{GM} (which we refer to for details),
and at the end, instead of (A9.14) of \cita{GM}, we obtain
$$ \eqalignno{
& \left| \Val'(\th_{0})
\left( g^{(h_{\ell_{j}})}_{\ell_{j}}(\tilde\o') -
g^{(h_{\ell_{j}})}_{\ell_{j}}(\tilde\o) \right)
\Val'(\th_{1})
\Big( \prod_{i=2}^{s}  \Val(\th_{i}) \Big) \right|
& \eqa(A4.9) \cr
& \qquad \qquad
\le C^{|V(T)|} e^{-\k K(T)/2}
\e^{|V_{w}(T)|+|V_{v}^{1}(T)|}
\left\| \tilde \o' - \tilde\o \right\|_{\io} \cr
& \qquad \qquad
\le C^{|V(T)|} \e^{|V_{w}(T)|+|V_{v}^{1}(T)|}
e^{-\k K(T)/4} e^{-\k 2^{(k-1)/\t}/4}
\left\| \tilde \o' - \tilde\o \right\|_{\io} , \cr} $$
where $K(T)$ is defined in \equ(4.7).

Therefore can bound $\left\| \tilde\o^{(p)}(\e)-
\tilde\o^{(p-1)}(\e) \right\|_{\io}$
with a constant times $\e$ times the same expression with
$p$ replaced with $p-1$, \ie
$\left\| \tilde\o^{(p-1)}(\e)-\tilde\o^{(p-2)}(\e)
\right\|_{\io}$, so that, by the inductive hypothesis,
the bound \equ(6.12) follows.




\rife{AS}{1}{
M. Abramowitz, I.A. Stegun:
{\it Handbook of mathematical functions},
Do\-ver, New York, 1965. }
\rife{Ba}{2}{
D. Bambusi,
{\it Lyapunov center theorem for some nonlinear PDE's: a simple proof},
Ann. Scuola Norm. Sup. Pisa Cl. Sci. (4)
{\bf 29} (2000), no. 4, 823--837. }
\rife{BP}{3}{
D. Bambusi, S. Paleari,
{\it Families of periodic solutions of resonant PDEs},
J. Nonlinear Sci.
{\bf 11}  (2001),  no. 1, 69--87. }
\rife{BB1}{4}{
M. Berti, Ph. Bolle,
{\it Periodic solutions of nonlinear wave
equations with general nonlinearities},
Comm. Math. Phys.
{\bf 243} (2003), no. 2, 315--328. }
\rife{BB2}{5}{
M. Berti, Ph. Bolle,
{\it Multiplicity of periodic solutions of nonlinear 
wave equations},
Nonlinear Anal.,
to appear.}
\rife{BL}{6}{
G.D. Birkhoff, D.C. Lewis,
{\it On the periodic motions near a given periodic motion
of a dynamical system},
Ann. Math.
{\bf 12} (1933), 117--133. }
\rife{Bo1}{7}{
J. Bourgain,
{\it Construction of quasi-periodic solutions for Hamiltonian
perturbations of linear equations and applications to nonlinear PDE},
Internat. Math. Res. Notices
{\bf 1994}, no. 11, 475ff., approx. 21 pp. (electronic). }
\rife{Bo2}{8}{
J. Bourgain,
{\it Construction of periodic solutions of nonlinear
wave equations in higher dimension},
Geom. Funct. Anal.
{\bf 5} (1995), 629--639. }
\rife{Bo3}{9}{
J. Bourgain,
{\it Quasi-periodic solutions of Hamiltonian perturbations of 2D
linear Schr\"o\-din\-ger equations},
Ann. of Math. (2)
{\bf 148} (1998), no. 2, 363--439. }
\rife{BCN}{10}{
H. Br\'ezis, J.-M. Coron, L. Nirenberg,
{\it Free vibrations for a nonlinear wave equation and a theorem
of P. Rabinowitz},
Comm. Pure Appl. Math.
{\bf 33} (1980), no. 5, 667--684. } 
\rife{BN}{11}{
H. Br\'ezis, L. Nirenberg,
{\it Forced vibrations for a nonlinear wave equation},
Comm. Pure Appl. Math.
{\bf 31} (1978), no. 1, 1--30. }
\rife{Cr}{12}{
W. Craig,
{\it Probl\`emes de petits diviseurs dans les \'equations aux
d\'eriv\'ees partielles},
Pano\-ra\-mas and Syntheses 9. Soci\'et\'e Math\'ematique de France,
Paris, 2000. }
\rife{CW1}{13}{
W. Craig and C.E. Wayne,
{\it Newton's method and periodic solutions of nonlinear
wave equations},
Comm. Pure Appl. Math.
{\bf 46} (1993), 1409--1498. }
\rife{CW2}{14}{
W. Craig and C.E. Wayne,
{\it Nonlinear waves and the $1:1:2$ resonance},
``Singular limits of dispersive waves'' (Lyon, 1991), 297--313,
NATO Adv. Sci. Inst. Ser. B Phys., 320, Plenum, New York, 1994. } %
\rife{FR}{15}{
E.R. Fadell, P.H. Rabinowitz,
{\it Generalized cohomological index theories for Lie group
actions with an application to bifurcation questions
for Hamiltonian systems},
Invent. Math.
{\bf 45}  (1978), no. 2, 139--174. }
\rife{GG}{16}{
G. Gallavotti, G. Gentile,
{\it Hyperbolic low-dimensional invariant tori and
summations of divergent series},
Comm. Math. Phys.
{\bf 227} (2002), no. 3, 421--460. }
\rife{Ge}{17}{
G. Gentile,
{\it Whiskered tori with prefixed frequencies and Lyapunov spectrum},
Dynam. Stability Systems
{\bf 10} (1995), no. 3, 269--308. }
\rife{GGM}{18}{
G. Gallavotti, G. Gentile, V. Mastropietro,
{\it A field theory approach to Lindstedt series for hyperbolic
tori in three time scales problems},
J. Math. Phys.
{\bf 40} (1999), no. 12, 6430--6472. } 
\rife{GM}{19}{
G. Gentile and V. Mastropietro,
{\it Construction of periodic solutions of the nonlinear wave equation
with Dirichlet boundary conditions by the Linsdtedt series method},
J. Math. Pures Appl. (9),
to appear. }
\rife{GoR}{20}{
C. Godsil, G. Royle,
{\it Algebraic graph theory},
Graduate Texts in Mathematics 207.
Spring\-er, New York, 2001. }
\rife{GrR}{21}{
I.S. Gradshteyn, I.M. Ryzhik,
{\it Table of integrals, series, and products},
Sixth edition, Academic Press, Inc., San Diego, 2000. }
\rife{HP}{22}{
F. Harary, E.M. Palmer,
{\it Graphical enumeration},
Academic Press, New York-London, 1973. }
\rife{K1}{23}{
S.B. Kuksin,
{\it Nearly integrable infinite-dimensional Hamiltonian systems},
Lecture Notes in Mathematics 1556, Springer, Berlin, 1994. }
\rife{K2}{24}{
S.B. Kuksin,
{\it Fifteen years of KAM for PDE},
Preprint, 2003. }
\rife{KP}{25}{
S.B. Kuksin, J.P\"oschel,
{Invariant Cantor manifolds of quasi-periodic oscillations
for a nonlinear Schr\"odinger equation},
Ann. of Math. (2)
{\bf 143} (1996), no. 1, 149--179. }
\rife{LS}{26}{
B.V. Lidski\u\i, E.I. Shul$^{\prime}$man,
{\it Periodic solutions of the equation $u_{tt}-u_{xx}+u^{3}=0$},
Funct. Anal. Appl.
{\bf 22} (1988), no. 4, 332--333. }
\rife{L}{27}{
A. M. Lyapunov,
{\it Probl\`eme g\'en\'eral de la  stabilit\'e du mouvement},
Ann. Sc. Fac. Toulouse
{\bf 2} (1907), 203--474. }
\rife{PBC}{28}{
S. Paleari, D. Bambusi, S. Cacciatori,
{\it Normal form and exponential stability for some
nonlinear string equations},
Z. Angew. Math. Phys.
{\bf 52}  (2001),  no. 6, 1033--1052. }
\rife{P}{29}{
J. P\"oschel,
{\it Quasi-periodic solutions for a nonlinear wave equation},
Comment. Math. Helv.
{\bf 71} (1996), no. 2, 269--296. }
\rife{R1}{30}{
H.P. Rabinowitz,
{\it Periodic solutions of nonlinear hyperbolic
partial differential equations},
Comm. Pure Appl. Math.
{\bf 20} (1967), 145--205. }
\rife{R2}{31}{
P.H. Rabinowitz,
{\it Free vibrations for a semilinear wave equation},
Comm. Pure Appl. Math.
{\bf 31} (1978), no. 1, 31--68. }
\rife{Wa}{32}{
C.E. Wayne,
{\it Periodic and quasi-periodic solutions of nonlinear wave
equations via KAM theory},
Comm. Math. Phys.
{\bf 127} (1990), no. 3, 479--528. }
\rife{We}{33}{
A. Weinstein,
{\it Normal modes for nonlinear Hamiltonian systems}, Invent. Math.
{\bf 20} (1973), 47--57. }
\rife{Wh}{34}{
H. Whitney,
{\it Analytic extensions of differential
functions defined in closed sets},
Trans. Amer. Math. Soc.
{\bf 36} (1934), no. 1, 63--89. }
\biblio
\ciao